\documentclass[review]{elsarticle}

\pdfoutput=1
\usepackage{graphicx}
\usepackage{amsfonts, amsmath, amssymb}
\usepackage{booktabs,siunitx}
\usepackage[svgnames,table]{xcolor}
\usepackage[tableposition=above]{caption}
\usepackage{pifont}
\usepackage{bm}
\usepackage{cancel}
\usepackage{caption}
\usepackage{color}
\usepackage{enumerate}
\usepackage{float}
\usepackage{hyperref}
\usepackage[english]{babel}
\usepackage[margin=1in]{geometry}
\usepackage{mathtools}
\usepackage{setspace}
\usepackage{subfigure}
\usepackage{lineno}
\newcommand \D [2]{\frac{\partial #1}{\partial #2}}

\renewcommand{\vec}[1]{\bm{\mathrm{#1}}}
\newcommand{\V}[1]{\bm{\mathrm{#1}}}

\def \div{\nabla \cdot \mbox{}}
\def \grad{\nabla}

\def \thetaFD{\theta^{\text{FD}}}

\def \x{\vec{x}}
\def \y{\vec{y}}

\def \n{\vec{n}}

\def \u{\vec{u}}

\def \I{\vec{I}}
\def \F{\vec{F}}

\def \U{\vec{U}}
\def \L{\vec{L}}

\def \cM{\vec{\mathcal{M}}}

\def \Sb{S_b}

\def \Vbt{V_b(t)}
\def \Vb{V_b}

\def \C{\vec{C}}

\def \cF{\vec{\mathcal{F}}}

\def \F{\vec{F}}

\def \g{\vec{g}}

\def \I{\vec{I}}
\def \Ib{\I_{\text{b}}}

\def \Mb{\text{M}_{\text{b}}}

\def \Nx{N_x}
\def \Ny{N_y}

\def \Omegal{\Omega_{\text{l}}}
\def \Omegag{\Omega_{\text{g}}}

\def \R{\vec{R}}
\def \U{\vec{U}}
\def \Ub{\U_{\text{b}}}

\def \Ur{\U_{\text{r}}}
\def \W{\vec{W}}

\def \Wr{\W_{\text{r}}}

\def \X{\vec{X}}
\def \Xcom{\X_{\text{com}}}

\def \cS{\vec{\mathcal{S}}}
\def \cJ{\vec{\mathcal{J}}}

\def \cO{\mathcal{O}}

\def \f{\vec{f}}

\def \fc{\f_{\text{c}}}

\def \half{\frac{1}{2}}
\def \3half{\frac{3}{2}}
\def \5half{\frac{5}{2}}

\def \mul{\mu_{\text{l}}}
\def \mus{\mu_{\text{s}}}
\def \mug{\mu_{\text{g}}}
\def \n{\vec{n}}

\def \ncells{n_{\text{cells}}}
\def \ncycles{n_{\text{cycles}}}

\def \rhol{\rho_{\text{l}}}
\def \rhos{\rho_{\text{s}}}
\def \rhog{\rho_{\text{g}}}

\def \s{\vec{s}}

\def \sgn{\textrm{sgn}}
\def \u{\vec{u}}

\def \ub{\u_{\text{b}}}

\def \x{\vec{x}}

\def \div{\nabla \cdot \mbox{}}
\def \grad{\nabla}

\def \dt{\Delta t}
\def \dx{\Delta x}
\def \dy{\Delta y}

\def \delU{\Delta \vec{U}}

\def \Ds{{\mathrm d}\s}

\def \Dx{{\mathrm d}\x}

\def \ds{\Delta s}
\def \dt{\Delta t}
\def \dx{\Delta x}

\newcommand{\upperRomannumeral}[1]{\uppercase\expandafter{\romannumeral#1}}

\newcommand{\REVIEW}[1]{{#1}}

\begin{document}

\begin{frontmatter}
	
\title{Simulating water-entry/exit problems using Eulerian-Lagrangian and fully-Eulerian fictitious domain methods within the open-source IBAMR library}
\author[SDSU]{Amneet Pal Singh Bhalla\corref{mycorrespondingauthor}}
\ead{asbhalla@sdsu.edu}
\author[Northwestern]{Nishant Nangia}
\author[Polito]{Panagiotis Dafnakis}
\author[Polito]{Giovanni Bracco}
\author[Polito]{Giuliana Mattiazzo}

\address[SDSU]{Department of Mechanical Engineering, San Diego State University, San Diego, CA, USA}
\address[Northwestern]{Department of Engineering Sciences and Applied Mathematics, Northwestern University, Evanston, IL, USA}
\address[Polito]{Department of Mechanical and Aerospace Engineering, Politecnico di Torino, Turin, Italy}
\cortext[mycorrespondingauthor]{Corresponding author}

\begin{abstract}
In this paper we employ two implementations of the fictitious domain (FD) method to simulate water-entry and 
water-exit problems and demonstrate their ability to simulate practical marine engineering problems. 
In FD methods, the fluid momentum equation is extended within the solid domain using an additional body force
that constrains the structure velocity to be that of a rigid body. Using this formulation, a single set of
equations is solved over the entire computational domain. The constraint force is calculated in two distinct ways: 
one using an Eulerian-Lagrangian framework of the immersed boundary (IB) method and another using a fully-Eulerian 
approach of the Brinkman penalization (BP) method. Both FSI strategies use the same multiphase flow algorithm that solves 
the discrete incompressible Navier-Stokes system in conservative form. A consistent transport scheme is employed to advect 
mass and momentum in the domain, which ensures numerical stability of high density ratio multiphase flows involved in 
practical marine engineering applications.  Example cases of a free falling wedge \REVIEW{(straight and inclined)} and cylinder are 
simulated, and the numerical results are compared against benchmark cases in literature.
%to validate both FD implementations
\end{abstract}


\end{frontmatter}

%%%%%%%%%%%%%%%%%%%%%%%%%%%%%
\section{Introduction}
Fluid-structure interaction (FSI) at the free water surface is a fundamental hydrodynamics problem that is 
of great importance to engineers working in the fields of naval architecture and marine engineering~\cite{Tveitnes08,Wu2004,Howison91,Zhu07}. 
These unsteady and nonlinear FSI problems can be further divided into two main categories of 
water-entry and water-exit of structures~\cite{Greenhow87,Greenhow88}. Some practical examples
of water-entry and exit problems include hydrodynamic impact on bow structures of ships during slamming and wave run-up 
effects on marine platforms~\cite{Hou18}. The water-entry of free falling marine structures produce large impact loads 
that can threaten their immediate and long-term safety. Therefore, it is important to estimate the impact loads on  
structures to ensure their safe design and operability. Several experimental~\cite{Wu2004,Yettou2006}, 
theoretical~\cite{von1929,Dobrovol69,Watanabe86}, and more recently computational fluid dynamics 
(CFD) techniques~\cite{Kleefsman2005,Zhang2010,Hou18,Nair2018} have been used in the literature to study these problems at a fundamental level. 
The latter approach is the subject matter of the current study. 

With the advancement of computing technology, it is now possible to simulate full three dimensional, unsteady FSI problems
involving complex structural geometries. Presently, both boundary element method (BEM) based~\cite{Takagi2004,Lee1995wamit,Lee2006wamit-user} 
and incompressible Navier-Stokes (INS) method based simulations~\cite{Shen2015,Bhushan2011} are routinely employed 
in the design process of marine structures. The BEM is based on potential flow equations, which ignore the full nonlinear convective and viscous dissipation terms found in the INS equations; hence, BEM based solvers are much faster in compute time than INS based solvers.
However, INS solvers are more general than BEM solvers and therefore 
can reliably model complicated phenomena like wave breaking, wave overtopping, and wave run-up over structures~\cite{Jacobsen2012,Kasem2010,Chen14}. Moreover, several commercial CFD codes such as  \REVIEW{Fluent}~\cite{Fluent2009} 
or \REVIEW{STAR-CCM+}~\cite{Star2012}, and open-source codes such as \REVIEW{OpenFOAM}~\cite{Higuera2013} now include support for modeling structures interacting with free water surface. In these codes, the most robust way of simulating large-displacement FSI is to use two sets of meshes: an underlying fluid mesh and an overlaying structure mesh for providing boundary conditions to the fluid 
solver. The fluid mesh can be block-structured, whereas the structure mesh is generally unstructured to represent complex geometries. 
The inside region of the solid is not meshed as it does not participate in the solution process~\footnote{Only a no-slip condition on the fluid-structure interface is required to simulate FSI of a rigid body.}. This dual mesh 
approach to FSI is also known as the Chimera or overset mesh method~\cite{Shen2015,Carrica2007} in the literature. These codes and the 
FSI strategy based on overset meshes have been adopted by several marine and ocean research groups to model complicated FSI problems~\cite{Shen2015,Carrica2007,Hou18,Facci2016}.     
 
In this work we discuss a different FSI strategy to simulate water-entry and water-exit problems and demonstrate that it can be efficiently 
and reliably used to simulate \REVIEW{these} complex  problems. Our FSI strategy is based on the \emph{fictitious domain} (FD) 
methodology~\cite{Patankar2000,Sharma2005}, in which the fluid equations are extended \emph{into} the solid domain. 
The velocity field inside the solid is obtained by solving the common momentum and incompressibility equations. 
The main advantage of this FSI strategy is that complex structures can be modeled without employing an unstructured mesh for 
the solid. Additionally, the momentum and continuity equations can be solved on structured Cartesian grids using fast linear solvers. Another benefit of this strategy is 
the implicit-coupling of the fluid and structure; stability-preserving sub-iterations between the fluid and solid domains are not required, in contrast with overset mesh based methods. \REVIEW{For non neutrally-buoyant structures that are frequently encountered in marine engineering applications, numerically stable solutions are obtained in a single iteration using the strongly-coupled FD formulations. This stability preserving characteristics of FD methods are attributed to resolution of a correct or physical density field in the inertial term of the momentum equation (including the inner region occupied by the immersed solid); this circumvents any numerical issues pertaining to added/reduced mass effects.}

There are various approaches to implement FD methods. Two of the most popular approaches are the \emph{immersed boundary} (IB) 
method~\cite{Peskin02} and the \emph{Brinkman penalization} (BP) method~\cite{Angot99,Carbou03}. In the FD/IB approach the structure is 
tracked in Lagrangian frame of reference, whereas in FD/BP method the immersed structure is tracked on the Eulerian grid 
itself (usually by some indicator function). Both methods solve the INS equations in the entire domain with an additional body forcing 
term in the solid domain. The additional force in the momentum equation acts like a constraint force, which imposes rigid body velocity in the region occupied by the structure. FD/IB methods estimate this constraint force in the Lagrangian form and transfers it back to the Eulerian grid using suitable IB kernels~\cite{Bhalla13,Vanella09}. FD/BP methods calculate the constraint force in Eulerian form directly~\cite{Bergmann11}. 
\REVIEW{An efficient time-splitting approach using the \emph{distributed Lagrange multiplier} method (DLM) of Sharma and Patankar~\cite{Sharma2005} 
is employed to calculate the Lagrangian constraint force.} In contrast, FD/BP methods apply the constraint force in a time-implicit manner 
while solving the momentum equation. Sec.~\ref{sec_sol_method} describes the full time-stepping algorithm for these two methods in detail. \REVIEW{We remark that there are various versions of the IB method described in the literature~\cite{Mittal05}; in this work the \emph{original} IB method machinery of Peskin~\cite{Peskin02} is employed.}

\REVIEW{More recently, fully-implicit FD methods have been proposed in the literature. In these methods, the \emph{unknown} constraint forces are solved for as a part of an extended system of INS equations; i.e. the fluid velocity, pressure, and rigidity-enforcing constraint forces are solved for simultaneously as a large block matrix system~\cite{Stein2016,Stein2017,Kallemov16,Usabiaga17,Feldman2016}. These developments have been enabled due to advances in linear algebra techniques such as physics-based preconditioning for iterative Krylov solvers~\cite{Kallemov16,Usabiaga17}, and large scale sparse direct solvers for systems involving dense Schur complements~\cite{Feldman2016}. Some of these methods also allow for higher-order spatiotemporal discretization schemes and general boundary condition treatments on the fluid-structure interface~\cite{Stein2016,Stein2017}.} 

Another class of hybrid overset mesh/IB methods, known as the sharp-interface approach~\cite{Borazjani13,RMittal08,YHTseng03,Calderer2014}, is often used to simulate complicated FSI problems.
These methods solve the fluid equations on regular Cartesian grids but 
zero-out solution inside the solid domain programmatically. The computational domain is divided into fluid nodes, solid nodes, 
and IB nodes. IB nodes are located near the solid surface in the fluid side and provide velocity boundary conditions 
to the fluid nodes. Sub-iterations are generally required to couple the two domains in order to maintain numerical stability. These methods also require complicated computational geometry algorithms to impose structure velocities at the IB nodes.
Sharp-interface methods have also been used to model marine engineering FSI problems with success~\cite{Calderer2014,Zhang2010}.
However, we do not consider them in this work as they are fundamentally different to FD methods.

%We discretize the equations of motion in conservative form and employ a consistent transport of mass and momentum 
%for all phases~\cite{Nangia2019}. The consistent transport of mass and momentum ensures numerical stability for high 
%density ratio multiphase flows~\cite{Rudman1998,Raessi2012,Nangia2018}. We use the level set method~\cite{Osher2006,Sussman1994} 
%to track air-water and fluid-solid interfaces in this work. Sec.~\ref{sec_cons_ins} describes the consistent conservative scheme in the context of the incompressible Navier-Stokes solver.  

The remainder of the paper is organized as follows. We first introduce the continuous and discrete system 
of equations in Secs.~\ref{sec_cont_eqs} and~\ref{sec_sol_method}, respectively. Next we describe the time-stepping schemes 
of the FD/BP and FD/IB methods in Sec.~\ref{sec_temporal_scheme}. Comparisons and salient features 
of the two FD implementations are described in Sec.~\ref{sec_comp_fd}. Software implementation is described in Sec.~\ref{sec_sfw}. 
\REVIEW{A 2D free-falling inclined wedge with three free degrees of freedom is simulated in Sec.~\ref{sec_inclined_wedge} to validate the 
implementation of the FD/BP method, which is relatively new compared to our more mature FD/IB implementations~\cite{Bhalla13,Nangia17,Nangia2019,Kallemov16,Usabiaga17}.} 
Simulations of water-entry and water-exit of a free falling wedge and a cylinder with both FD methods are shown in the remainder of
Sec.~\ref{sec_examples} and the results are compared with literature. \REVIEW{Finally, computational costs of the two FD implementations  
are compared in Sec.~\ref{sec_cost_comparison}.}

%%%%%%%%%%%%%%%%%%%%%%%%%%%%%%
\section{The continuous equations of motion} \label{sec_cont_eqs}
\subsection{Multiphase fictitious domain formulation}
\label{sec_cont_eqns}
We begin by describing the continuous governing equations for a coupled
multiphase fluid-structure system occupying a fixed region of space $\Omega \subset \mathbb{R}^d$, for $d = 2$ or $3$
spatial dimensions. In fictitious domain FSI formulations, 
the momentum and divergence-free condition for the domain occupied by the fluid and structure are described
in a fixed Eulerian coordinate system $\x = (x_1, \ldots, x_d) \in \Omega$.
For the FD/IB method, a Lagrangian description of the immersed body configuration is employed, in which 
$\s = (s_1, \ldots s_d) \in B$ denotes the fixed material coordinate system 
attached to the structure and $B \subset \mathbb{R}^d$ is the Lagrangian curvilinear 
coordinate domain. The position of the immersed structure is denoted by $\X (\s,t)$ in the Lagrangian frame,
with the body occupying the volumetric Eulerian region $\Vbt \subset \Omega$ at time $t$.
In contrast, the FD/BP method 
uses an indicator function $\chi(\x,t)$ defined on the Eulerian grid to describe the location of the body. 
The indicator function is non-zero only in the structure domain $\Vbt$. We use spatially and temporally varying density 
$\rho(\x,t)$ and dynamic viscosity $\mu(\x,t)$ fields to model not only multiple fluids occupying the domain, 
\REVIEW{but also non neutrally-buoyant structures}.
The equations of motion for the coupled fluid-structure system for the fictitious domain formulation read as
\begin{align}
\D{\rho \u(\x,t)}{t} + \div \rho\u(\x,t)\u(\x,t) &= -\grad p(\x,t) + \div \left[\mu \left(\grad \u(\x,t) + \grad \u(\x,t)^T\right) \right]+ \rho\g + \fc(\x,t) , \label{eqn_momentum}\\
  \div \u(\x,t) &= 0. \label{eqn_continuity} 
\end{align}
Eqs.~\eqref{eqn_momentum} and \eqref{eqn_continuity} are the incompressible 
Navier-Stokes momentum and continuity equations written in conservative form for the fixed region in space $\Omega$.
Here, $\u(\x,t)$ is the fluid velocity, 
$p(\x,t)$ is the pressure, and $\fc(\x,t)$ is the Eulerian constraint force density that is non-zero 
only in the region occupied by the structure.
The gravitational acceleration is denoted by $\g = (g_1, \ldots, g_d)$.
\REVIEW{In the present study, we choose to directly work with the conservative form of the momentum equation, since
it has been shown that methods based on the non-conservative form exhibit numerical instabilities
for problems involving air-water interfaces~\cite{Raessi2008, Raessi2012, Desjardins2010, Ghods2013, Nangia2018}.} 

The specific form of the constraint force $\fc(\x,t)$ depends on the particular fluid-structure interaction algorithm employed. For the FD/IB methodology, the constraint forces are first calculated on the Lagrangian mesh and later 
transferred to the background Eulerian grid. Conversely, the Lagrangian structure is displaced by interpolating 
the background Eulerian velocity onto the Lagrangian domain. The interactions between 
Eulerian and Lagrangian quantitates are mediated by integral transformations using a Dirac delta function, usually 
defined as a tensor product of one-dimensional singular kernels $\delta(\x) = \Pi_{i=1}^{d}\delta_i(x_i)$. The Lagrangian-Eulerian 
interaction equations are written as
\begin{align}  
\fc(\x,t)  &= \int_{B} \F(\s,t) \, \delta(\x - \X(\s,t)) \, \Ds, \label{eqn_F_f} \\
  \U(\s,t) &= \int_{\Omega} \u(\x,t) \, \delta(\x - \X(\s,t)) \, \Dx, \label{eqn_u_interpolation} \\
   \D{\X}{t} (\s,t) &= \U(\s,t). \label{eqn_body_motion} 
\end{align}
Eq.~\eqref{eqn_F_f} relates the Lagrangian force density $\F(\s,t)$ to a corresponding Eulerian density $\fc(\x,t)$, 
which is commonly known as \emph{force spreading} operation in the IB literature~\cite{Peskin02}. Eq.~\eqref{eqn_u_interpolation} 
relates the physical velocity of each Lagrangian material point $\U(\s,t)$ to the background Eulerian velocity field 
 $\u(\x,t)$, hence defining the \emph{velocity interpolation} operation. The velocity interpolation operation ensures that the immersed 
structure moves according to the local fluid velocity $\u(\x,t)$ (see Eq.~\eqref{eqn_body_motion}), and that the
no-slip condition is implicitly satisfied at the fluid-solid interfaces. \REVIEW{The standard discretization of these operators are described later in Sec.~\ref{sec_lagIB}, and we refer readers to Peskin~\cite{Peskin02} for a detailed analysis of their properties.}
\REVIEW{In the FD/IB formulation the appearance of the Eulerian constraint force density $\fc(\x,t)$ is due to
a rigidity constraint imposed on the Lagrangian velocity field $\frac{1}{2}\left[\grad \U(\s,t) + \grad \U(\s,t)^T\right] = 0$,
which is continuously enforced through a \emph{distributed Lagrange multiplier} force field (see Patankar et al.~\cite{Patankar2000} and Sharma and Patankar~\cite{Sharma2005}). Discretely, an approximation to the Lagrangian force density $\F(\s,t)$ is computed and spread onto the background Eulerian grid; this process is described briefly in Sec.~\ref{sec_ts_ib}, and in more detail by Shirgaonkar et al.~\cite{Shirgaonkar2009} and Bhalla et al~\cite{Bhalla13}.}

In the FD/BP method, the constraint force $\fc(\x,t)$ is defined as a (Brinkman) penalization force that enforces a desired 
motion $\ub(\x,t)$ in the spatial location occupied by the body.
More specifically,  the immersed structure is treated as a porous
body with a vanishing permeability $K \ll 1$ \REVIEW{(effectively making the region impenetrable and translate with the desired rigid body velocity) and the penalization force is formulated as~\cite{Angot99,Bergmann11,Gazzola2011}}
\begin{align}  
\fc(\x,t)  &= \frac{\chi(\x,t)}{K}\left(\ub(\x,t) - \u(\x,t)\right).
\end{align}
In contrast with the immersed boundary method, the FD/BP method is a purely Eulerian approach to modeling the fluid-structure system.
Sec.~\ref{sec_ts_bp} describes the numerical algorithm for computing the rigid body velocity $\ub(\x,t)$ from the fluid-structure interaction.

\subsection{Interface tracking}
\label{sec_cont_ls}
To prescribe the material properties (i.e. density and viscosity) for the three phases on the Eulerian grid, we use 
two scalar level set functions $\phi(\x,t)$ and $\psi(\x,t)$. 
The zero-contour of $\phi$ function implicitly defines the liquid-gas interface, whereas the zero-contour of $\psi$ function defines 
the structure boundary. The level set function $\phi$ conveniently allows prescription of the
liquid density $\rhol$ and viscosity $\mul$ in the spatial region $\Omegal(t) \subset \Omega$ occupied by the liquid, and 
the gas density $\rhog$ and viscosity $\mug$ in the spatial region $\Omegag(t) \subset \Omega$ occupied by the gas. 
The codimension-1 interface between these two fluids is denoted as $\Gamma(t) = \Omegal \cap \Omegag$.
The complex topological changes in the liquid-gas interface due to fluid-fluid and fluid-structure interactions can be easily handled within
the level set framework without employing any remeshing procedures~\cite{Osher1988,Sethian2003,Sussman1994}. 
Level set methods are also relatively simple to implement on locally-refined meshes. Similar to the $\phi$ level set, the $\psi$ level set function allows prescription of the solid density $\rhos$ and viscosity $\mus$ in the region occupied by the structure $\Vbt$. 
The codimension-1 boundary of the immersed structure is denoted as $\Sb(t) = \partial \Vb(t)$.

As the simulation progresses in time, all three phases are advected by the incompressible Eulerian velocity field. 
This phase transport is governed by the conservative, linear level set advection equations
\begin{align}
\D{\phi}{t} + \div \phi \u &= 0, \label{eq_ls_fluid_advection} \\
\D{\psi}{t} + \div \psi \u &= 0. \label{eq_ls_solid_advection}
\end{align}
The density and viscosity in the three phases are determined as a function of these auxiliary fields as 
\begin{align}
\rho (\x,t) &= \rho(\phi(\x,t), \psi(\x,t)), \label{eq_rho_ls}\\
\mu (\x,t) &= \mu(\phi(\x,t), \psi(\x,t)) \label{eq_mu_ls}.
\end{align}
In practice, regularized Heaviside functions (see  Sec.~\ref{sec_reinit}) are used to obtain the discretized form of 
Eqs.~\eqref{eq_rho_ls} and~\eqref{eq_mu_ls}.

Although various functional forms can be used to define the level sets in multiphase flow applications,
the most practical form is the \emph{signed distance function}. At time $t = 0$, the distances to $\Gamma(t)$ and
$\Sb(0)$ are computed and set as initial conditions to the level set advection 
equations~\eqref{eq_ls_fluid_advection} and~\eqref{eq_ls_solid_advection}:

\begin{align}
\phi\left(\x, 0\right) &= 
\begin{cases} 
       \min\limits_{\y \in \Gamma(0)} \|\x-\y\|,  & \x \in \Omegag(0), \\
        -\min\limits_{\y \in \Gamma(0)} \|\x-\y\|,  & \x \in \Omegal(0), \label{eq_signed_phi}
\end{cases} \\
\psi\left(\x, 0\right) &= 
\begin{cases} 
       \min\limits_{\y \in \Sb(0)} \|\x-\y\|,  & \x \not \in \Vb(0), \\
        -\min\limits_{\y \in \Sb(0)} \|\x-\y\|,  & \x \in \Vb(0). \label{eq_signed_psi}
\end{cases}
\end{align}
Notice that under advection, there is no guarantee that $\phi$ and $\psi$ will remain signed distance
 functions~\cite{Chopp1993}.
At each time step, a reinitialization procedure is used to maintain the signed distance properties.
Note that by using this formulation, initial conditions are only required for $\phi$ and $\psi$ and not for
$\rho$ and $\mu$.

\section{Spatial discretization} \label{sec_discrete_eqs}
 \label{sec_discrete_eqs}
This section describes the discretization of the governing equations for the coupled
FSI system for both fictitious domain formulations. For the FD/BP method,
we use only Eulerian quantities that are discretized on a staggered Cartesian grid, whereas for the FD/IB method,
additional Lagrangian quantities are approximated on a collection of immersed markers.
The Lagrangian markers can be arbitrarily positioned on the background Eulerian grid without conforming to 
the grid lines. Regularized versions of the Dirac delta function are used to define discrete grid transfer operations for the FD/IB method.
To simplify the treatment of the two methods, we focus on describing the $d = 2$ spatial dimensions case; the discretization in three spatial dimensions is analogous.
We refer readers to prior studies~\cite{Bhalla13,Nangia2019} for a description of the FD/IB method in $3$D.

\subsection{Eulerian discretizaton for FD methods}
\label{eulerian_discretization}
We employ a staggered Cartesian grid discretization for quantities described in the Eulerian frame; 
see Fig.~\ref{fig_discrete_diagram}. A Cartesian grid \emph{made up of $\Nx \times \Ny$ cells} covers the physical, 
rectangular domain $\Omega$ with mesh spacing $\dx$ and $\dy$ in each direction.
Assuming that the bottom left corner of the domain
is situated at the origin $(0,0)$, each cell center of the grid has
position $\x_{i,j} = \left((i + \half)\dx,(j + \half)\dy\right)$
for $i = 0, \ldots, \Nx - 1$ and $j = 0, \ldots, \Ny - 1$.
For a given cell $(i,j)$, $\x_{i-\half,j} = \left(i\dx,(j + \half)\dy\right)$ is the physical location of the cell face 
that is half a grid space away from $\x_{i,j}$ in the $x$-direction,
and $\x_{i,j-\half,k} =\left((i + \half)\dx,j\dy\right)$ is the physical location of the cell 
face that is half a grid cell away from  $\x_{i,j}$ in the $y$-direction.
The flow and structure level sets, and pressure degrees of freedom are approximated at cell centers and are
denoted by $\phi_{i,j}^{n} \approx \phi\left(\x_{i,j}, t^n\right)$,
$\psi_{i,j}^{n} \approx \psi\left(\x_{i,j}, t^n\right)$, and $p_{i,j}^{n} \approx p\left(\x_{i,j},t^{n}\right)$, respectively.
Here, $t^n$ denotes the time at time step $n$.
The material properties are also approximated at cell centers, $\rho_{i,j}^{n} \approx \rho\left(\x_{i,j},t^{n}\right)$
and  $\mu_{i,j}^{n} \approx \mu\left(\x_{i,j},t^{n}\right)$;
these quantities are interpolated onto the required degrees of freedom as needed (see ~\cite{Nangia2018} for further details).
Velocity components are staggered and are defined on their respective cell faces:
$u_{i-\half,j}^{n} \approx u\left(\x_{i-\half,j}, t^{n}\right)$, and
$v_{i,j-\half}^{n} \approx v\left(\x_{i,j-\half}, t^{n}\right)$.
The components of the \REVIEW{gravitational and constraint forces} on the right-hand side of the momentum equation
are also approximated on respective faces of the staggered grid.

Classic second-order finite differences are used to discretize spatial derivative
operators and are denoted with $h$ subscripts; i.e. $\grad \approx \grad_h$. We refer readers 
to prior studies~\cite{Nangia2018,Cai2014,Griffith2009,Bhalla13} for a full description
of these staggered grid discretizations.

\subsection{Lagrangian discretization for the FD/IB method}
\begin{figure}[]
  \centering
  \subfigure[\REVIEW{Continuous domain}]{
    \includegraphics[scale = 0.3]{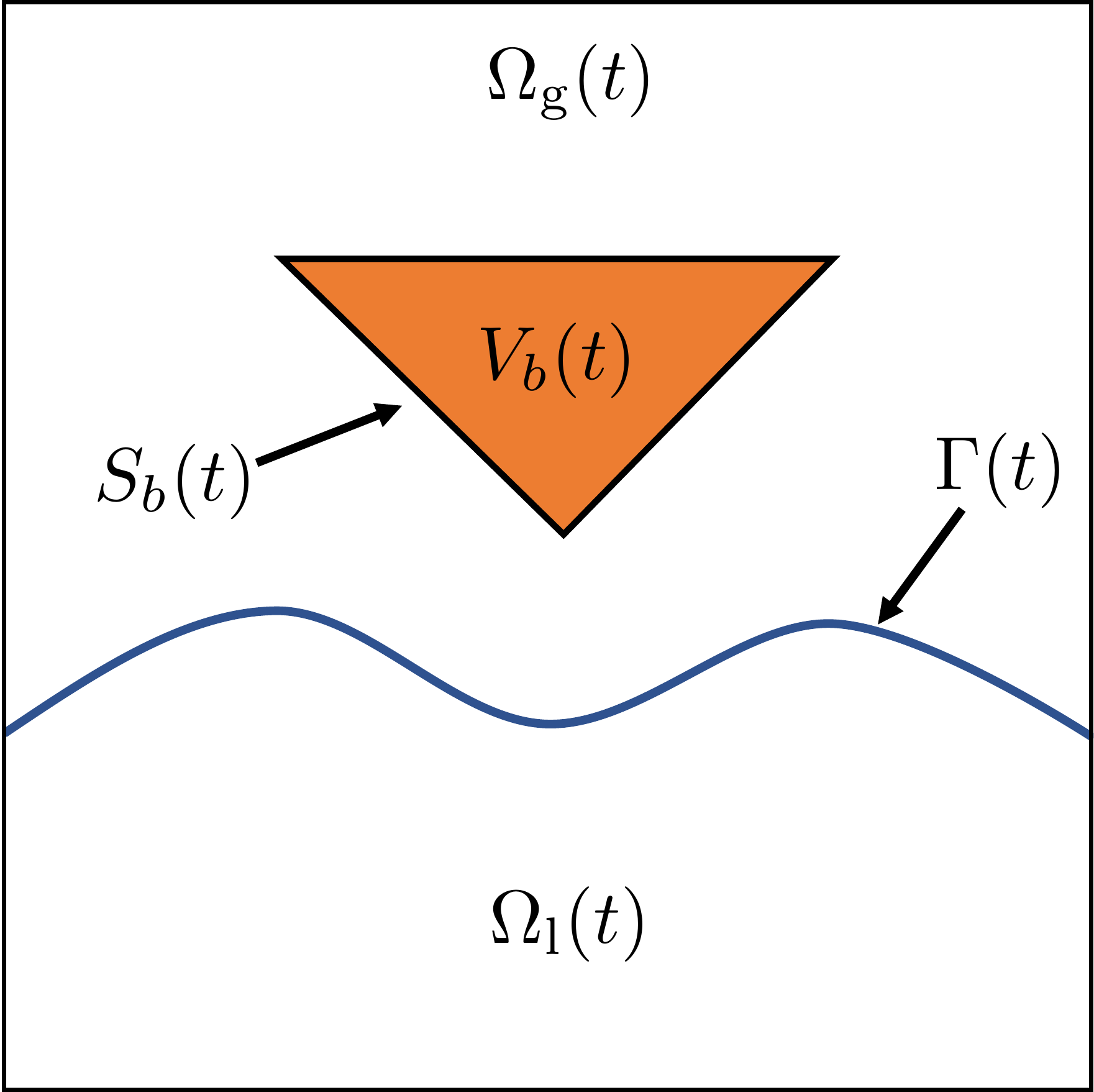}
    \label{ib_continuous}
  }
   \subfigure[\REVIEW{FD/BP discretized domain}]{
    \includegraphics[scale = 0.3]{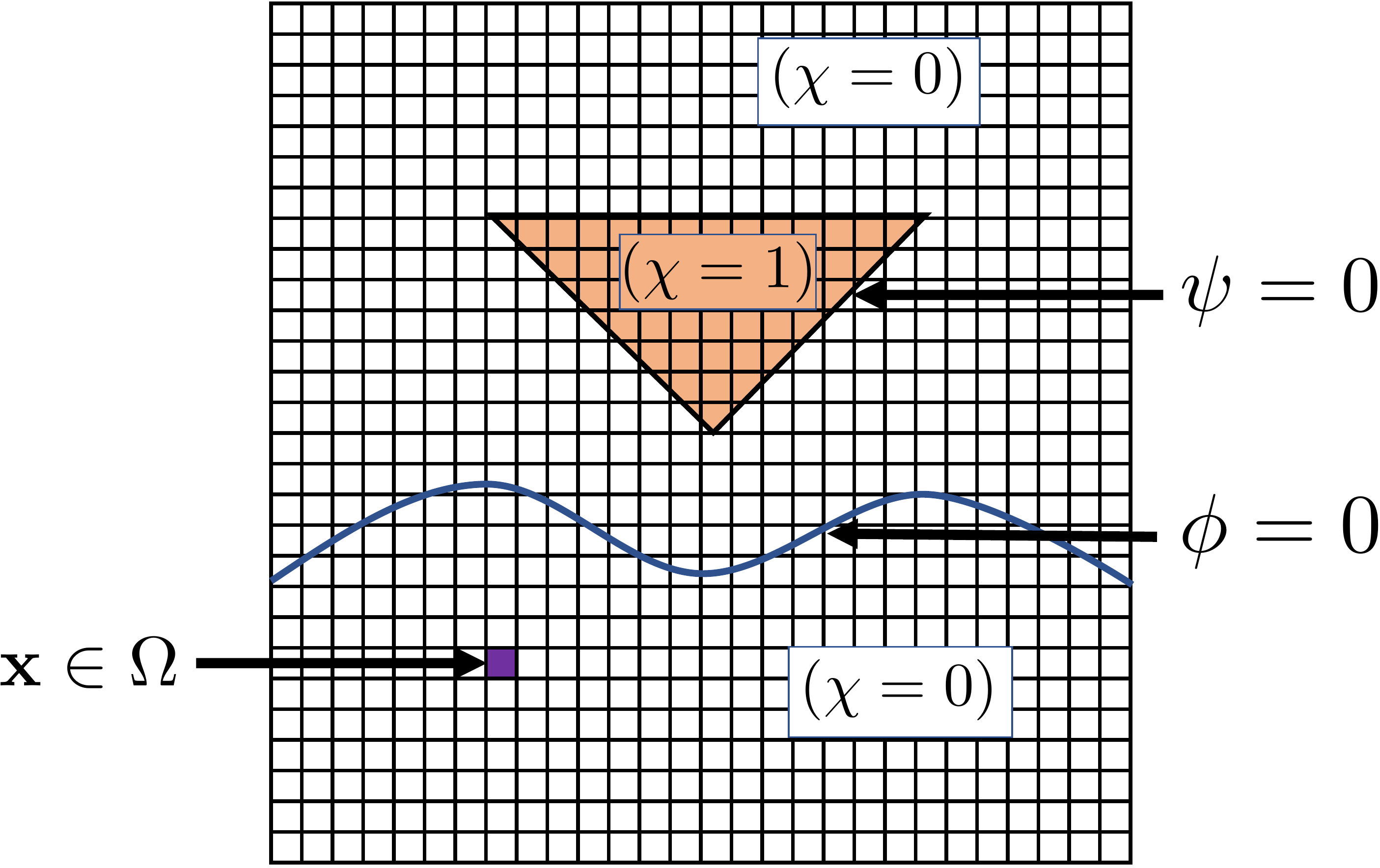}
    \label{bp_discrete}
  }
     \subfigure[\REVIEW{FD/IB discretized domain}]{
    \includegraphics[scale = 0.3]{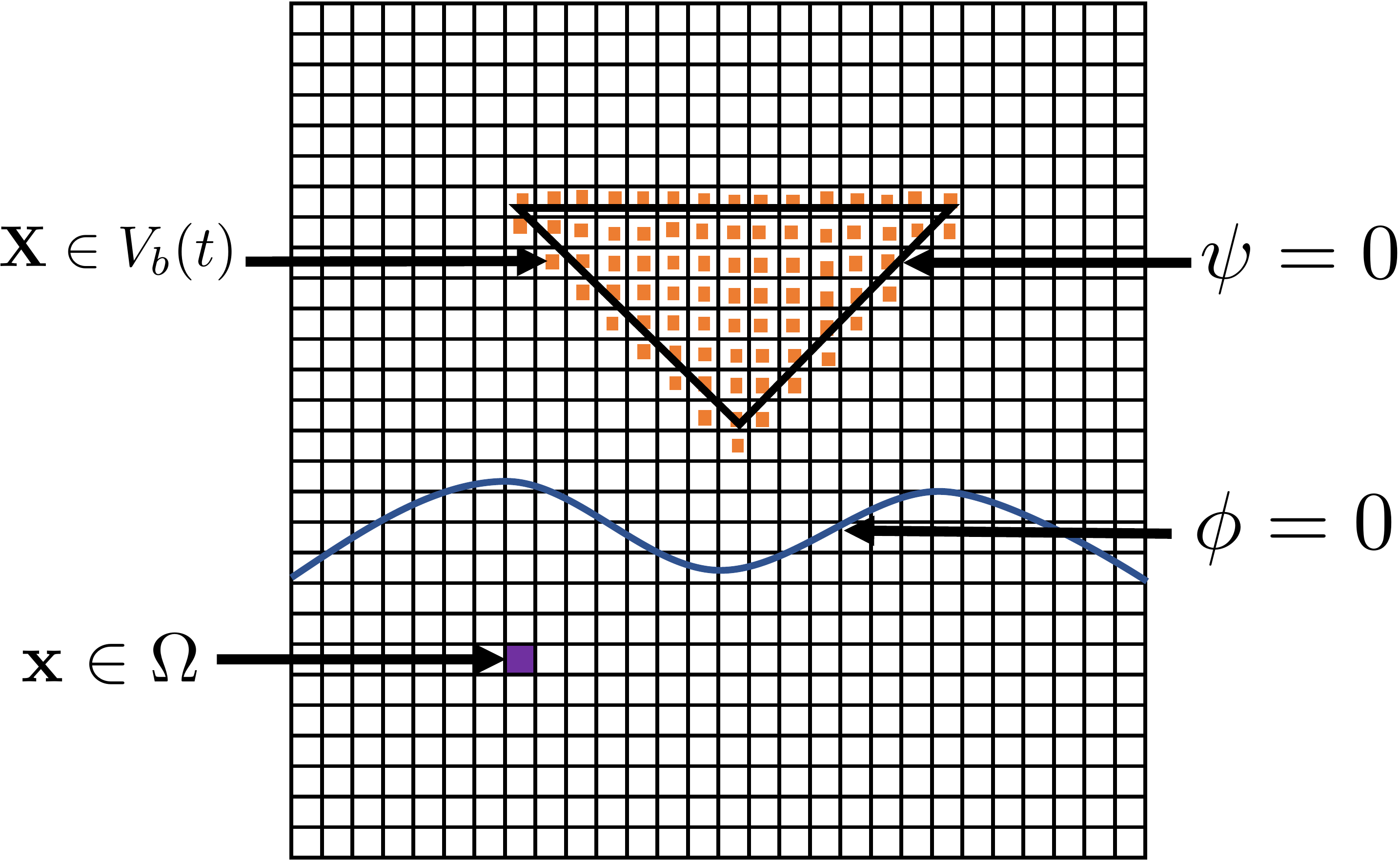}
    \label{ib_discrete}
  }
  \subfigure[\REVIEW{Two grid cells \& Lagrangian markers}]{
    \includegraphics[scale = 0.3]{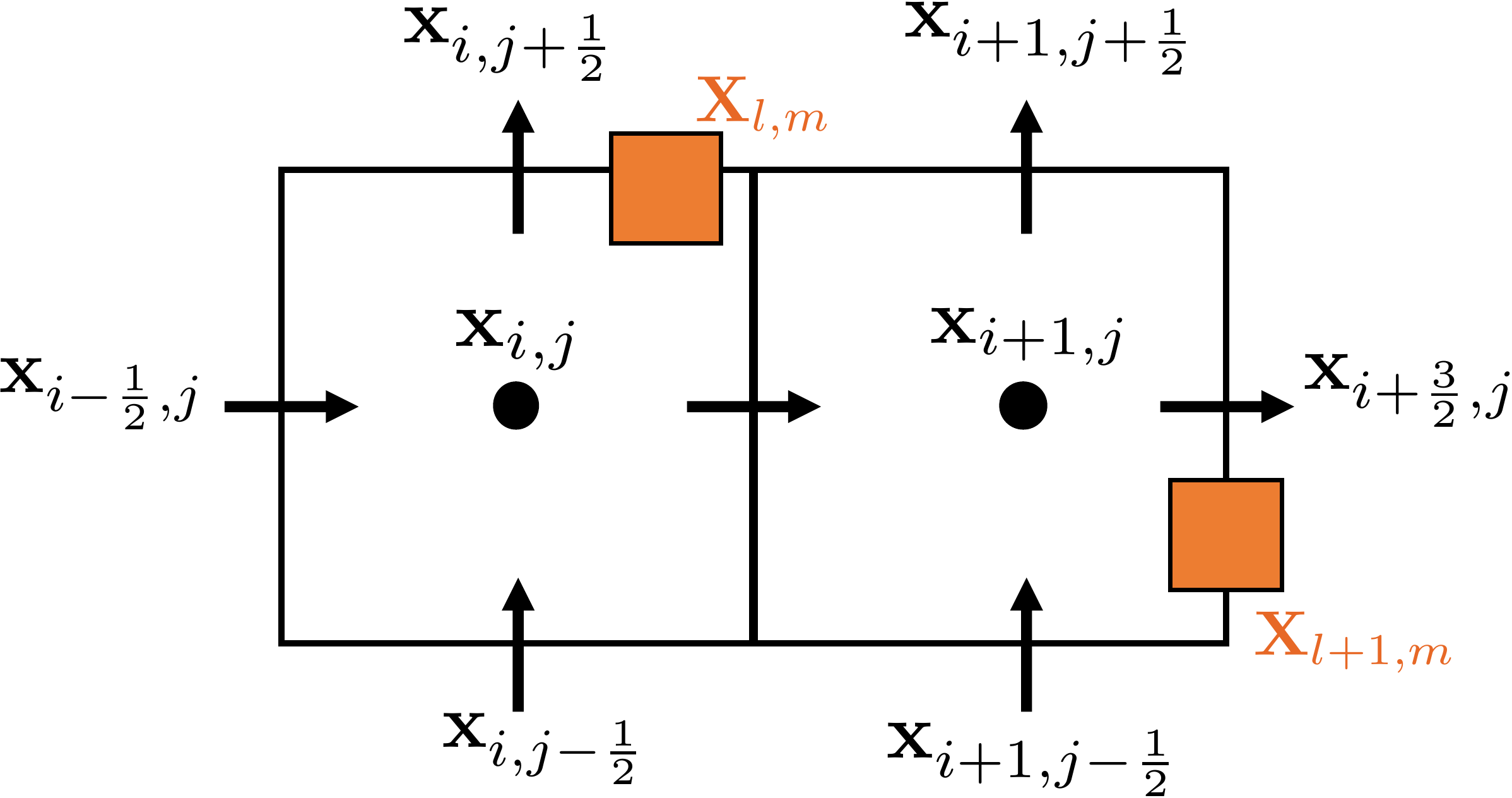}
    \label{grid_cell}
  }
  \caption{\REVIEW{
  \subref{ib_continuous} Sketch of the immersed structure interacting with liquid and 
  gas phases in a rectangular domain. 
  \subref{bp_discrete} Numerical discretization of the domain $\Omega$ into
  Eulerian grid cells ($\blacksquare$, purple) and the indicator function $\chi$ used in the FD/BP method to differentiate between the fluid and solid regions; $\chi = 1$ inside the 
  structure domain and $\chi = 0$ in liquid and gas domains.
  \subref{ib_discrete} Numerical discretization of the domain $\Omega$ into
  Eulerian grid cells ($\blacksquare$, purple) and Lagrangian markers ($\blacksquare$, orange) for the 
  FD/IB method.
  \subref{grid_cell} Two Cartesian grid cells on which the components of the velocity field $\u$
  are approximated on the cell faces ($\rightarrow$, black); the pressure $p$ and level sets $\phi$ and $\psi$ 
  are approximated on the cell center ($\bullet$, black);
  and the Lagrangian quantities are approximated on the marker points ($\blacksquare$, orange), which can
  arbitrarily cut through the Eulerian grid. }
  }
  
  \label{fig_discrete_diagram}
\end{figure}

Lagrangian quantities such as positions, velocities, and forces are defined on immersed markers 
that are allowed to arbitrarily cut through the background Cartesian mesh (see Fig.~\ref{grid_cell}).
These markers are indexed by $(l,m)$ with curvilinear mesh spacings $(\ds_1,\ds_2)$.

A discrete approximation to any general quantity defined on marker points is described by
$\Phi^n_{l,m} \approx \Phi(\s_{l,m}, t^n) = \Phi(l\ds_1, m\ds_2, t^n)$ at time $t^n$.
More specifically, the position, velocity, and force of a marker point are denoted as 
$\X_{l,m}$, $\U_{l,m}$, and $\F_{l,m}$, respectively. 
Fig.~\ref{ib_discrete} shows a sketch of Lagrangian-Eulerian discretization in two spatial dimensions.

\subsection{Lagrangian-Eulerian interaction for the FD/IB method}
\label{sec_lagIB}
The transfer of quantities between the Eulerian and Lagrangian grids
requires discrete approximations to the velocity interpolation and force spreading integrals described by Eqs.~\eqref{eqn_F_f} and~\eqref{eqn_u_interpolation}.
It is convenient to use short-hand notation to denote these integrals. More specifically, the force spreading 
integral of Eq.~\eqref{eqn_F_f} is denoted by $\f = \cS_h[\X] \, \F$, in which $\cS_h[\X]$ is the discrete 
version of the force-spreading operator. The velocity interpolation integral of Eq.~\eqref{eqn_u_interpolation}  
is denoted by $\U = \cJ_h[\X] \, \u$, in which $\cJ_h[\X]$ is the discrete version of velocity-interpolation operator.
It can be shown that if $\cS_h$ and $\cJ_h$ are taken to be adjoint operators, i.e. $\cS_h = \cJ_h^{*}$, then 
the Lagrangian-Eulerian coupling conserves energy~\cite{Peskin02}.

The discrete velocity interpolation of the staggered grid fluid velocity
onto a specific configuration of Lagrangian markers (i.e. $\U = \cJ_h[\X] \u)$ reads
\begin{align}
U_{l,m} & = \sum_{\x_{i-\half,j} \in \Omega} u_{i-\half,j} \delta_h\left(\x_{i-\half,j} - \X_{l,m}\right)  \dx\dy, \\
V_{l,m} & = \sum_{\x_{i,j-\half} \in \Omega} v_{i,j-\half} \delta_h\left(\x_{i,j-\half} - \X_{l,m}\right)  \dx\dy,
\end{align}
The discrete spreading of a force density (defined on Lagrangian markers) onto faces of the staggered grid
(i.e. $\f = \cS_h[\X] \F$) reads
\begin{align}
(f_{1})_{i-\half,j} & = \sum_{\X_{l,m} \in \Vb} (F_{1})_{l,m} \delta_h\left(\x_{i-\half,j} - \X_{l,m}\right) \ds_1\ds_2, \\
(f_{2})_{i,j-\half} & = \sum_{\X_{l,m} \in \Vb} (F_{2})_{l,m} \delta_h\left(\x_{i,j-\half} - \X_{l,m}\right) \ds_1\ds_2.
\end{align}
In the above expressions, $ \delta_h(\x)$ denotes a regularized version of the two-dimensional Dirac delta function
based on a four-point kernel function~\cite{Peskin02}.
We use the same discrete Dirac delta function for both force-spreading and velocity interpolation operators, which ensures that 
$\cS_h = \cJ_h^{*}$. We refer readers to~\cite{Bhalla13,Peskin02} for more details on various properties \REVIEW{(including the spatial invariance property)} and implementation of the grid transfer operations.

\section{Solution methodology} \label{sec_sol_method}
In this section, we describe the full time-stepping scheme and the fluid-structure interaction algorithms employed for 
the FD/BP and FD/IB methods. We first describe the numerical elements common to both implementations, such as 
material property specification, level set advection and reinitialization, and incompressible Navier-Stokes solver for high density 
ratio multiphase flows. The main difference between the two FD methods is the fluid-structure coupling algorithm, which is detailed thereafter. 

\subsection{Material property specification}
\label{sec_reinit}
As described earlier in Sec.~\ref{sec_cont_ls}, the zero isocontours of $\phi(\x,t)$ and $\psi(\x,t)$ represent the liquid-air interface $\Gamma(t)$ and the boundary of the 
immersed structure $\Sb(t)$, respectively.
Using the signed distance property of $\phi$ and $\psi$, we define smoothed Heaviside functions that are regularized over $\ncells$ grid cells on either 
side of the interfaces (assuming $\dx = \dy)$,

\begin{align}
\widetilde{H}^{\text{flow}}_{i,j} &= 
\begin{cases} 
       0,  & \phi_{i,j} < -\ncells \dx,\\
        \frac{1}{2}\left(1 + \frac{1}{\ncells \dx} \phi_{i,j} + \frac{1}{\pi} \sin\left(\frac{\pi}{ \ncells \dx} \phi_{i,j}\right)\right) ,  & |\phi_{i,j}| \le \ncells \dx,\\
        1,  & \textrm{otherwise},   \label{eq_heaviside_flow}
\end{cases} \\
\widetilde{H}^{\text{body}}_{i,j} &= 
\begin{cases} 
       0,  & \psi_{i,j} < -\ncells \dx,\\
        \frac{1}{2}\left(1 + \frac{1}{\ncells \dx} \psi_{i,j} + \frac{1}{\pi} \sin\left(\frac{\pi}{ \ncells \dx} \psi_{i,j}\right)\right) ,  & |\psi_{i,j}| \le \ncells \dx,\\
        1,  & \textrm{otherwise},  \label{eq_heaviside_body}
\end{cases}
\end{align}
A given material property $\zeta$ (such as $\rho$ or $\mu$) is then prescribed in the whole domain using a two-step process.
First, the material property in the ``flowing" phase is set via the liquid-gas level set function
\begin{equation}
\label{eq_ls_flow}
\zeta^{\text{flow}}_{i,j} = \zeta_\text{l} + (\zeta_\text{g} - \zeta_\text{l}) \widetilde{H}^{\text{flow}}_{i,j}.
\end{equation}
Next, the material property is set on cell centers throughout the computational domain, taking into account the solid 
phase~\footnote{For solid viscosity we use $\mu_s = \mu_l$ following the recommendations described in~\cite{Patel2018,Nangia2019}.}
\begin{equation}
\label{eq_ls_solid}
\zeta_{i,j}^{\text{full}} = \zeta_\text{s} + (\zeta^{\text{flow}}_{i,j} - \zeta_\text{s}) \widetilde{H}^{\text{body}}_{i,j}.
\end{equation}
Without the loss of generality, the liquid phase is represented by the negative values of $\phi$ and the solid phase is represented by 
the negative $\psi$ values.
\REVIEW{Note that in the above equations, we have assumed that the number of transition cells is the same across $\Gamma$ and $\Sb$. This is not a strict requirement of the numerical method, but it is true for all the cases considered in the present work.}

\begin{figure}[]
  \centering
  \subfigure[\REVIEW{Material properties in the ``flowing" phases}]{
    \includegraphics[scale = 0.35]{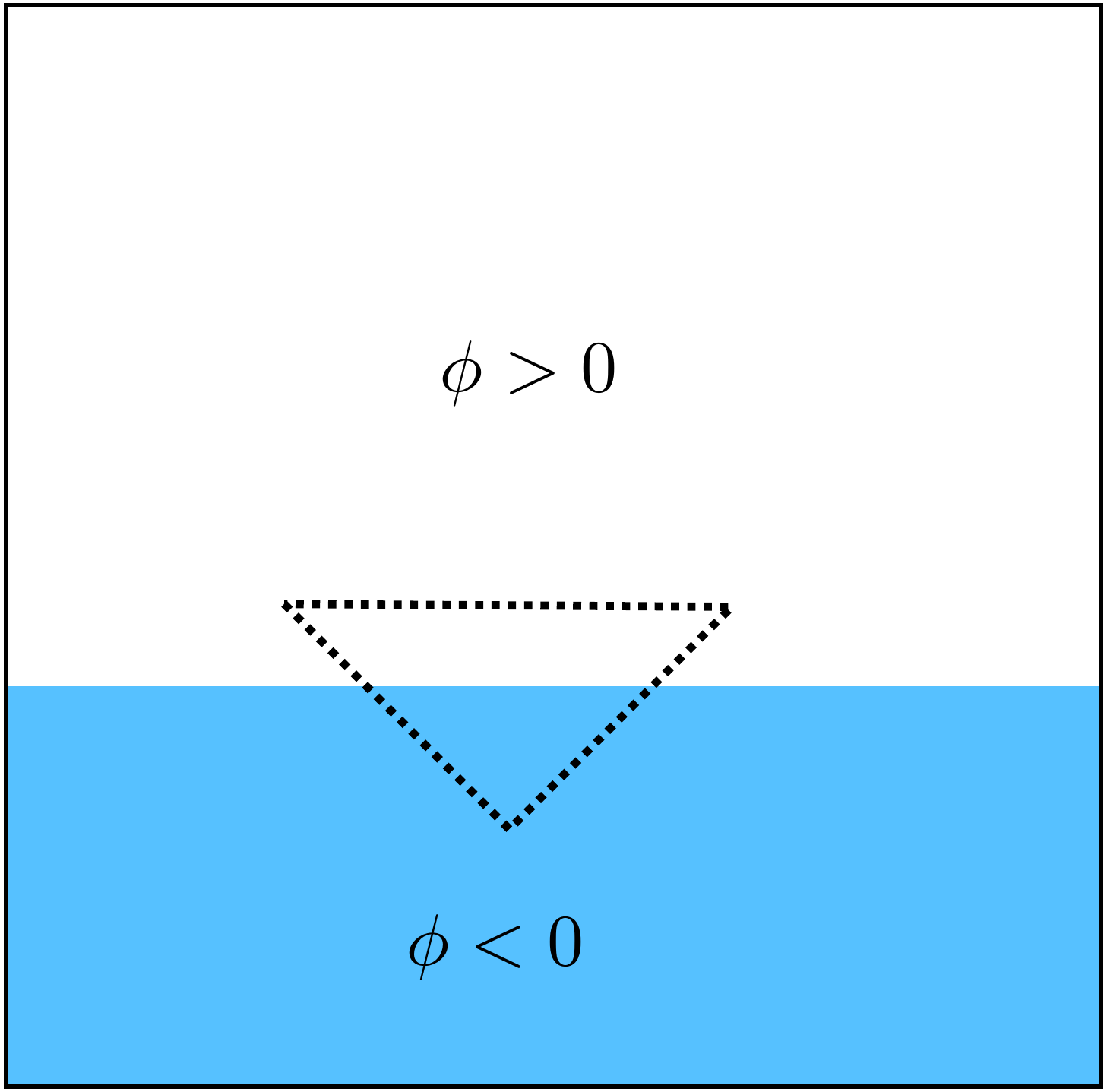}
    \label{ls_flow_diagram}
  }
   \subfigure[\REVIEW{Material properties in the entire domain}]{
    \includegraphics[scale = 0.35]{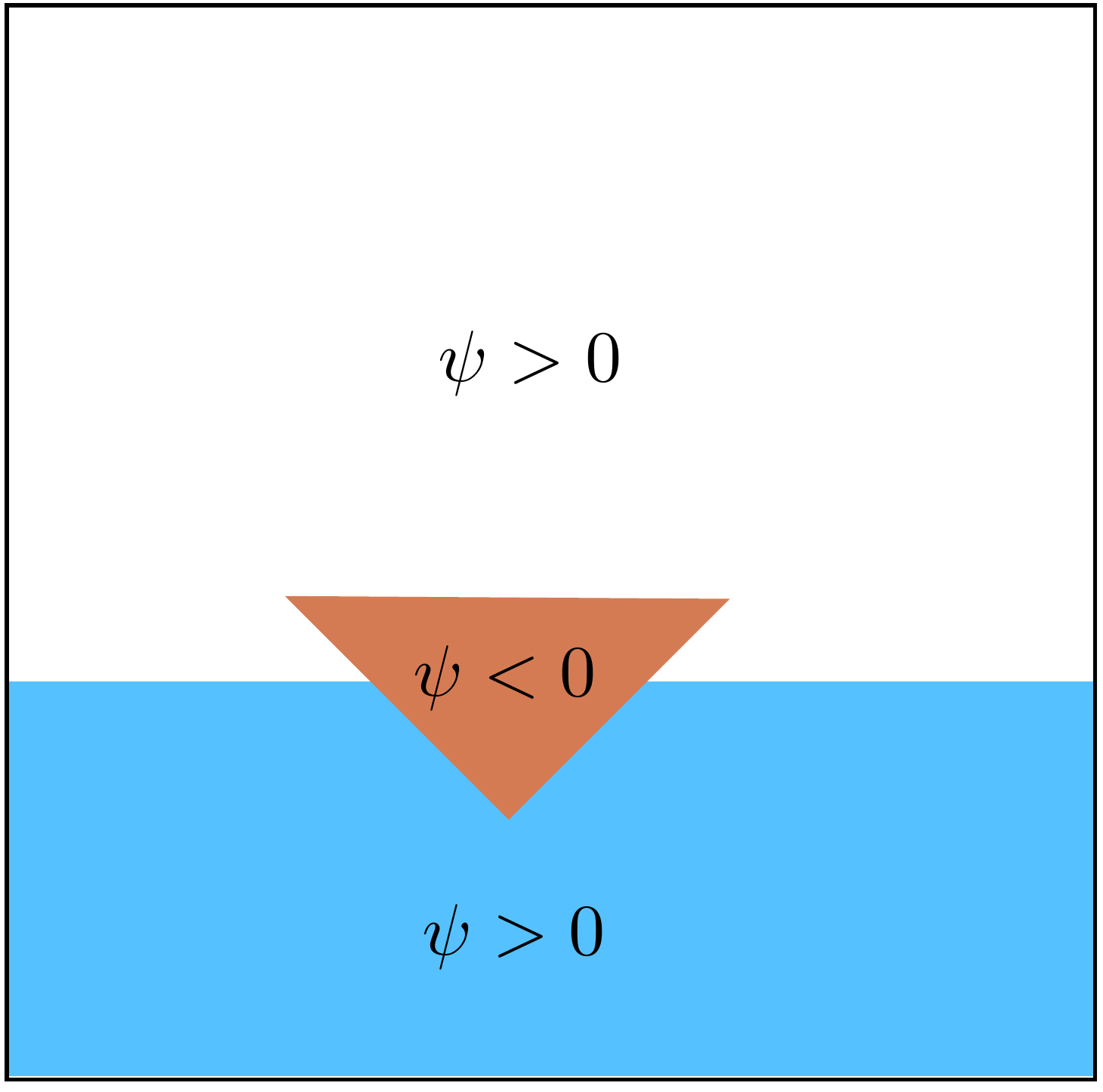}
    \label{ls_solid_diagram}
  }
  \caption{\REVIEW{Sketch of the two-stage process for prescribing the material properties in the computational domain.
  \subref{ls_flow_diagram} Density and viscosity are first prescribed in the ``flowing" phase based on the liquid-gas level set
  function $\phi$ (---, black) and ignoring the body's level set function $\psi$ (\texttt{---}, orange).
  \subref{ls_solid_diagram} Density and viscosity are then corrected in the solid phase}}
  \label{fig_ls_diagram}
\end{figure}

\REVIEW{In general, the signed distance property of $\phi$ and $\psi$ is not preserved under advection governed by
 Eqs.~\eqref{eq_ls_fluid_advection} and~\eqref{eq_ls_solid_advection}. Therefore, a reinitialization process is carried out 
 at the end of each time step such that $\phi^{n+1}$ and $\psi^{n+1}$ represent a signed distance to their respective interfaces. This is described briefly in Appendix~\ref{app_reinit}, and in more detail by Nangia et al.~\cite{Nangia2019}.}  

\subsection{Full time stepping scheme}
\label{sec_temporal_scheme}
We now describe the general time stepping scheme employed over the time interval
$\left[t^n, t^{n} + \dt\right] = \left[t^n, t^{n+1}\right]$; within each time step, $\ncycles$ cycles of fixed-point iteration are used.
In the present work, we always use $\ncycles = 2$. Note that $k$ appears
as a superscript to distinguish the cycle number.
At the beginning of each time step we set $k = 0$, 
with $\u^{n+1,0} = \u^{n}$, $p^{n+\half,0} = p^{n-\half}$,
$\phi^{n+1,0} = \phi^{n}$, $\psi^{n+1,0} = \psi^{n}$, and $\X^{n+1,0} = \X^{n}$. At the first time step $n = 0$,
these quantities are prescribed initial conditions. The midpoint, time-centered approximation
to Lagrangian positions is given by $\X^{n+\half,k} = \half\left(\X^{n+1,k} + \X^{n}\right)$.

\subsubsection{Scalar advection}
\label{sec_scalar_adv}
The level set Eqs.~\eqref{eq_ls_fluid_advection}
and~\eqref{eq_ls_solid_advection} are discretized using a standard time-stepping approach as
\begin{align}
\frac{\phi^{n+1,k+1} - \phi^{n}}{\dt} + Q\left(\u^{n+\half,k}, \phi^{n+\half,k}\right) &= 0, \label{eq_dis_ls_fluid}\\
\frac{\psi^{n+1,k+1} - \psi^{n}}{\dt} + Q\left(\u^{n+\half,k}, \psi^{n+\half,k}\right) &= 0,  \label{eq_dis_ls_solid}
\end{align}
in which $Q(\cdot,\cdot)$ represents a discretization of the linear advection term on cell centers via
an explicit piecewise parabolic method. More specifically, the xsPPM7-limited version described in~\cite{Griffith2009,Rider2007} is employed.
Homogenous Neumann boundary conditions are enforced for $\phi$ and $\psi$ on $\partial \Omega$
using a standard ghost cell treatment~\cite{Harlow1965}.

\subsubsection{Incompressible Navier-Stokes solver: Conservative and consistent transport formulation} \label{sec_cons_ins}
The \REVIEW{conservative form of the} incompressible Navier-Stokes equations Eqs.~\eqref{eqn_momentum}
and~\eqref{eqn_continuity} are discretized as
\begin{align}
	&\frac{\breve{\V \rho}^{n+1,k+1} \u^{n+1,k+1} - { \V \rho}^{n} \u^n}{\dt} + \C^{n+1,k} = -\grad_h p^{n+\half, k+1}
	+ \left(\L_{\mu} \u\right)^{n+\half, k+1}
	+  \V \wp^{n+1,k+1}\g + \thetaFD \fc^{n+1,k+1}, \label{eq_c_discrete_momentum}\\
	& \grad_h \cdot \u^{n+1,k+1} = 0 \label{eq_c_discrete_continuity},
\end{align}
in which $\thetaFD = 1$ for the FD/BP method and $\thetaFD = 0$ for the FD/IB method, i.e., the Eulerian constraint forces are \emph{included} when solving the INS momentum equation for the FD/BP method and \emph{ignored} for the FD/IB method.
The reason for ignoring constraint forces in the FD/IB method will become apparent later when we discuss its FSI algorithm. Similarly,
the specific value of density $\V \wp$ field used to compute the gravitational body force $\V \wp \g$ will be explained in the context
of each fluid-structure coupling algorithm.

\REVIEW{Note that $\left(\L_{\mu} \u\right)^{n+\half, k+1} =  \half\left[\left(\L_{\mu} \u\right)^{n+1,k+1} + \left(\L_{\mu} \u\right)^n\right]$
is a semi-implicit approximation to the viscous strain rate with
$\left(\L_{\mu}\right)^n = \grad_h \cdot \left[\mu^{n} \left(\grad_h \u + \grad_h \u^T\right)^n\right]$, making the above time-stepping scheme with $\ncycles = 2$ resemble a combination of
Crank-Nicolson for the viscous terms and explicit midpoint rule for the convective term.
The newest approximation to viscosity $\mu^{n+1,k+1}$ is obtained via the two-stage process
described in Eqs.~\eqref{eq_ls_flow} and~\eqref{eq_ls_solid}. The newest approximation to density 
$\breve{\V \rho}^{n+1,k+1}$ and the discretization of the convective term $\C^{n+1,k}$ are computed such that
they satisfy consistent mass/momentum transport, which is required to maintain numerical stability
for air-water density ratios. We briefly describe this approach in Appendix~\ref{app_convective} and refer
the reader to Nangia et al.~\cite{Nangia2019} for more details.}
	
\subsubsection{Fluid-structure coupling: FD/BP method} \label{sec_ts_bp}

In the fictitious domain Brinkman penalization formulation, we retain the constraint force in the momentum equation by
taking $\thetaFD = 1$. The constraint/penalization force enforcing the rigid body motion is proportional to the difference 
between the desired structure velocity and the fluid velocity. For the time-stepping scheme it reads
\begin{align}
\fc^{n+1,k+1} = \frac{\widetilde{\chi}}{K}\left(\ub^{n+1,k+1} - \u^{n+1,k+1}\right),  \label{eq_bp_discrete}
\end{align}
in which  $\widetilde{\chi} = 1 - \widetilde{H}^{\text{body}}$, \REVIEW{$ \widetilde{H}^{\text{body}}$ is the regularized structure Heaviside 
function (Eq.~\eqref{eq_heaviside_body}) and $K \sim \cO(10^{-8})$; this is sufficiently small to enforce the rigidity constraint 
in the structure domain, as described by prior studies~\cite{Bergmann11,Verma2017,Gazzola2011}}.  The rigid body velocity $\ub$ in Eq.~\eqref{eq_bp_discrete} 
can be expressed in terms of the translational $\Ur$ and rotational $\Wr$ center of mass velocities
\begin{equation}
\ub^{n+1,k+1} = \Ur^{n+1,k+1} + \Wr^{n+1,k+1} \times \left(\x - \Xcom^{n+1,k+1}\right).
\end{equation} 
The center of mass velocities can be obtained in two distinct ways:

\begin{enumerate}
	\item \emph{Fully prescribed motion}: \\
        For some specified rigid body motion of the structure, i.e, the translational and rotational velocities of the body are known \emph{a priori}, 
        we can directly prescribe the velocity field at time step $n+1$ as
	\begin{equation}
		\label{eq_prescribed_velocity}
		\ub^{n+1,k+1} = \Ur^{n+1} + \Wr^{n+1} \times \left(\x - \Xcom^{n+1}\right).
	\end{equation}
	This algorithm can be used to simulate one-way FSI problems such as flows over stationary bluff bodies
	or structures entering or exiting fluid-gas interfaces with known velocity.

	\item \emph{Free-body motion}: \\
        The rigid body velocity in this case can be obtained by integrating Newton's second law of motion
	\begin{align}
		\Mb \frac{\Ur^{n+1,k+1} - \Ur^n}{\dt} &=  \cF^{n+1,k} + \Mb \g,  \label{eq_newton_u} \\
		\Ib \frac{\Wr^{n+1,k+1} - \Wr^n}{\dt} &=  \cM^{n+1,k},  \label{eq_newton_w}
	\end{align}
	in which $\Mb$ is the mass, $\Ib$ is the moment of inertia, $\cF$ is the net hydrodynamic force, $\cM$ is the net 
	hydrodynamic torque and $\Mb \g$ is the net gravitational force acting on the body. Eqs.~\eqref{eq_newton_u} and~\eqref{eq_newton_w} 
	are integrated using a forward-Euler scheme to compute $\Ur^{n+1,k+1}$, $\Wr^{n+1,k+1}$ and $\Xcom^{n+1,k+1}$. 
	In practice we employ quaternions to integrate Eq.~\eqref{eq_newton_w} in the initial reference frame, which avoids 
	recomputing $\Ib$ as the body rotates \REVIEW{in a complex manner in three spatial dimensions}. 
\end{enumerate}
 
We remark that since the gravitational force is included in the rigid body equation of motion~\eqref{eq_newton_u}, it is not necessary to 
include the volumetric gravity term $\V{\rho}_\text{s} \g$ in the momentum equation~\eqref{eq_c_discrete_momentum}. 
In fact, it is advantageous to just use $\V{\rho}^{\text{flow}} \g$ to avoid spurious velocity currents near the fluid-solid interface 
due to high density gradients~\cite{Nangia2019}. Similar arguments hold for the prescribed motion case. Therefore, we 
use $\V\wp \g = \V{\rho}^{\text{flow}} \g$ in Eq.~\eqref{eq_c_discrete_momentum} for the FD/BP method.

\begin{figure}[]
  \centering
    \includegraphics[scale = 0.4]{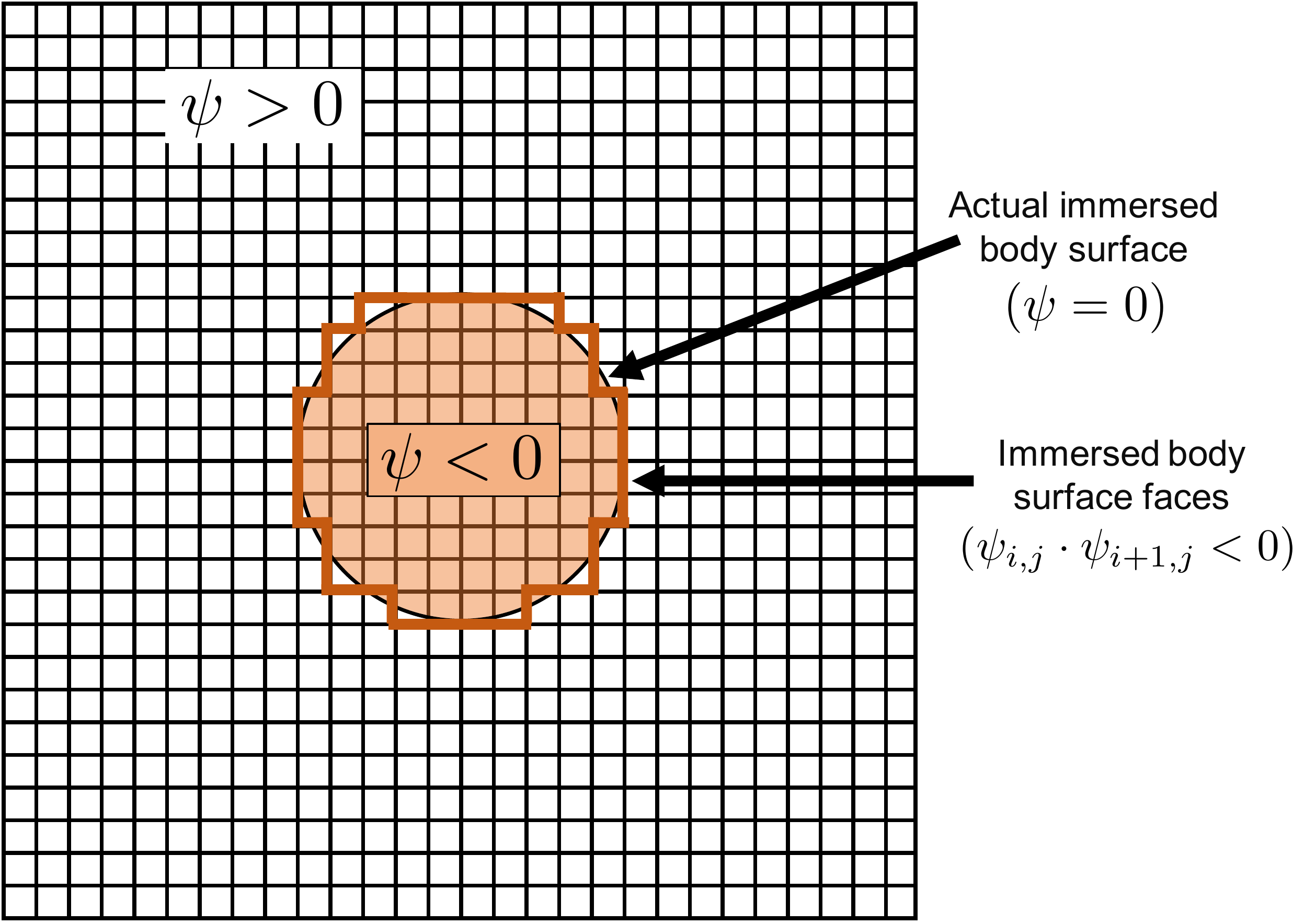}
    \label{hydro_force_cells}
  \caption{Discrete, stair step representation of the body's surface $\Sb$ on a Cartesian grid.
  For two adjacent cells with structure level set $\psi$ values of opposite sign,
  the common face with normal vector $\n_f$ and surface area $\Delta A_f$ is used to evaluate the hydrodynamic force and torque integrals.
}
  \label{hydro_force_cells_diagram}
\end{figure}

The hydrodynamic forces $\cF$ and torques $\cM$ acting on the body are calculated by directly summing pressure and viscous forces from the surrounding fluid on the areal elements of the body surface
\begin{align}
\cF^{n+1,k}  &= \sum_f    \left(-p^{n+1,k} \n_f + \mu_f \left(\grad_h \u^{n+1,k} +  \left(\grad_h \u^{n+1,k}\right)^T \right)\cdot \n_f \right) \Delta A_f,   \label{eq_int_fh} \\
\cM^{n+1,k}  &= \sum_f  \left(\x - \Xcom^{n+1,k}\right) \times  \left(-p^{n+1,k} \n_f + \mu_f \left(\grad_h \u^{n+1,k} +  \left(\grad_h \u^{n+1,k}\right)^T \right)\cdot \n_f \right) \Delta A_f.  \label{eq_int_mh}
\end{align}
In general, the surface boundary $\Sb$ does not conform to the underlying Eulerian grid. Therefore, the above sums are evaluated by representing 
the body's surface in a stair step manner using the grid cells adjacent to $\Sb$. Fig.~\ref{hydro_force_cells_diagram} depicts such a representation. More specifically, a grid face is considered to be a part of the structure boundary if the two grid cells
containing it have structure level set $\psi$ values of \emph{opposite sign}; the summation $\sum_f$ shown in Eqs~\eqref{eq_int_fh} and~\eqref{eq_int_mh} are over these particular grid faces.
This simplifies the hydrodynamic force and torque computations significantly, since all of the required
quantities are readily available or can be interpolated (by simple averaging) onto the Cartesian cell faces
(e.g. face-centered pressure or viscosity $\mu_f$). This is one of the main advantages of fictitious domain 
methods over sharp interface methods: the solution variables are valid on either side of the interface. The latter methods require one-sided interpolations using computational 
geometry constructs.  

\subsubsection{Fluid-structure coupling: FD/IB method} \label{sec_ts_ib}
In the fictitious domain immersed boundary formulation, we ignore the constraint forces in the momentum equation by
\REVIEW{setting} $\thetaFD = 0$.  Therefore, the velocity field computed by the flow solver using 
Eqs.~\eqref{eq_c_discrete_momentum} and~\eqref{eq_c_discrete_continuity} will not satisfy the 
rigid body motion constraints placed in the structure domain. However, the velocity field will be correct in the fluid domain.
If $\widetilde{\u}^{n+1,k+1}$ denotes the velocity solution obtained by ignoring the 
constraint forces, then to correct the velocity in $\Vb(t)$ to $\u^{n+1,k+1}$, we carry out the following \emph{projection} step~\cite{Bhalla13}
\begin{equation}
\label{eq_constraint_projection}
\breve{\V \rho}^{n+1,k+1} \left(\frac{\u^{n+1,k+1} - \widetilde{\u}^{n+1,k+1}}{\dt}\right) = \fc^{n+1,k+1}.
\end{equation}
Similar to Brinkman penalization, the constraint force can be computed using the difference between 
two velocity fields: the desired body velocity and the interpolated uncorrected fluid velocity on the Lagrangian mesh 
$\delU^{n+1,k+1}_{l,m}$ 
\begin{align}
\fc^{n+1,k+1} &= \frac{\breve{\V \rho}^{n+1,k+1}}{\dt}\cS_h\left[\X^{n+\half, k}\right] \delU^{n+1,k+1} \nonumber \\
& = \frac{\breve{\V \rho}^{n+1,k+1}}{\dt}\cS_h\left[\X^{n+\half, k}\right] \left(\Ub^{n+1,k+1} - \cJ_h\left[\X^{n+\half,k}\right]\widetilde{\u}^{n+1,k+1}\right) \label{eq_constraint_force},
\end{align}
which vanishes outside the structure domain.
By correcting the fluid velocity in this way, we ensure that the Eulerian velocity in $\Vb(t)$ approximately matches that
of the solid's Lagrangian velocity $\Ub^{n+1,k+1}$.
Combining the above two equations
yields a simplified update equation for the Eulerian velocity field
\begin{equation}
\label{eq_velocity_correction}
\u^{n+1,k+1} = \widetilde{\u}^{n+1,k+1} + 
\cS_h\left[\X^{n+\half, k}\right] \left(\Ub^{n+1,k+1} - \cJ_h\left[\X^{n+\half,k}\right]\widetilde{\u}^{n+1,k+1}\right).
\end{equation}
We note that there is no guarantee that this corrected velocity will satisfy the divergence-free condition discretely;
it is likely that $\grad_h \cdot \u^{n+1,k+1} \ne 0$. However, we have found that an additional divergence-free velocity projection is \emph{not} necessary to obtain physically accurate results,
corroborating previous investigations by Bhalla et al~\cite{Bhalla13}.

To compute $\fc^{n+1,k+1}$, we first determine $\Ub^{n+1,k+1}$ in the Lagrangian frame.
The rigid body velocity of each Lagrangian marker can be written as 
(omitting the time superscripts)
\begin{equation}
(\Ub)_{l,m} = \Ur + \Wr \times \R_{l,m},
\end{equation}
in which $\R_{l,m} = \X_{l,m} - \Xcom$ is the radius vector pointing from the center of mass to the
Lagrangian marker position. Again considering the two FSI scenarios:
\begin{enumerate}
	\item \emph{Fully prescribed motion}: \\
	For problems in which the motion of the body is known \emph{a priori} as a function of time, we can directly
	prescribe the Lagrangian velocity field at time step $n+1$ as
	\begin{equation}
		\label{eq_prescribed_velocity}
		(\Ub)^{n+1,k+1}_{l,m} = \Ur^{n+1} + \Wr^{n+1} \times \R^{n+\half,k}_{l,m},
	\end{equation}
	which can be used to update the positions of the Lagrangian markers
	\begin{equation}
		\label{eq_prescribed_position}
		\X_{l,m}^{n+1,k+1} = \X_{l,m}^n + \dt (\Ub)^{n+\half,k+1}_{l,m}.
	\end{equation}

	\item \emph{Free-body motion}: \\
	For fully coupled problems in which the body moves as a result of the fluid-structure interaction,
	the Lagrangian velocity field at time step $n+1$ is determined by \emph{redistributing} the
	linear and angular momentum~\cite{Patankar2000,Shirgaonkar2009,Bhalla13} in the structure domain
	\begin{align}
		\Mb \Ur^{n+1,k+1} &= \sum_{\X_{l,m} \in \Vb} \rhos \left(\cJ_h\left[\X^{n+\half,k}\right]\widetilde{\u}^{n+1,k+1}\right)_{l,m} \ds_1\ds_2, \label{eq_linear_conservation} \\
		\Ib \Wr^{n+1,k+1} &= \sum_{\X_{l,m} \in \Vb} \rhos \R^{n+\half,k}_{l,m} \times \left(\cJ_h\left[\X^{n+\half,k}\right]\widetilde{\u}^{n+1,k+1}\right)_{l,m} \ds_1\ds_2. \label{eq_angular_conservation}
	\end{align}
	The structure's velocity
	and position are then updated via Eqs.~\eqref{eq_prescribed_velocity} and~\eqref{eq_prescribed_position}.
\end{enumerate}

Note that for fully prescribed motion, the gravitational force does 
not affect the body's (specified) velocity. To avoid spurious flow currents near the fluid-structure interface (due to 
 sharp density gradients) we simply use $\V \wp \g = \V{\rho}^{\text{flow}} \g$ for prescribed motion FSI problems. The free-body 
 motion scenario requires special consideration. Because a momentum redistribution procedure is employed to obtain 
 $\Ur$ and $\Wr$, the algorithm assumes that the uncorrected fluid momentum $\widetilde{\u}$ is obtained by including all 
 (including gravitational) body forces in the momentum equation~\eqref{eq_c_discrete_momentum}. Therefore, the gravitational 
body force should account for the solid density and we use $\V \wp \g = \V{\rho}^{\text{full}} \g$ for 
free-body motion FSI problems. See~\cite{Nangia2019} for more discussion.
 
The net hydrodynamic force $\cF$ and torque $\cM$ for the FD/IB method can be computed as a post-processing step 
using the Lagrangian quantities~\cite{Nangia17}
\begin{align}
\cF^{n+1} &=   \sum_{\X_{l,m} \in \Vb} \rhos \left[\frac{\left(\Ub\right)^{n+1}_{l,m} - \left(\Ub\right)^{n}_{l,m}}{\dt} -  \frac{\delU^{n+1}_{l,m}}{\dt}\right] \ds_1 \ds_2,  \label{eq_lmforce} \\ 
\cM^{n+1} &=   \sum_{\X_{l,m} \in \Vb} \rhos \R^{n+1}_{l,m} \times \left[\frac{\left(\Ub\right)^{n+1}_{l,m} - \left(\Ub\right)^{n}_{l,m}}{\dt} -  \frac{\delU^{n+1}_{l,m}}{\dt}\right] \ds_1 \ds_2, \label{eq_lmtorque}
\end{align}
in which the discrete approximations of the quantities on the right-hand side are readily available during each time step. 

Finally, we remark that the above methodologies assume that all rigid-body degrees of freedom either are \emph{all}
fully prescribed (locked) or \emph{all} undergoing free-body motion (unlocked). In our practical implementation, we are able to mix and match which degrees of freedom are locked and unlocked.
We make use of this flexibility in the numerical examples presented in this work .

\section{Comparison of the two fictitious domain methods} \label{sec_comp_fd}

Below, we list some of the similarities and differences of the two previously described FD algorithms:

\begin{itemize}

\item Both methods extend the fluid momentum equation into the solid domain, which results in a valid solution 
on both sides of the structure interface. 

\item Both methods formulate the constraint force in terms of a difference between the desired body and fluid velocities. For a 
specific value of permeability $K = \Delta t /  \breve{\V \rho}$, identical forms of the constraint force $\fc$ are 
obtained~\footnote{For practical water-entry and water-exit problems, $K  = \Delta t/ \rhol$ is $\cO(10^{-8})$.}. 

\item The FD/BP method treats the constraint force implicitly, whereas the FD/IB method treats it explicitly. The explicit treatment of $\fc$ allows the use of an existing fluid solver without any modifications. The implicit Brinkman penalization method necessarily 
requires changes to an existing fluid solver infrastructure~\footnote{It is also possible to treat $\fc$ explicitly in FD/BP method without requiring modification of an existing fluid solver. We have not yet analyzed the accuracy and stability of explicit FD/BP methods, however.}.  

\item The fully Eulerian nature of the FD/BP method is an attractive feature from a domain decomposition perspective, which can enable parallel scalability. In contrast, special 
load-balancing techniques must be employed to efficiently distribute Lagrangian and Eulerian data used in the FD/IB method.

\item  The momentum redistribution step of the FD/IB method avoids the need to compute hydrodynamic forces and torques on the immersed 
surface explicitly (which is a requirement of the FD/BP algorithm). The former approach requires Eulerian-Lagrangian interpolation and 
spreading routines, which may become expensive for large volumetric bodies. 

\item The FD/IB method works best for volumetric forces that are defined throughout the interior region of the structure. Incorporating point forces 
and torques that act only at certain points of the body (e.g. hinge forces or spring/damper forces) is not straightforward. Such forces 
and torques can be easily incorporated in Newton's law of motion used in the FD/BP method.    

\end{itemize}

%%%%%%%%%%%%%%%%%%%%%%%%%%%%%%
\section{Software implementation} \label{sec_sfw}
The numerical algorithm described here is implemented
within the IBAMR library~\cite{IBAMR-web-page}, which is an open-source C++
simulation software focused on immersed boundary methods with
adaptive mesh refinement. All of the numerical examples presented here
are publicly available via \url{https://github.com/IBAMR/IBAMR}.
IBAMR relies on SAMRAI \cite{HornungKohn02, samrai-web-page} for Cartesian grid 
management and the AMR framework. Linear and nonlinear solver support in IBAMR is provided by 
the PETSc library~\cite{petsc-efficient, petsc-user-ref, petsc-web-page}.
All of the example cases in the present work made use of distributed-memory
parallelism using the Message Passing Interface (MPI) library.
All of the example cases described in this section were carried out using 72 processors on the Fermi cluster at SDSU.

%%%%%%%%%%%%%%%%%%%%%%%%%%%%%%
\section{Numerical examples}
\label{sec_examples}
\REVIEW{We begin by simulating the challenging case of a two-dimensional, freely falling inclined wedge with three free degrees of freedom 
to validate the FD/BP method. We have extensively validated FD/IB implementations in previous work~\cite{Bhalla13,Bhalla13FDO,Nangia17,Kallemov16,Usabiaga17,Nangia2019,Bhalla14EHD,Dombrowski19} in the context
of both two and three phase flows.}  

Next, we simulate water-entry and exit of a freely falling wedge and cylinder in two-spatial 
dimensions, and compare the fluid-structure dynamics obtained from FD/BP and FD/IB methods.  
For these cases, the only unlocked degree of freedom is the vertical ($y$) direction. 

We use $\ncycles = 2$ for both methods. Two grid cells of 
smearing $\ncells = 2$ are used to transition between different material properties on either side of the interfaces. 
Water and air densities are taken to be $1000$ kg/m$^3$ and $1.2$ kg/m$^3$, respectively, and their respective 
viscosities are taken to be $10^{-3}$ Pa$\cdot$s and $1.8 \times 10^{-5}$ Pa$\cdot$s. Surface tensions effects are neglected. 
No-slip boundary conditions are imposed along $\partial \Omega$.
\REVIEW{
\subsection{Water-entry of a free falling inclined wedge} \label{sec_inclined_wedge}
In this section, we consider the case of an inclined, 2D wedge-shaped object impacting an air-water interface. The wedge is 
initially rotated counterclockwise through a heel angle of $5^\circ$ as shown in Fig.~\ref{fig_incl_wedge_schematic}. 
The isosceles triangle body has length $L = 0.61$ m and 
a deadrise angle of $20^\circ$. Its mass and moment of inertia are $124$ kg and $8.85$ kg $\cdot$ m$^2$, respectively.
The structure's dimensions and material properties are chosen to match the experimental study conducted by 
Xu et al.~\cite{Xu1999}. Additionally, this case has been studied numerically using a weakly compressible smoothed particle hydrodynamics (SPH) 
method by Oger et al.~\cite{Oger2006}, and an artificial compressibility method combined with Chimera grid-based Navier-Stokes solver by Nguyen et al.~\cite{Nguyen2016}. That is, the prior numerical studies~\cite{Oger2006,Nguyen2016} simulate the compressible version of the Navier-Stokes equation in the low Mach number regime in contrast to the \emph{incompressible} Navier-Stokes solver used in the present study.  

\begin{figure}[t!]
  \centering
  
    \includegraphics[scale = 0.4]{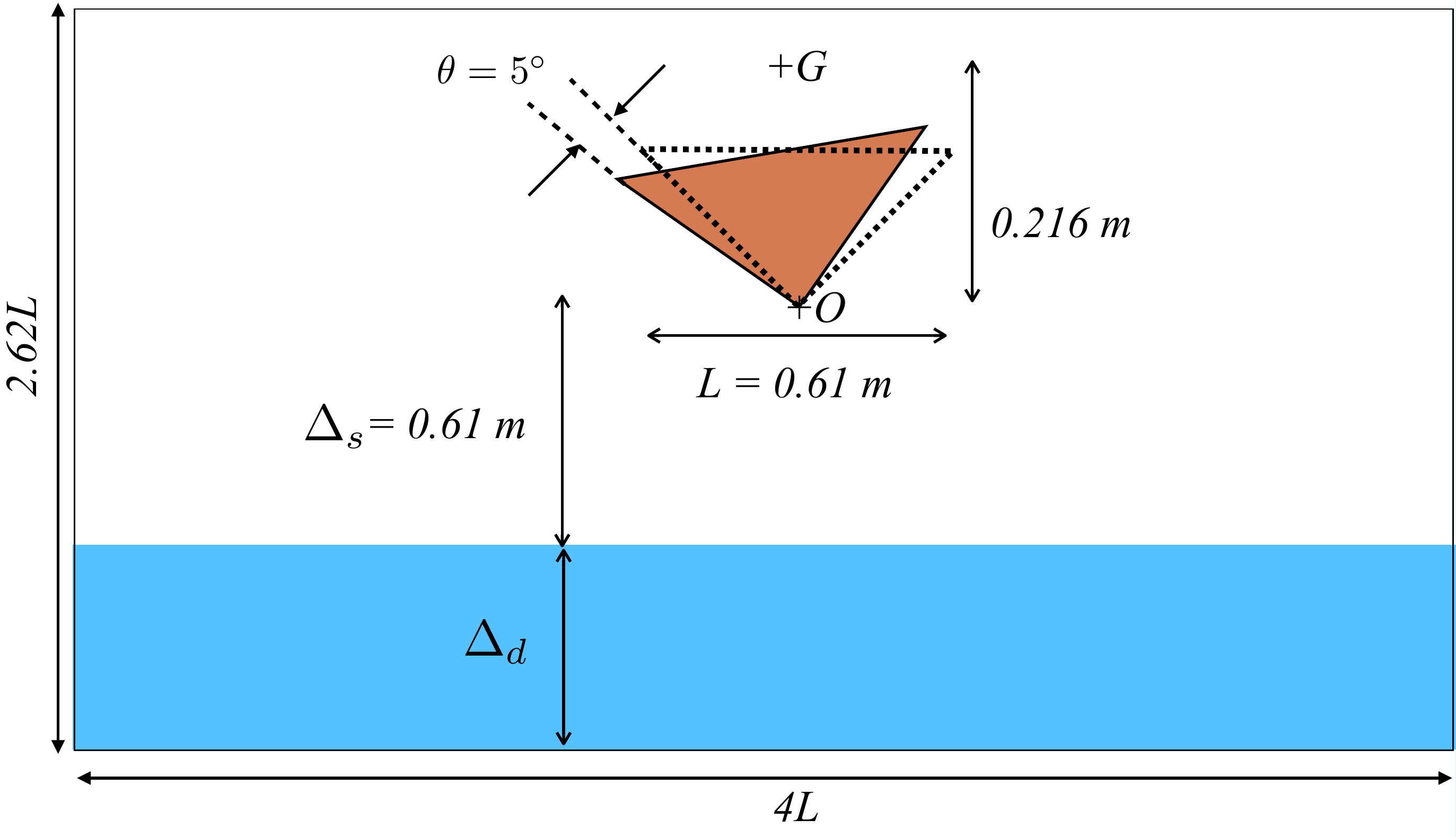}
    
  \caption{\REVIEW{Schematic of a free falling wedge with an initial heel angle of $5^\circ$. $G$ is the center of mass location for the unrotated wedge, while $O$ is the initial location of the bottom vertex of the wedge.
 The distance $O-G = 0.216$ m is the distance between the wedge tip and its COM point. Sketch is not to scale.}
}
  \label{fig_incl_wedge_schematic}
\end{figure}

The computational domain is taken to be $\Omega = [0, 4L] \times [0, 2.62L]$, which is discretized by uniform 
computational grids. The initial distance between the bottom vertex of the wedge (point O in the Schematic~\ref{fig_incl_wedge_schematic}) 
and the air-water interface is $\Delta_s = 0.61$ m.
We note that the domain dimensions and the depth of the initially quiescent water are not mentioned in prior experimental and numerical~\footnote{Only the domain length of $2.62 L$  is mentioned in~\cite{Oger2006,Nguyen2016}.} investigations~\cite{Xu1999,Oger2006,Nguyen2016}. 
Therefore as a preliminary test case, four different water depths
are considered: $\Delta_d = 0.13$ m, $0.225$ m, $0.35$ m, and $0.45$ m and the results are compared to the experimental data
of Xu et al.~\cite{Xu1999}. The domain is discretized by a $488 \times 320$ grid (medium resolution),
and a constant time step size of $\dt = 2.5 \times 10^{-5}$ is used.
Fig.~\ref{fig_incl_wedge2d_h_test} shows the time evolution of vertical acceleration (normalized by $g = 9.81$ m/s$^2$)
and angular acceleration for varying water depths. Note that the peak linear acceleration decreases with increased
water depth~\cite{Ghazizade2013}.
Based on Fig.~\ref{incl_wedge2d_vert_acc_h}, it is evident that the best agreement with mean experimental data is achieved when simulating this problem using a water depth of $\Delta_d = 0.13$ m. However, the wedge itself has a height
of $0.11$ m, making this value of water depth inadequate for simulating long-term dynamics (in particular, experimental data are available up to a final time of around $t = 0.45$ s). Therefore, we choose to use $\Delta_d = 0.225$ m for the remaining cases in this section. Finally we note that although the 
peaks in our simulated angular accelerations match well with experimental data (Fig.~\ref{incl_wedge2d_ang_acc_h}), there are some differences 
over the time interval $t = 0.37$ s to $t = 0.4$ s. We attribute these differences to three possibilities:
\begin{enumerate}
\item Compressibility effects during initial impact. Chen et al.~\cite{Chen2019} found that it can be important to consider compressibility of the water phase for problems involving substantial impact forces. Since our solver assumes incompressibility of liquid phases, fluid oozes from the sides of the object more 
quickly (compared to a compressible fluid), which may explain the faster angular dynamics observed in our simulations.
\item Minor differences in problem set up (such as domain boundaries or initial water depth) could also explain the deviations from the experimental data.
\item Unlike the linear acceleration data, which is reported as a mean of two experiments with a significant standard deviation, the angular acceleration data is reported from a single trial in Xu et al.~\cite{Xu1999}. It is unclear how repeatable these results are (i.e. multiple trials could exhibit large variability in results). In spite of some differences in the angular acceleration, we obtain an excellent match with the experimental data for angle of heel (see Fig.~\ref{incl_wedge2d_heel_angle}).
\end{enumerate}

\begin{figure}[t!]
  \centering
   \subfigure[\REVIEW{Dimensionless vertical acceleration ($\ddot{y}/g)$}]{
    \includegraphics[scale = 0.3]{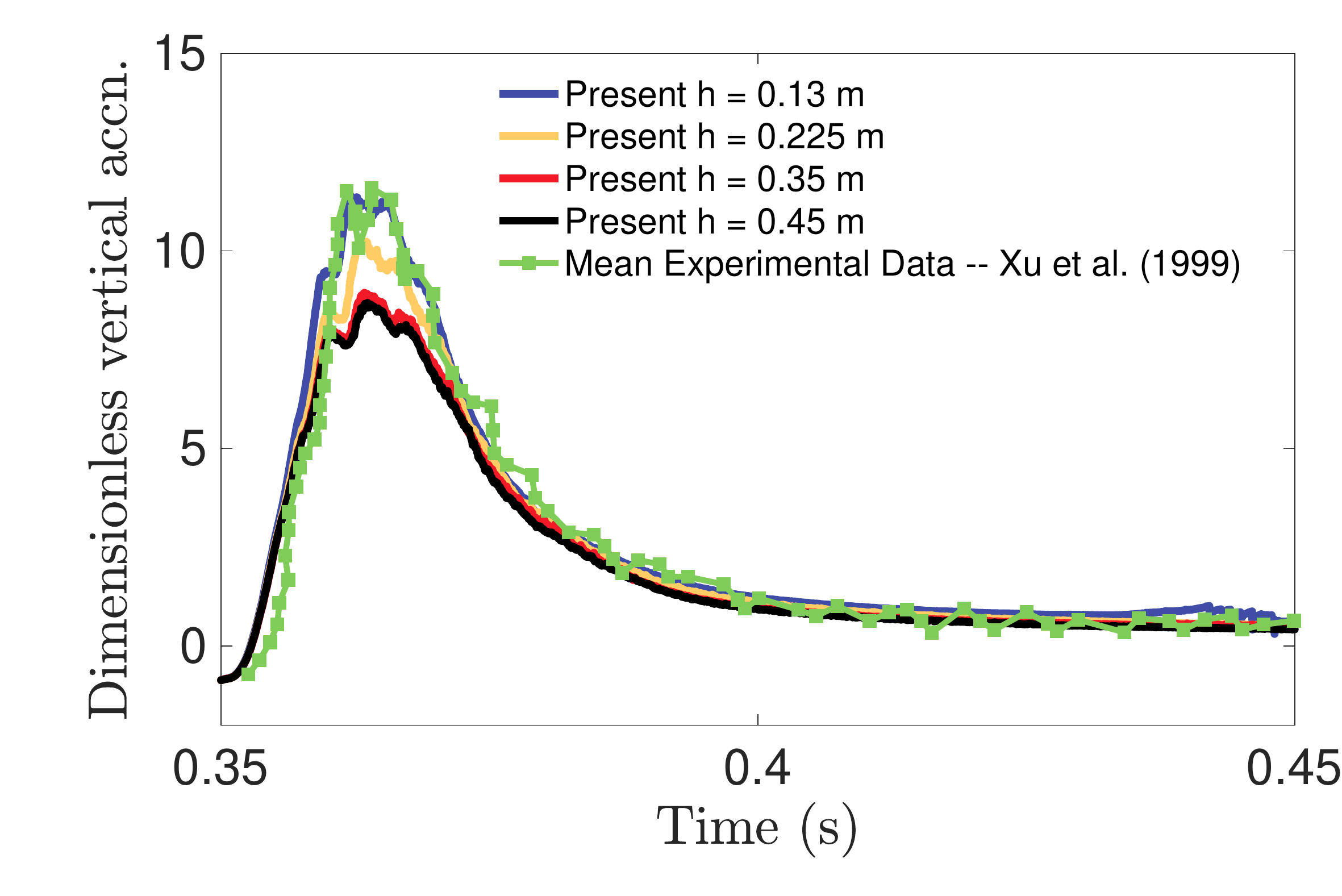}
    \label{incl_wedge2d_vert_acc_h}
  }
     \subfigure[\REVIEW{Angular acceleration}]{
    \includegraphics[scale = 0.3]{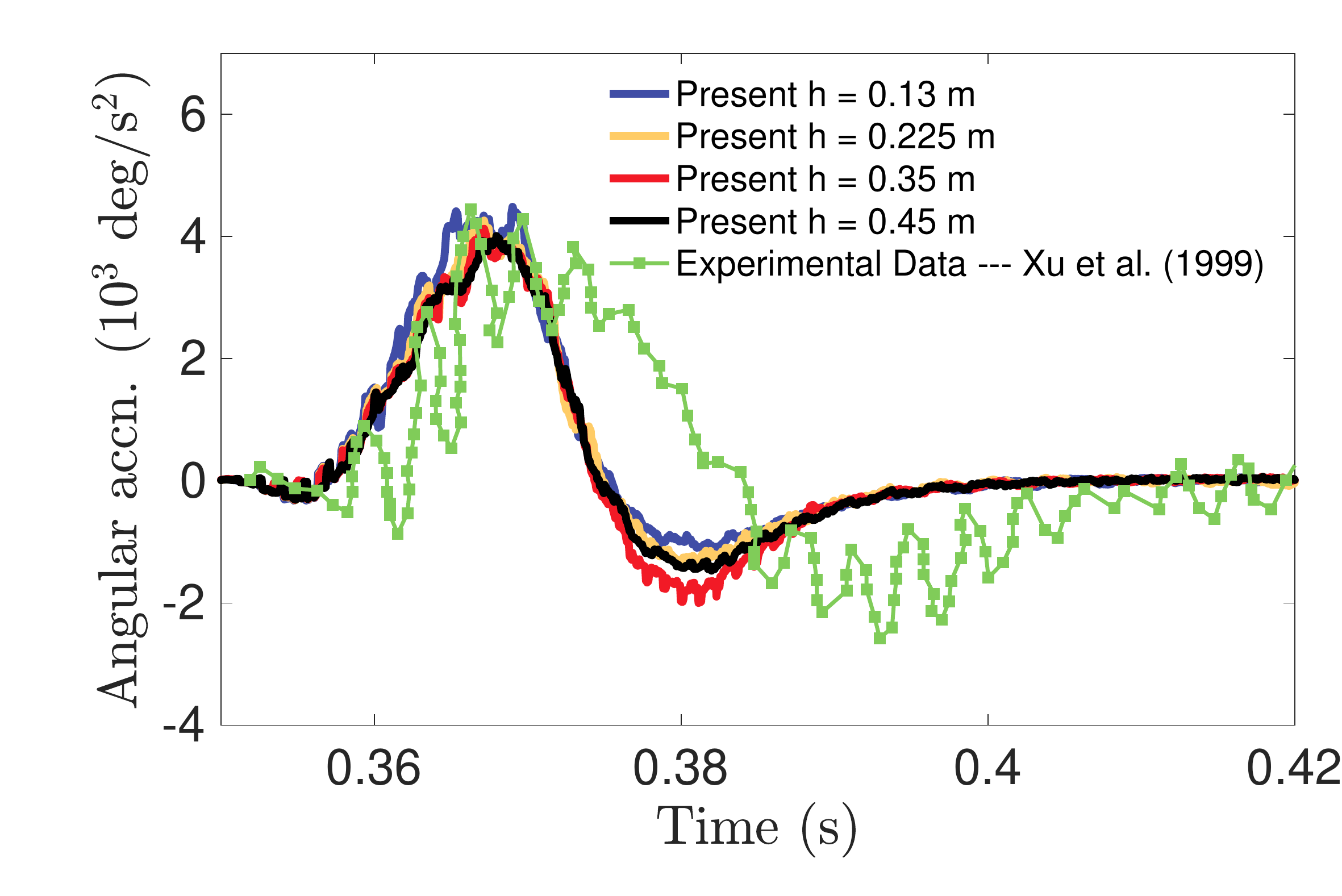}
    \label{incl_wedge2d_ang_acc_h}
  }
  \caption{\REVIEW{Temporal evolution of
  \subref{incl_wedge2d_vert_acc_h} dimensionless vertical acceleration, and 
  \subref{incl_wedge2d_ang_acc_h} angular acceleration for a 2D inclined wedge freely falling into water.
  (---) Present FD/BP simulation data for varying water depths;
  (--$\blacksquare$--, green) experimental data from Xu et al.~\cite{Xu1999}.}
}
  \label{fig_incl_wedge2d_h_test}
\end{figure}

For a fixed water depth $\Delta_d = 0.225$ m, we now simulate this problem with three different mesh resolutions;
a coarse grid ($244 \times 160$ with $\dt = 5 \times 10^{-5})$, a medium grid
($488 \times 320$ with $\dt = 2.5 \times 10^{-5}$) and a fine grid ($814 \times 534$ with $\dt = 1 \times 10^{-5}$)
are considered. Figs.~\ref{incl_wedge2d_heel_angle} and~\ref{incl_wedge2d_velocity} show time evolution of the wedge's
heel angle and vertical velocity, respectively, for all three mesh resolutions. For both quantities the solutions on the finest
grid are in excellent agreement with the experimental data of Xu et al.~\cite{Xu1999}, and grid convergence towards
these data is also seen. In particular, a medium grid resolution that corresponds to approximately 120 grid cells per wedge length is 
adequate to resolve the FSI dynamics of a free-falling wedge. 
In Figs.~\ref{incl_wedge2d_vert_acc} and~\ref{incl_wedge2d_ang_acc}, we show the dimensionless
vertical acceleration and angular acceleration for the finest grid resolution as a function of time.
Although the maximum vertical acceleration from the FD/BP method is slightly smaller than seen in experiments, the overall 
trend and values away from this peak match reasonably well. It is evident that the numerical simulations presented in
Oger et al.~\cite{Oger2006} and Nguyen et al.~\cite{Nguyen2016} suffer from the same mismatch in peak linear acceleration.
The peaks in angular acceleration are in decent agreement with the experimental study, although we again observe 
some differences in the trend over the time interval $t = 0.37$ s to $t = 0.4$ s even at the finest grid resolution. This implies that 
a converged solution for angular acceleration has been achieved for our simulations. The simulations from Oger et al. and Nguyen et al.
are in much better agreement with the experiment, which we attribute to the fact that these authors considered
compressibility effects in their numerical schemes. Finally in Fig.~\ref{fig_incl_wedge2d_moments}, we show the convergence of 
hydrodynamic vertical force and torque for three different grid resolutions. As expected, large impulses are seen as the 
wedge slams into the water just before $t = 0.4$ s.

\begin{figure}[t!]
  \centering
  \subfigure[\REVIEW{Heel angle}]{
    \includegraphics[scale = 0.3]{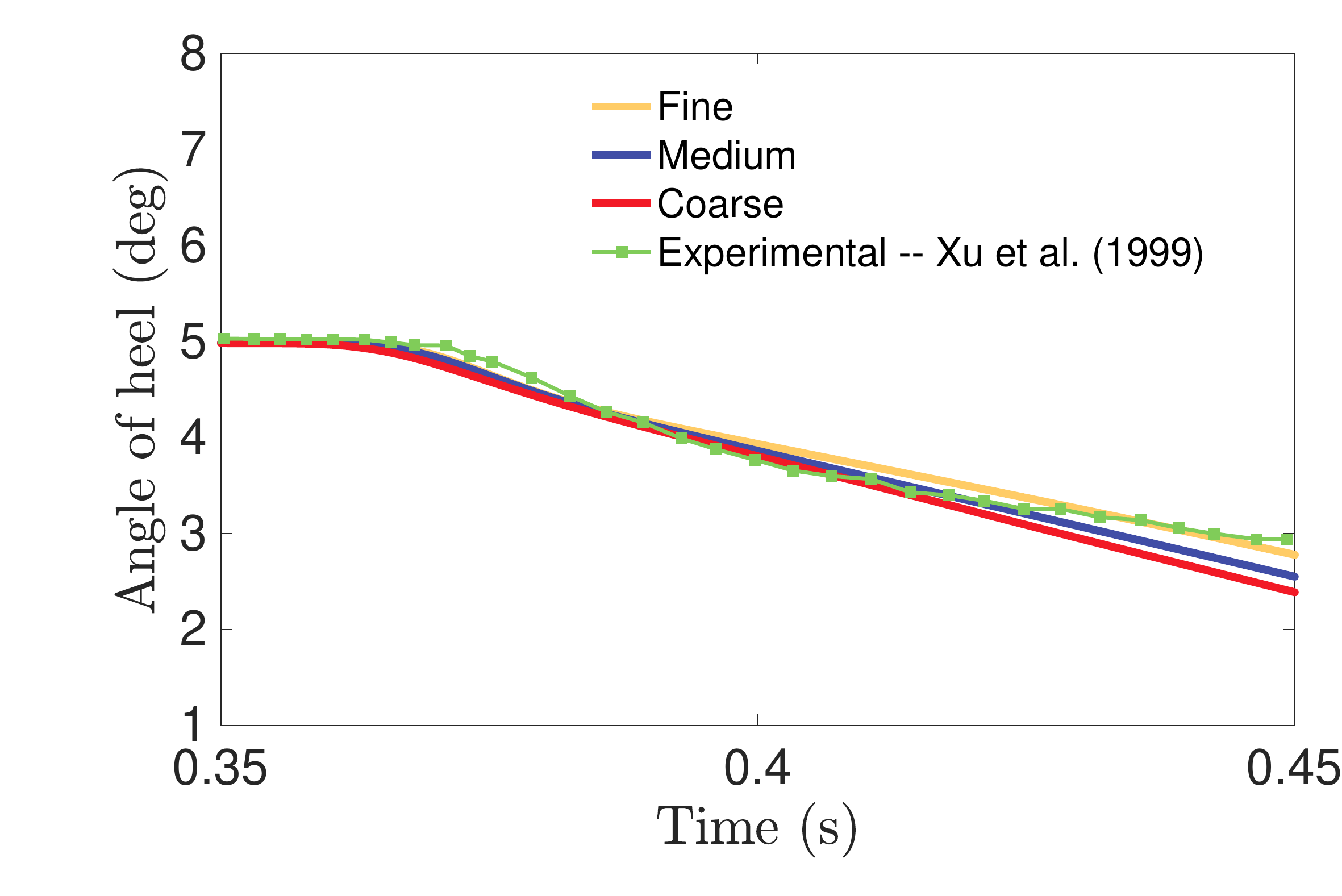}
    \label{incl_wedge2d_heel_angle}
  }
  \subfigure[\REVIEW{Vertical velocity}]{
    \includegraphics[scale = 0.3]{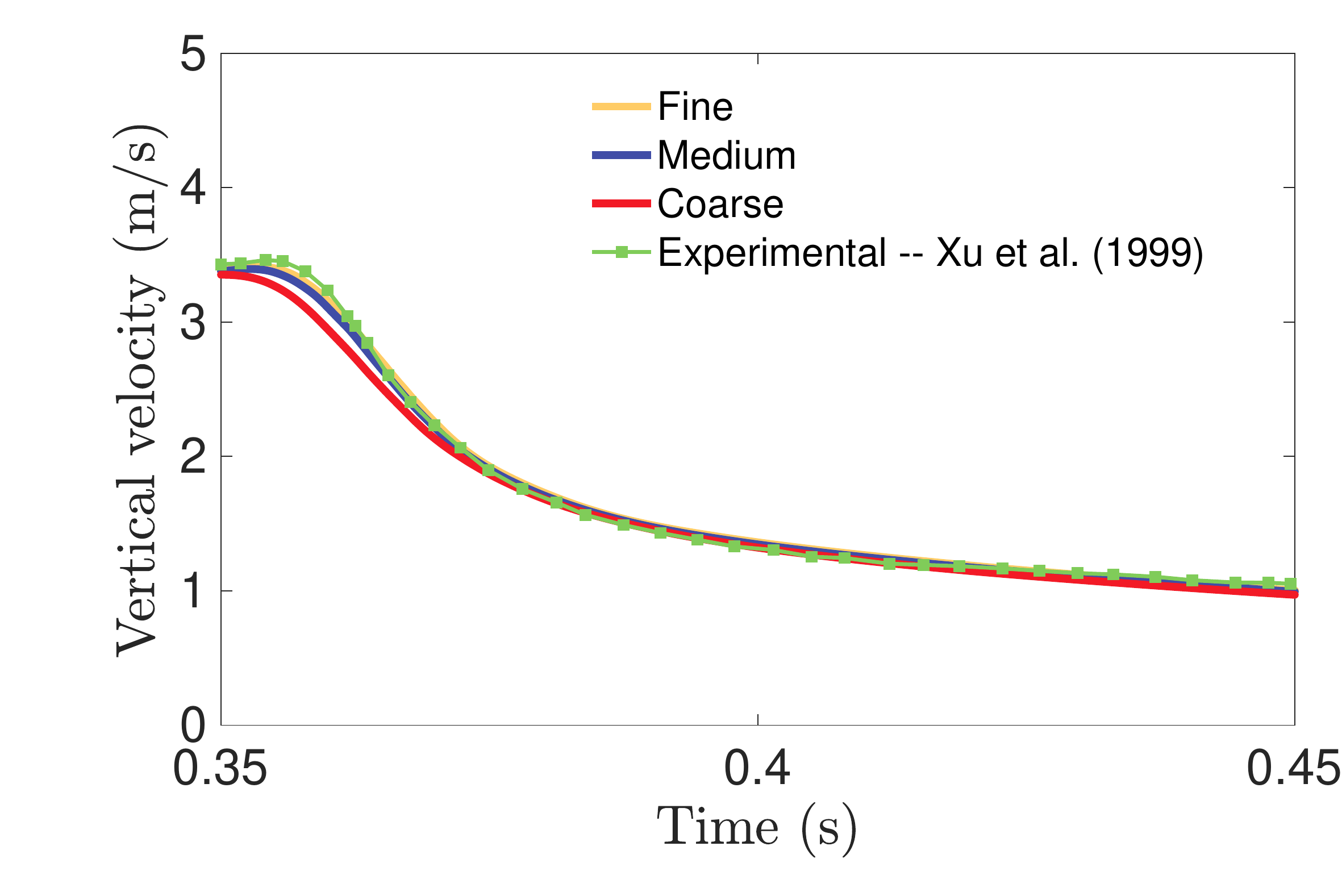}
    \label{incl_wedge2d_velocity}
  }
    \subfigure[\REVIEW{Dimensionless vertical acceleration ($\ddot{y}/g)$}]{
    \includegraphics[scale = 0.3]{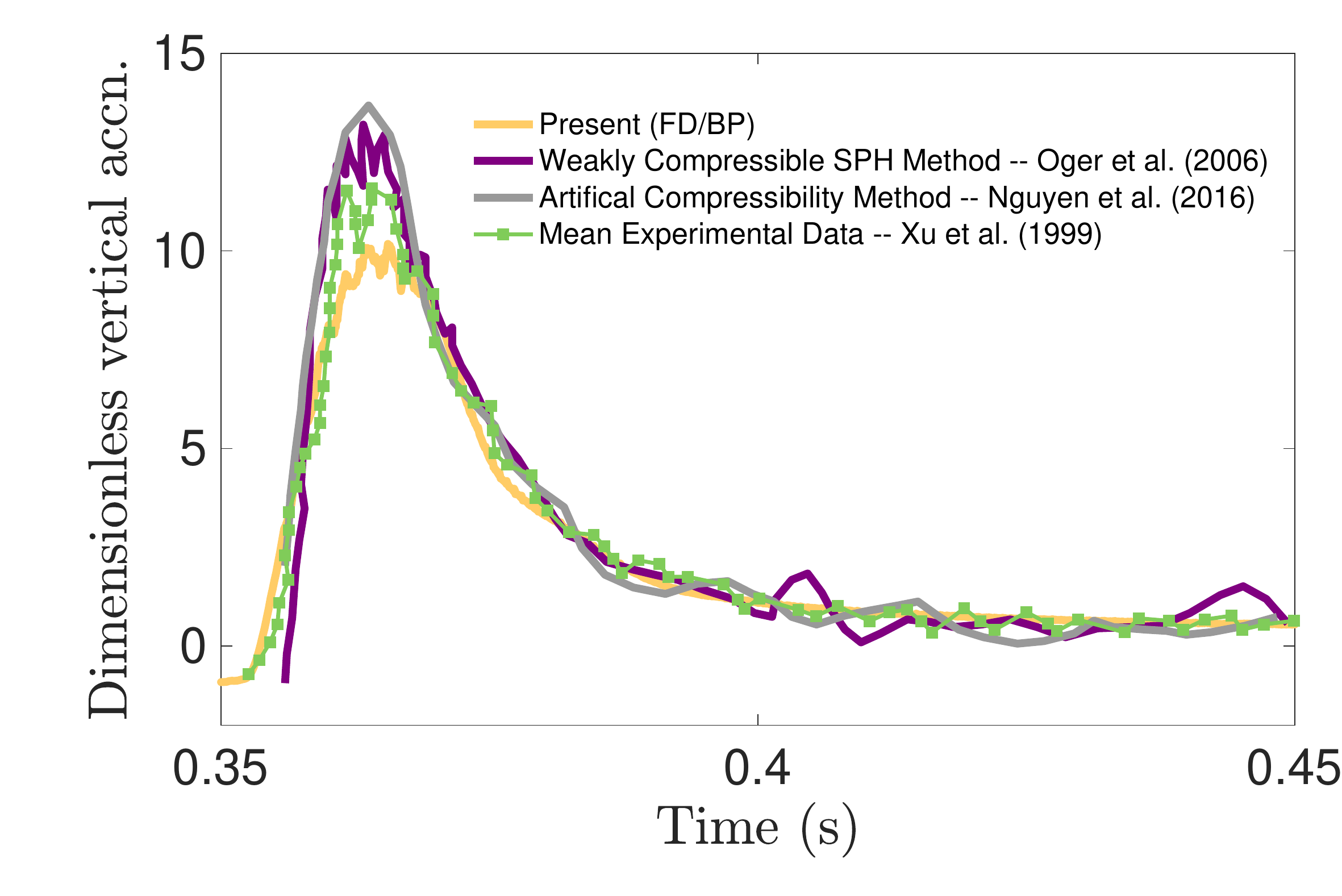}
    \label{incl_wedge2d_vert_acc}
  }
     \subfigure[\REVIEW{Angular acceleration}]{
    \includegraphics[scale = 0.3]{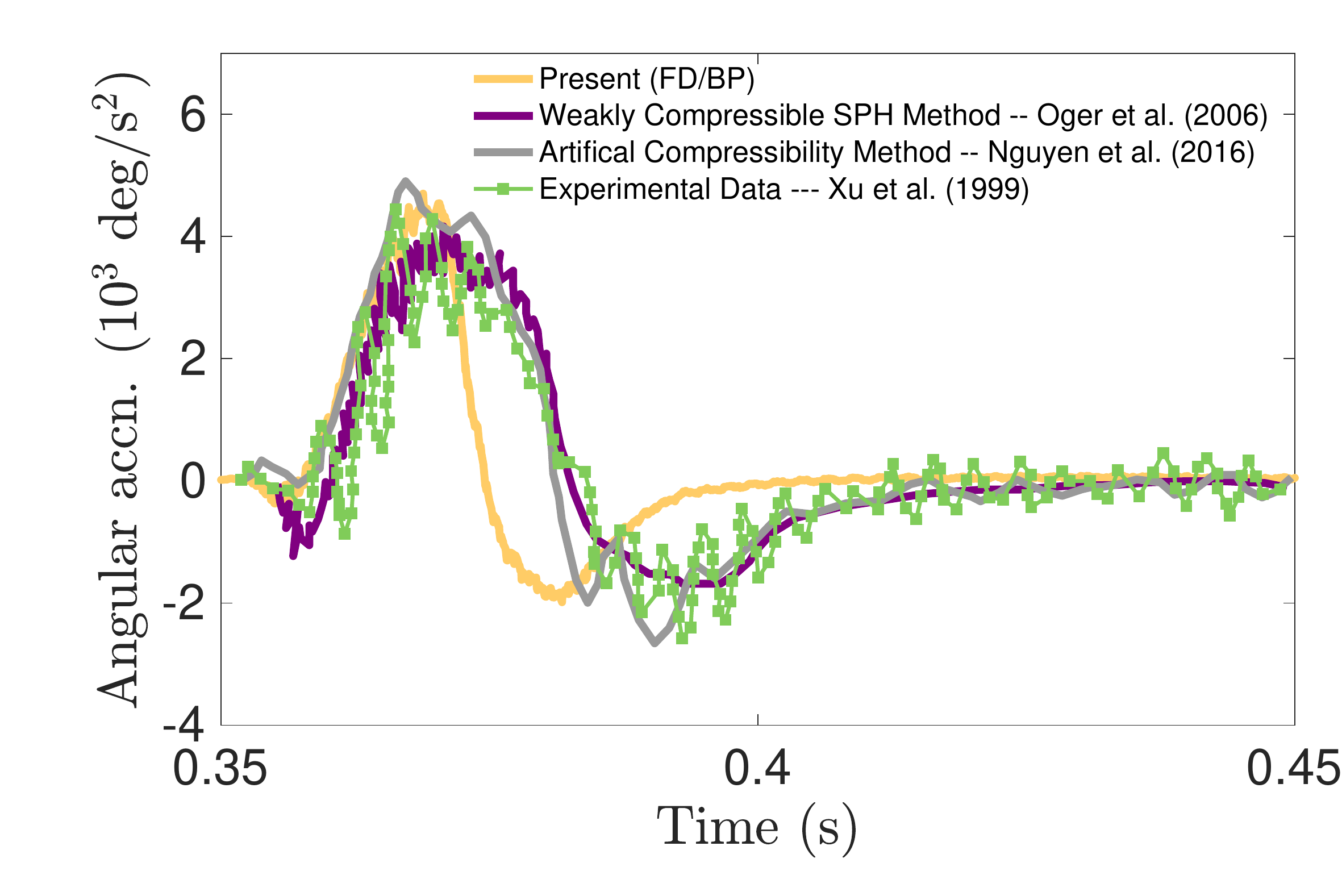}
    \label{incl_wedge2d_ang_acc}
  }

  \caption{\REVIEW{Temporal evolution of
  \subref{incl_wedge2d_heel_angle} heel angle,
  \subref{incl_wedge2d_velocity} vertical velocity,
   \subref{incl_wedge2d_vert_acc} dimensionless vertical acceleration, and
  \subref{incl_wedge2d_ang_acc} angular acceleration, for a 2D inclined wedge freely falling into water.
  (---, yellow) Present FD/BP simulation data for a fine grid resolution $814 \times 534$;
  (---, blue) Present FD/BP simulation data for a medium grid resolution $488 \times 320$;
  (---, red) Present FD/BP simulation data for a coarse grid resolution $244 \times 160$;
  (--$\blacksquare$--, green) experimental data from Xu et al.~\cite{Xu1999};
  (---, purple) 2D simulation data from Oger et al.~\cite{Oger2006};
  (---, grey) 2D simulation data from Nguyen et al.~\cite{Nguyen2016}.}
}
  \label{fig_incl_wedge2d}
\end{figure}

\begin{figure}[t!]
  \centering
   \subfigure[\REVIEW{Hydrodynamic vertical force}]{
    \includegraphics[scale = 0.3]{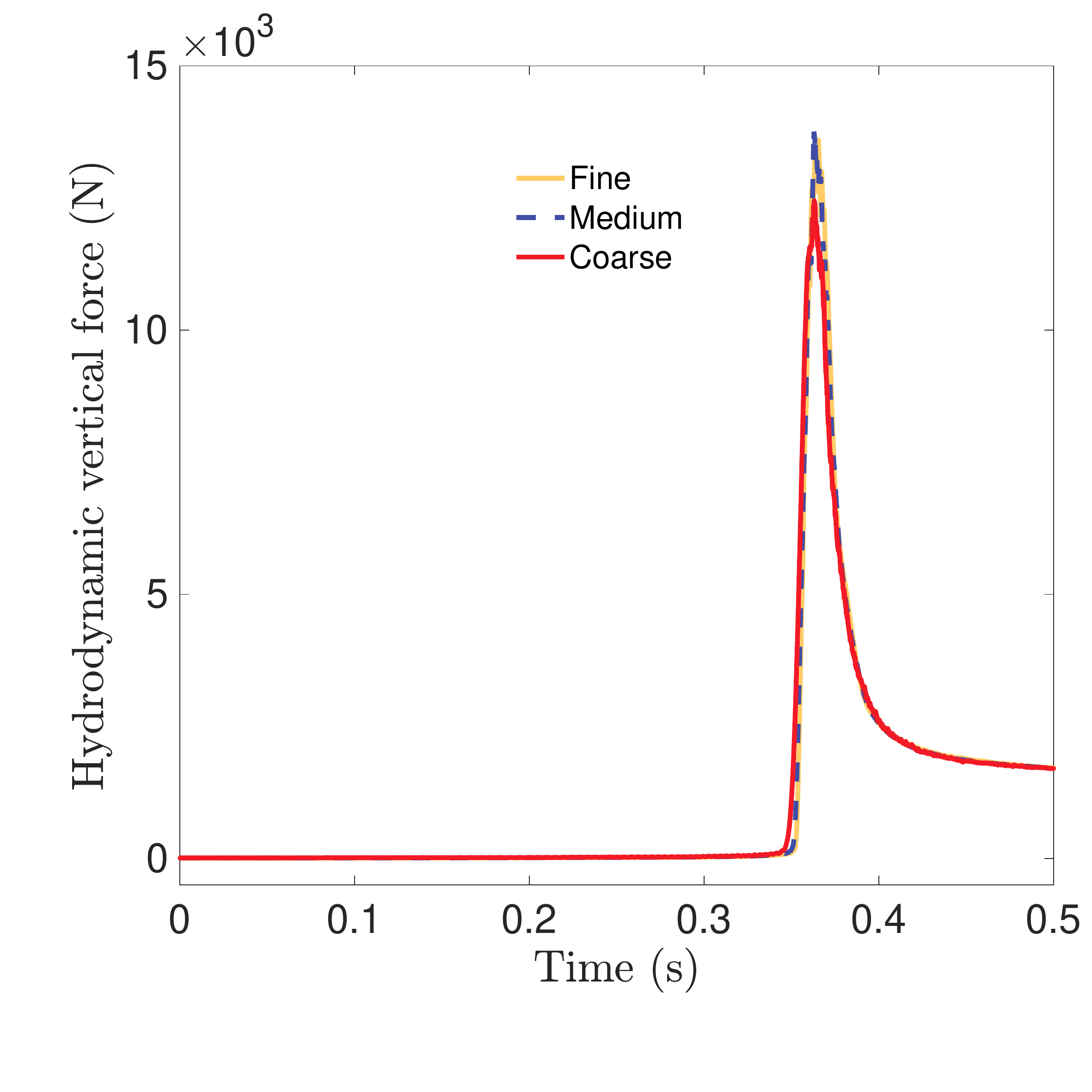}
    \label{incl_wedge2d_force}
  }
  \subfigure[\REVIEW{Hydrodynamic torque}]{
    \includegraphics[scale = 0.3]{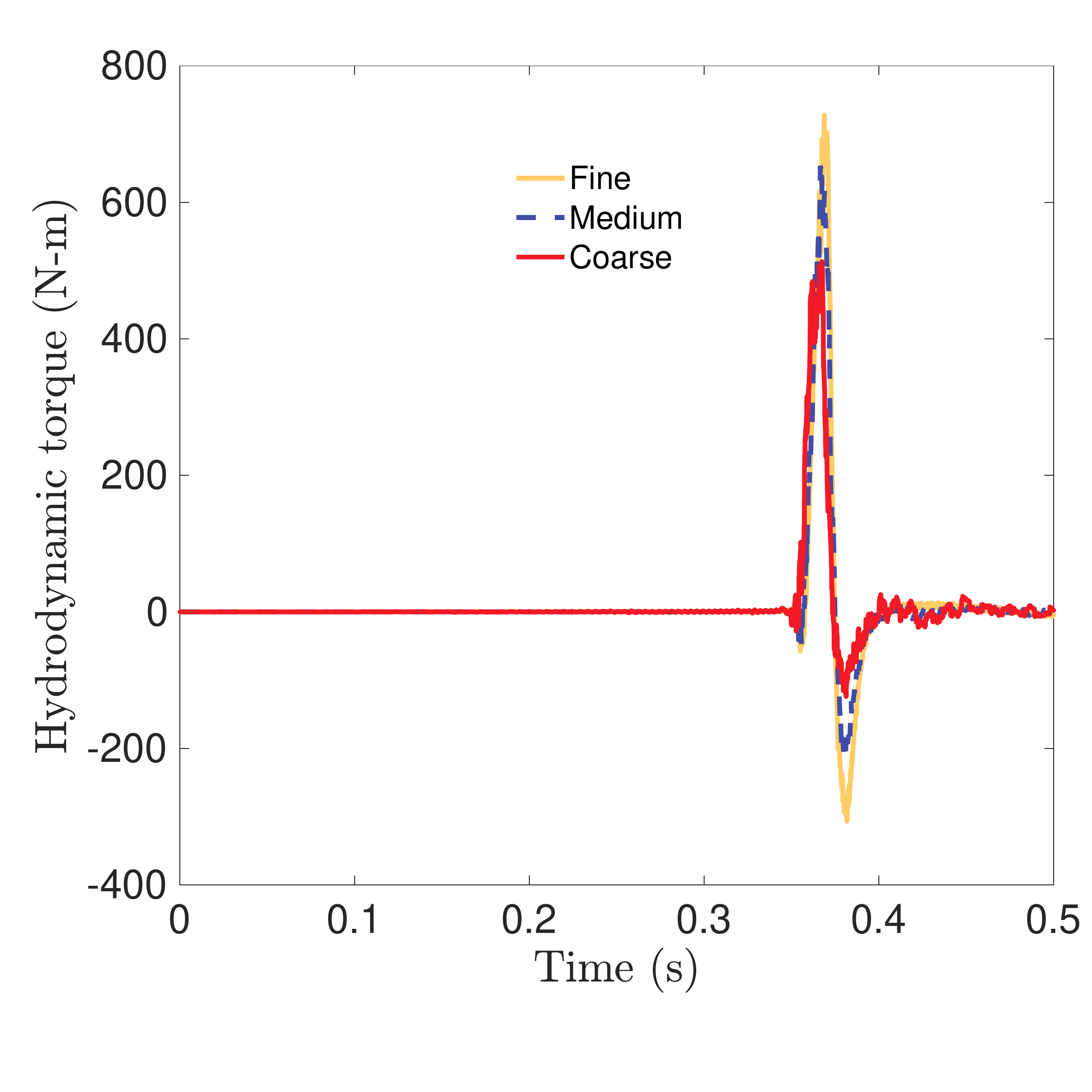}
    \label{incl_wedge2d_torque}
  }
  \caption{\REVIEW{Temporal evolution of hydrodynamic
  \subref{incl_wedge2d_force} vertical force, and 
  \subref{incl_wedge2d_torque} torque for a 2D inclined wedge freely falling into water.
  (---, yellow) Present FD/BP simulation data for a fine grid resolution $814 \times 534$;
  (\texttt{---}, blue) Present FD/BP simulation data for a medium grid resolution $488 \times 320$;
  (---, red) Present FD/BP simulation data for a coarse grid resolution $244 \times 160$.}
}
  \label{fig_incl_wedge2d_moments}
\end{figure}

Fig.~\ref{fig_inclined_wedge2d_viz} shows the evolution of fluid-structure interaction
along with the vorticity generated by the inclined wedge
for the medium and fine grid resolutions. As the wedge falls through the air phase, vortical
structures are shed from the top corners. Upon impact these large scale vortices retain their structure on the medium
grid, while they break down into smaller, satellite vortices on the fine grid. However, the overall trend of the vortex dynamics 
remains the same for the two grid resolutions. Asymmetric splashes emanate from the air-water
interface at later times, as the fine grid resolves the small scale droplets. Note that the fine grid is able to capture emanating spray droplets 
better than the medium grid, which tends to dissipate them. The measured
FSI quantities (heel angle, vertical velocity, and hydrodynamic moments) do not vary significantly between the grids, however.
The overall dynamics are in relatively good agreement with the other numerical results in literature~\cite{Oger2006,Nguyen2016}.

\begin{figure}[]
  \centering
     \subfigure[\REVIEW{Fluid \& structure medium grid}]{
    \includegraphics[scale = 1.0]{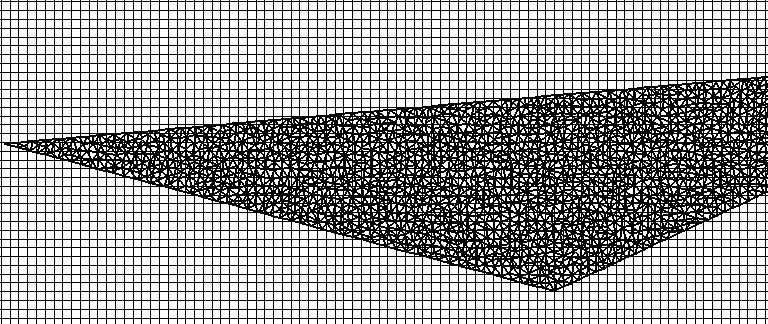}
    \label{iw_med_res}
  }
     \subfigure[\REVIEW{Fluid \& structure fine grid}]{
    \includegraphics[scale = 0.97]{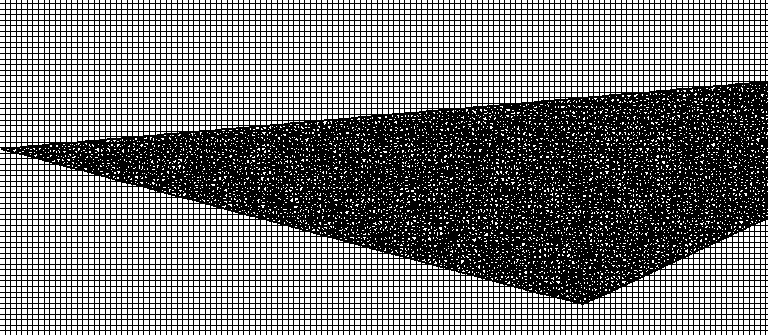}
    \label{iw_high_res}
  } 
   \subfigure[\REVIEW{Medium grid, $t = 0.220$ s}]{
    \includegraphics[scale = 0.3]{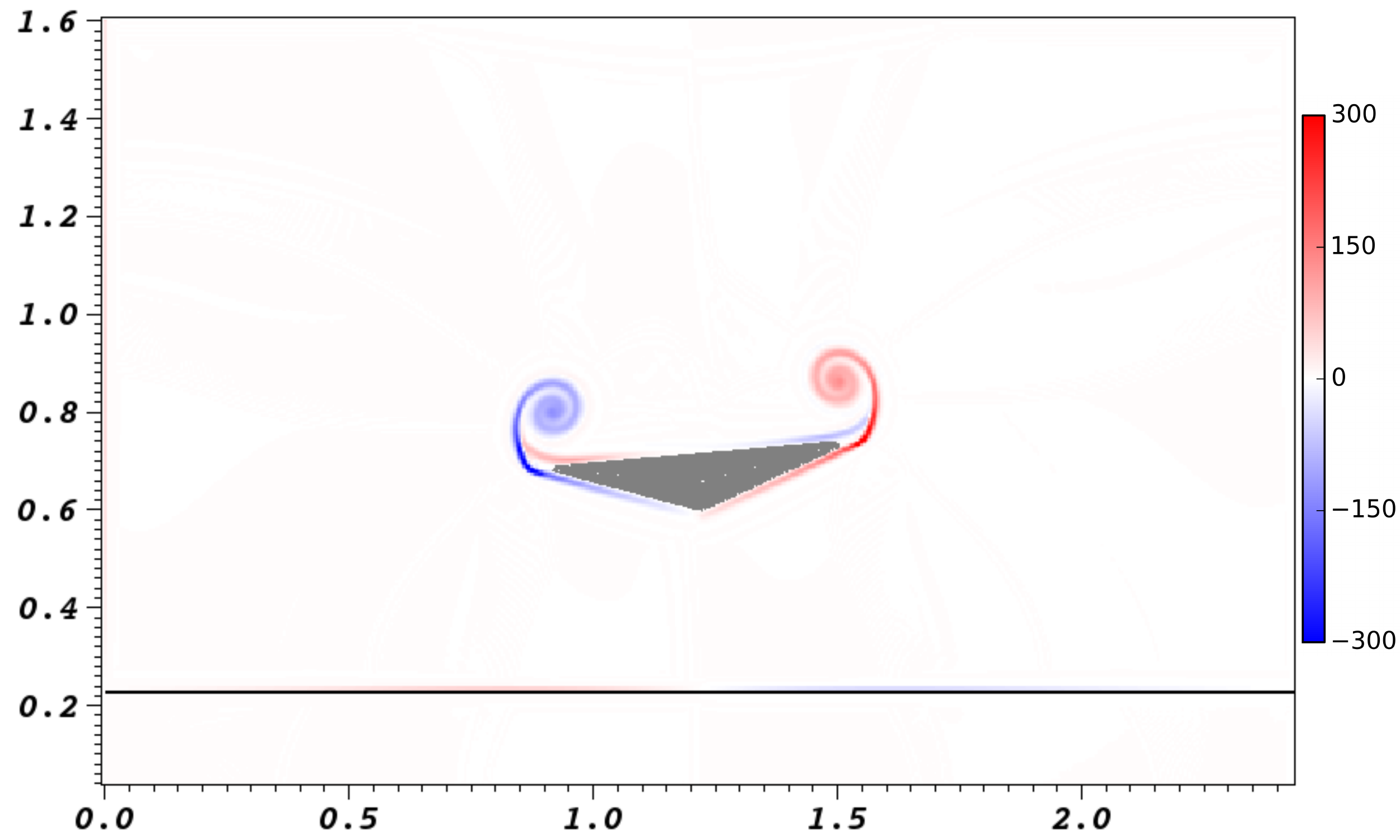}
    \label{iw_med_wedge2d_t0p22}
  }
    \subfigure[\REVIEW{Fine grid, $t = 0.220$ s}]{
    \includegraphics[scale = 0.3]{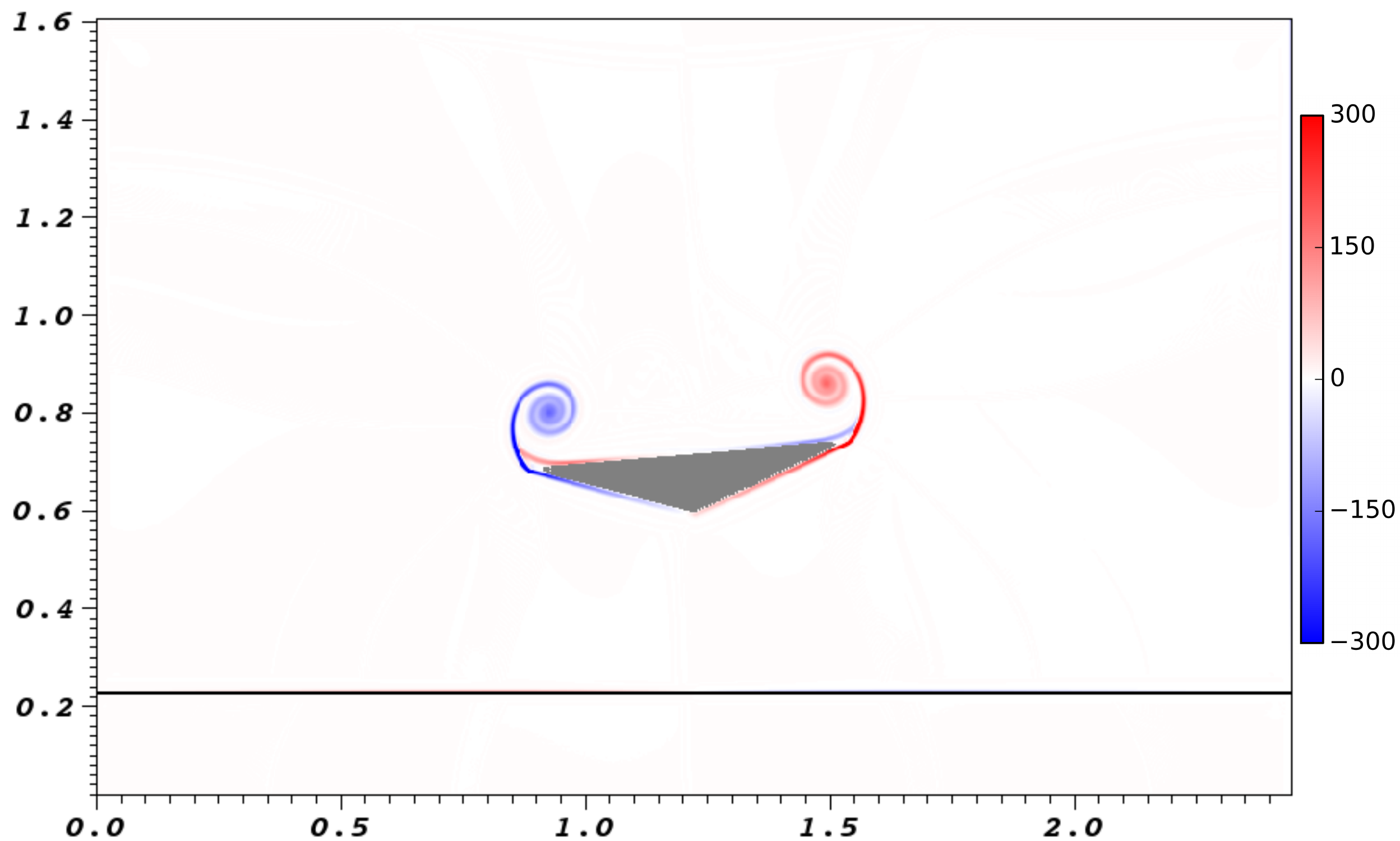}
    \label{iw_hi_wedge2d_t0p22}
  }
    \subfigure[\REVIEW{Medium grid, $t = 0.371$ s}]{
    \includegraphics[scale = 0.3]{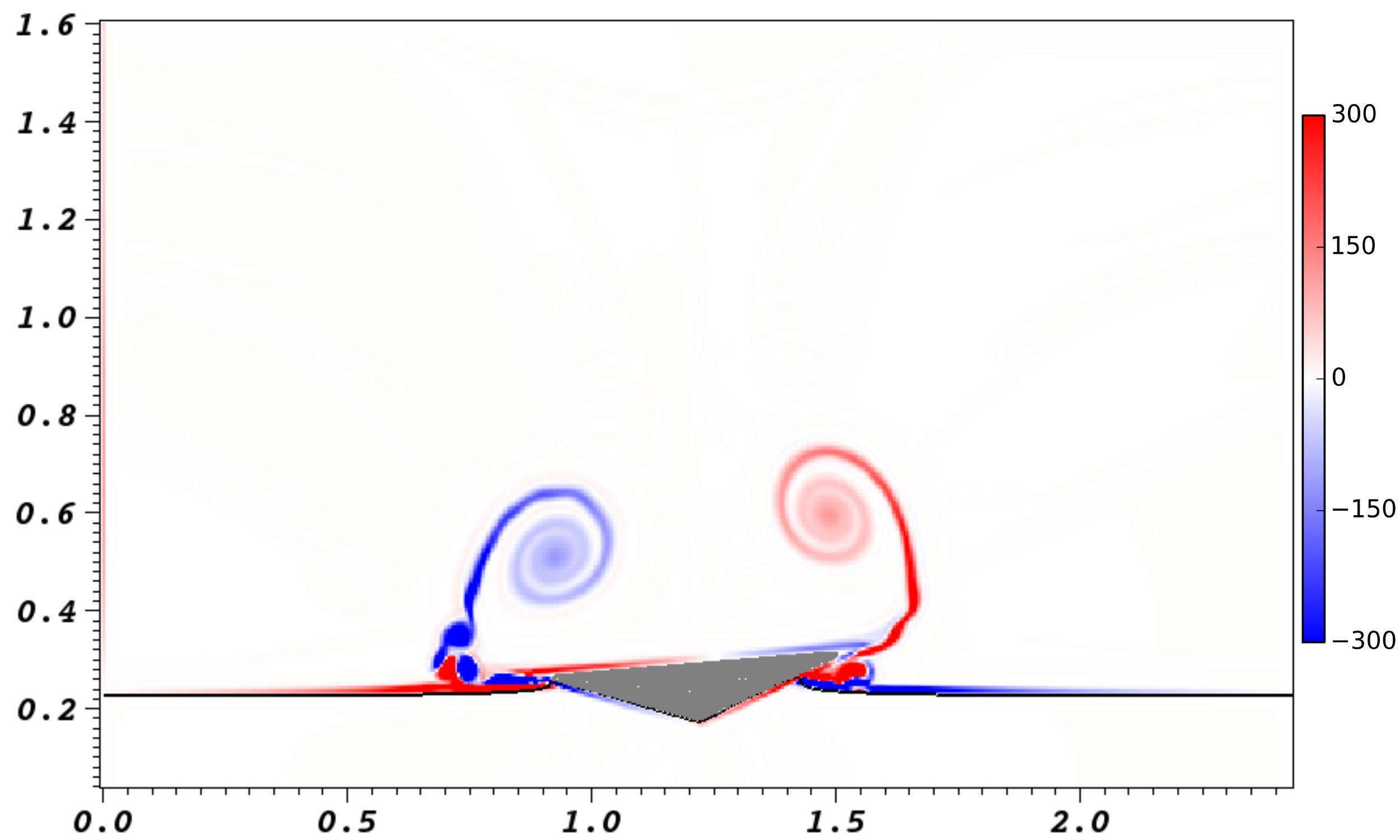}
    \label{iw_med_wedge2d_t0p371}
  }
  \subfigure[\REVIEW{Fine grid, $t = 0.371$ s}]{
    \includegraphics[scale = 0.3]{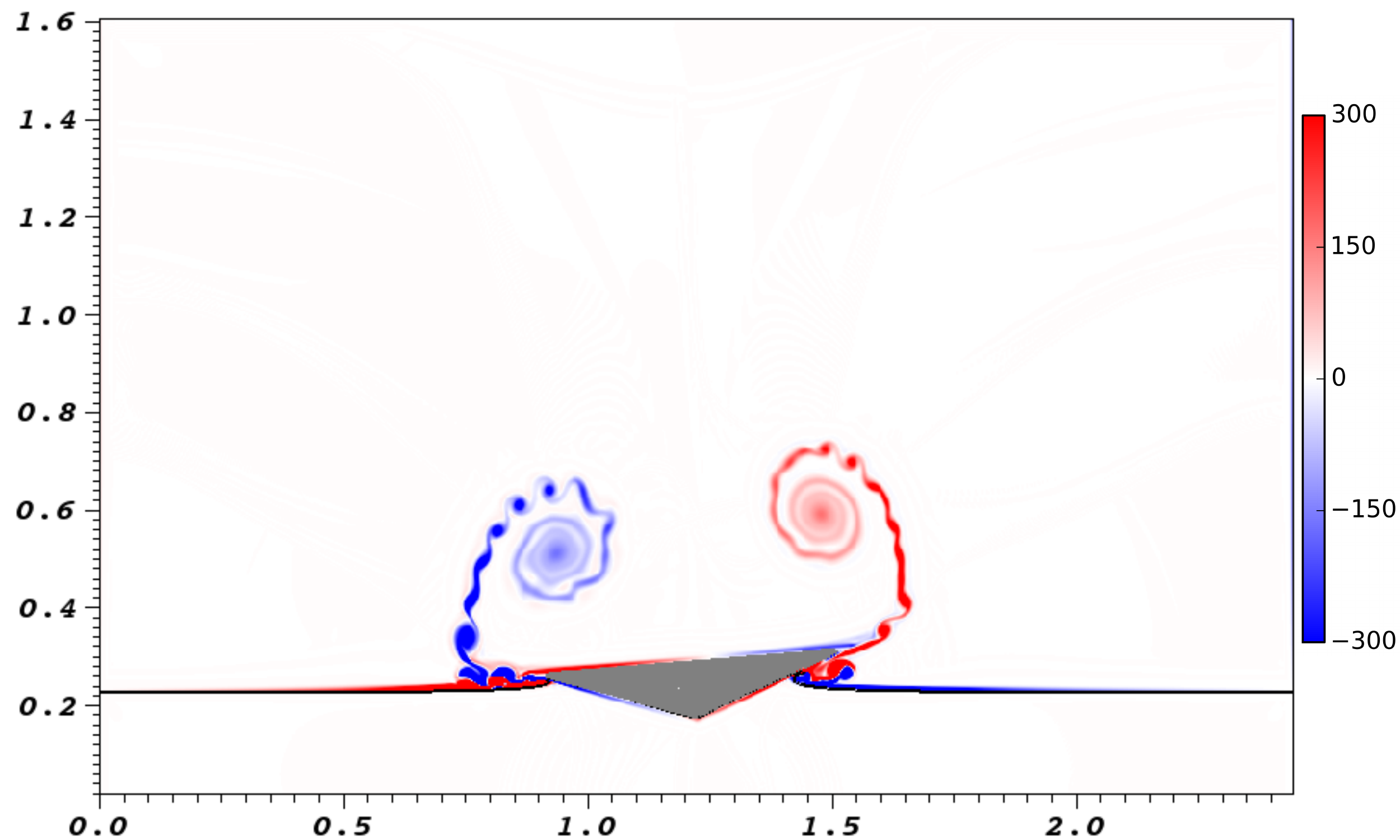}
    \label{iw_hi_wedge2d_t0p371}
  }
   \subfigure[\REVIEW{Medium grid, $t = 0.411$ s}]{
    \includegraphics[scale = 0.3]{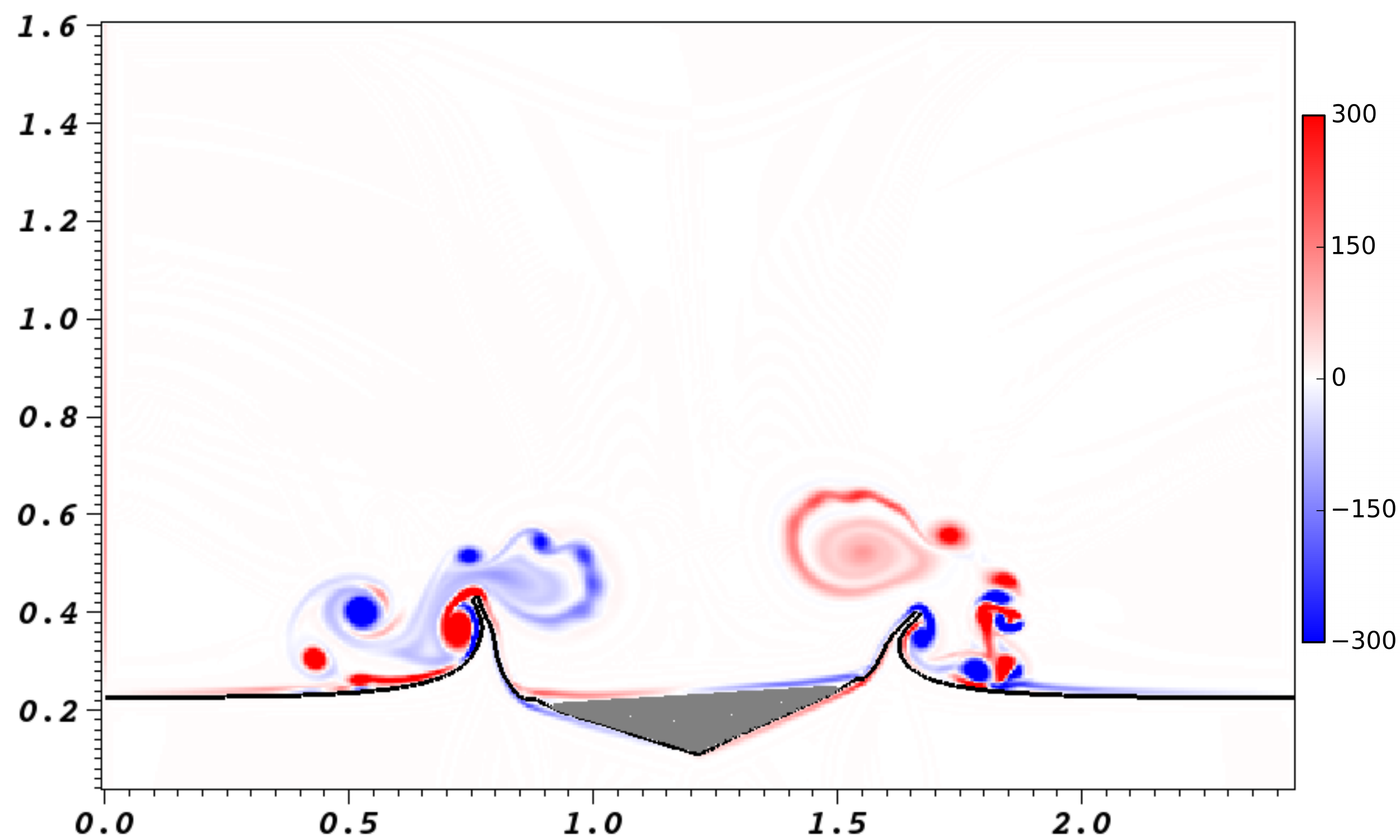}
    \label{iw_med_wedge2d_t0411}
  }
    \subfigure[\REVIEW{Fine grid, $t = 0.411$ s}]{
    \includegraphics[scale = 0.3]{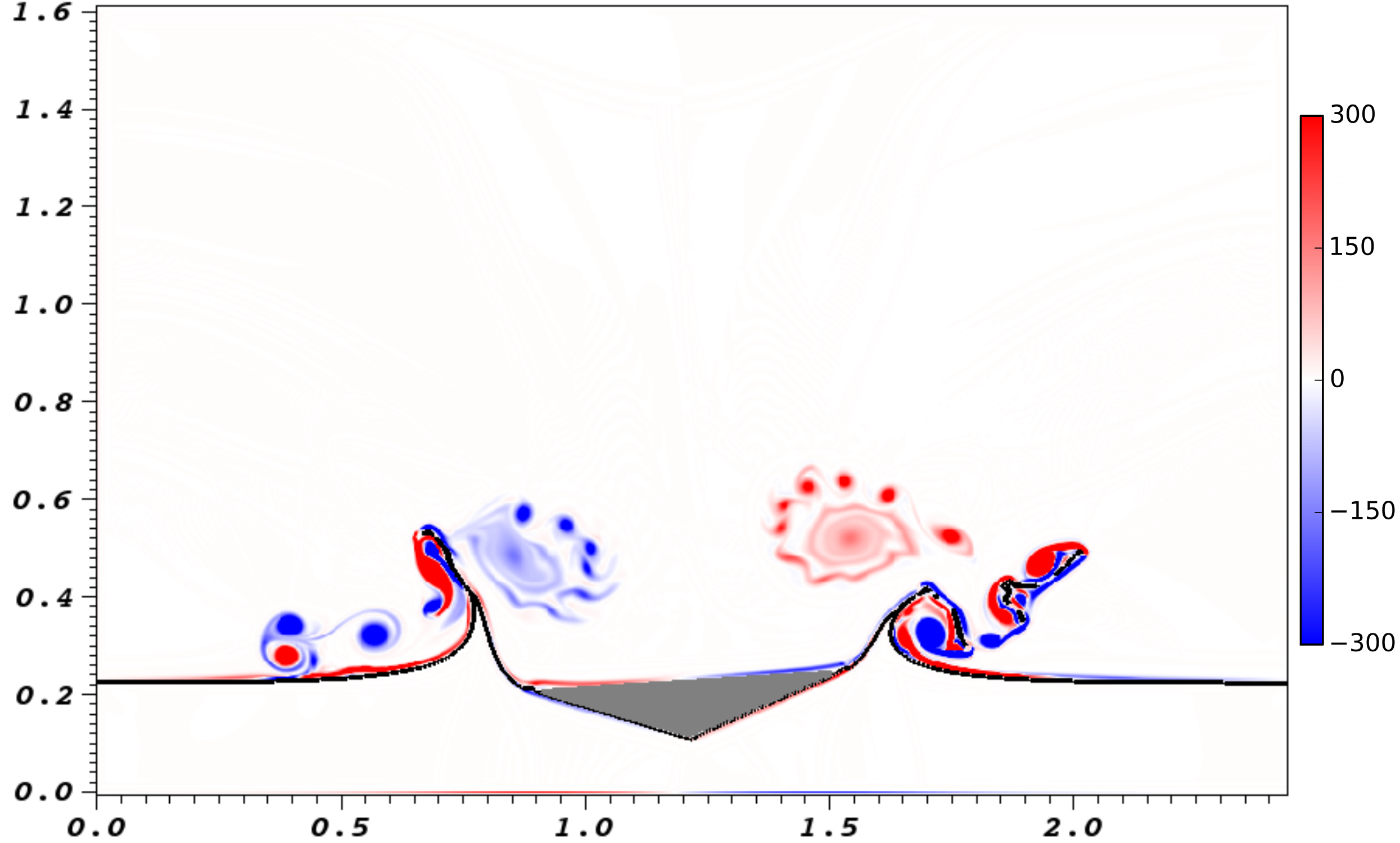}
    \label{iw_hi_wedge2d_t0411}
  }
  \caption{\REVIEW{Vorticity generated by a 2D inclined wedge freely falling into water at three different time 
  instances using the FD/BP method on medium and fine grids. The plotted vorticity is in the range $-300$ to $300$ s$^{-1}$.}
}
  \label{fig_inclined_wedge2d_viz}
\end{figure}

Finally, Fig.~\ref{fig_inclined_wedge2d_pressure} shows the dimensionless pressure field $C_p = p/(\rhol g L)$ as the wedge
slams into the air-water interface. Immediately following impact high pressures are seen at the bottom tip of the wedge. This
high pressure region shifts towards the left side of the body, which has more surface area covered with water. The results
are in excellent agreement with the simulations shown in Nguyen et al.~\cite{Nguyen2016}. With this example, we have 
validated the accuracy the FD/BP method for simulating complex, high inertia FSI. Hereafter, we focus our attention on
comparing and contrasting the FD/BP and FD/IB methodologies.

\begin{figure}[]
  \centering
   \subfigure[\REVIEW{$t = 0.3570$ s}]{
    \includegraphics[scale = 0.3]{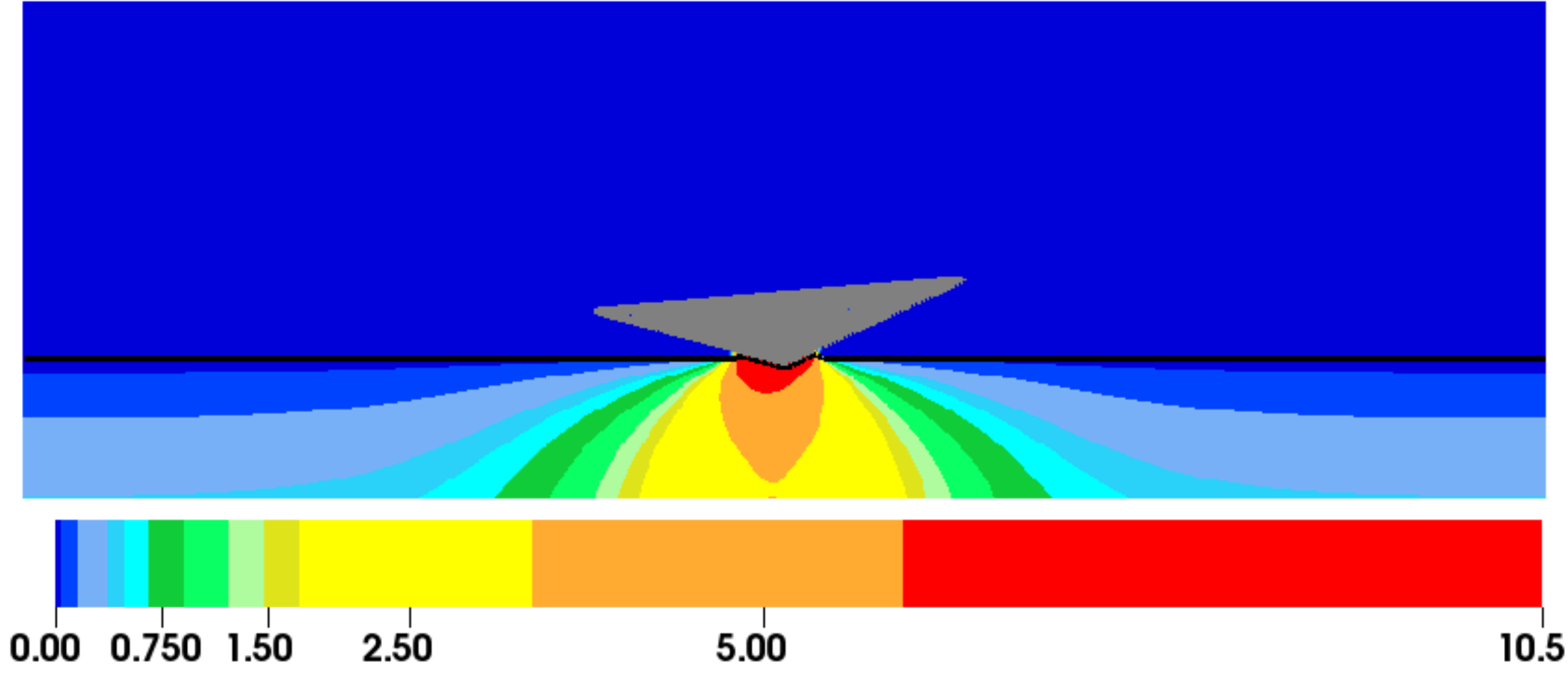}
    \label{iw_pressure_t0p357}
  }
    \subfigure[\REVIEW{$t = 0.3608$ s}]{
    \includegraphics[scale = 0.3]{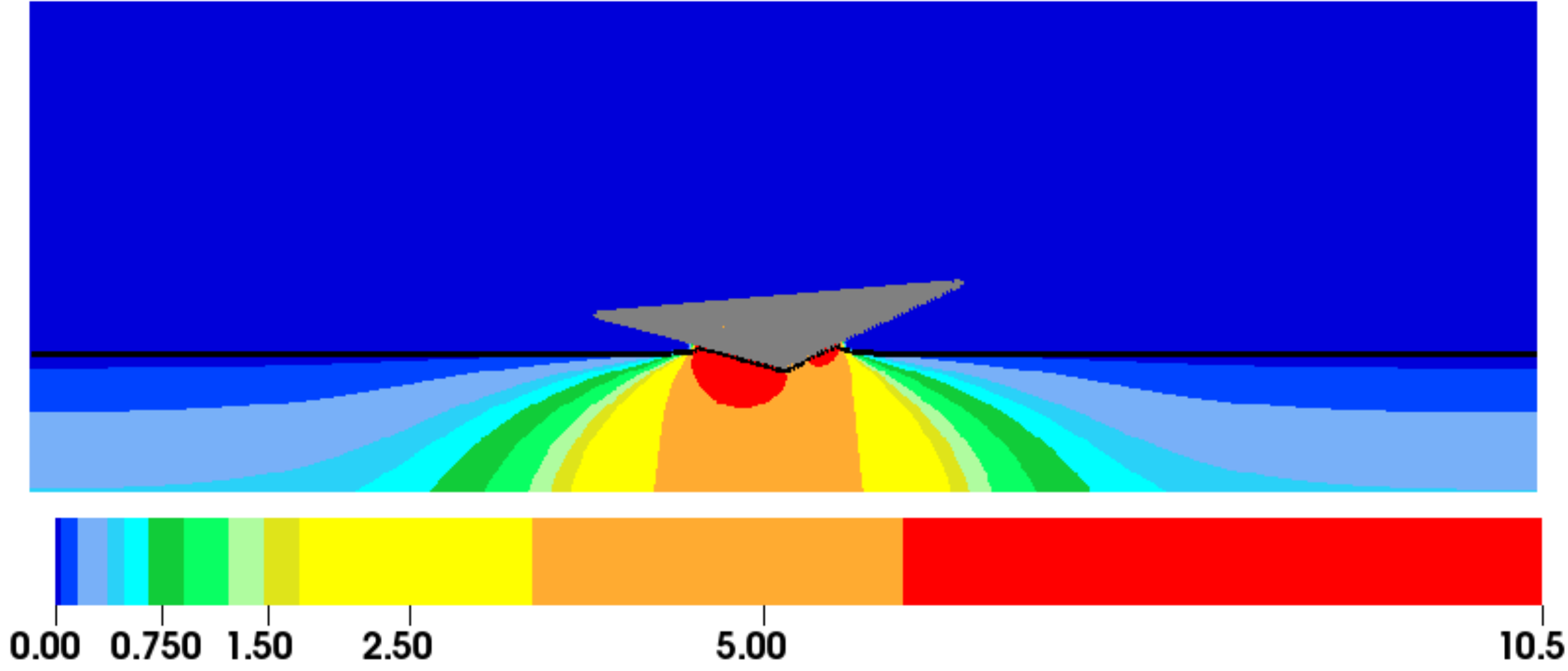}
    \label{iw_pressure_t0p3608}
  }
    \subfigure[\REVIEW{$t = 0.3676$ s}]{
    \includegraphics[scale = 0.3]{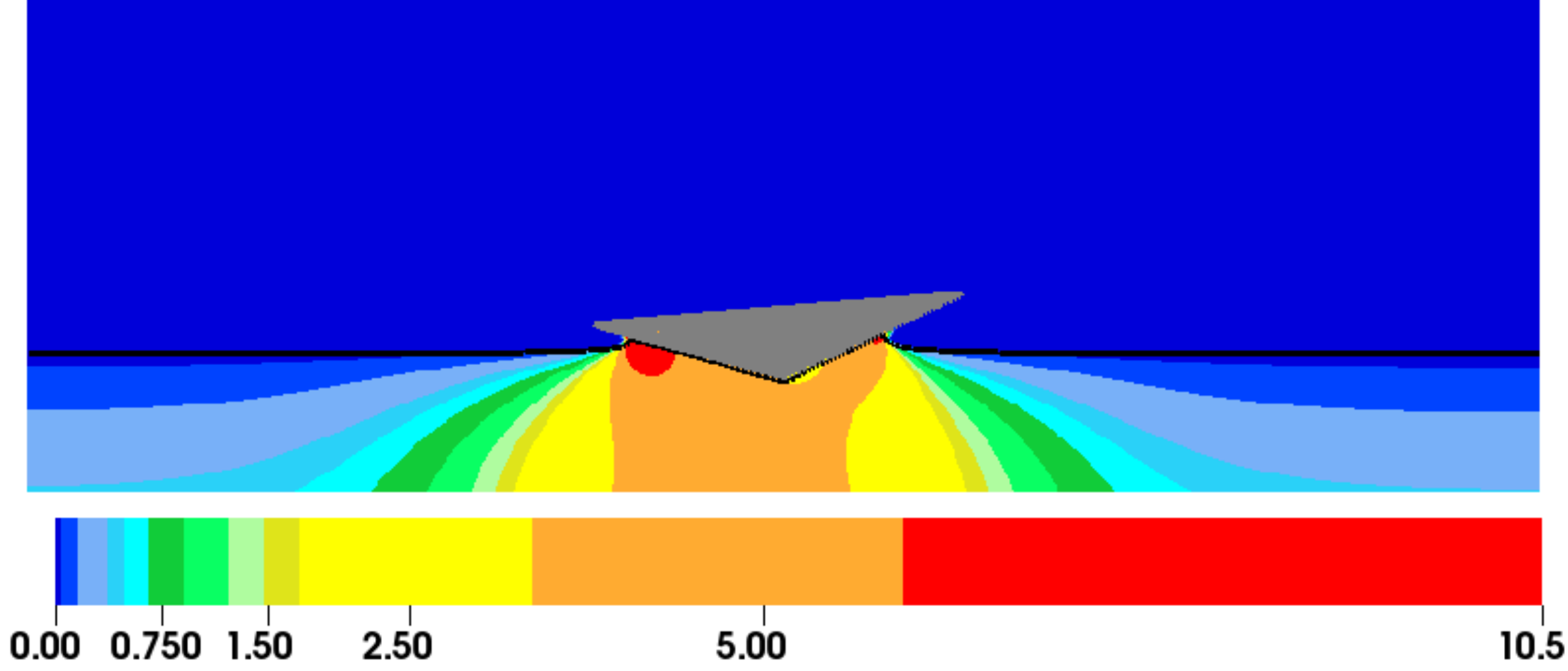}
    \label{iw_pressure_t0p3676}
  }
   \subfigure[\REVIEW{$t = 0.3726$ s}]{
    \includegraphics[scale = 0.3]{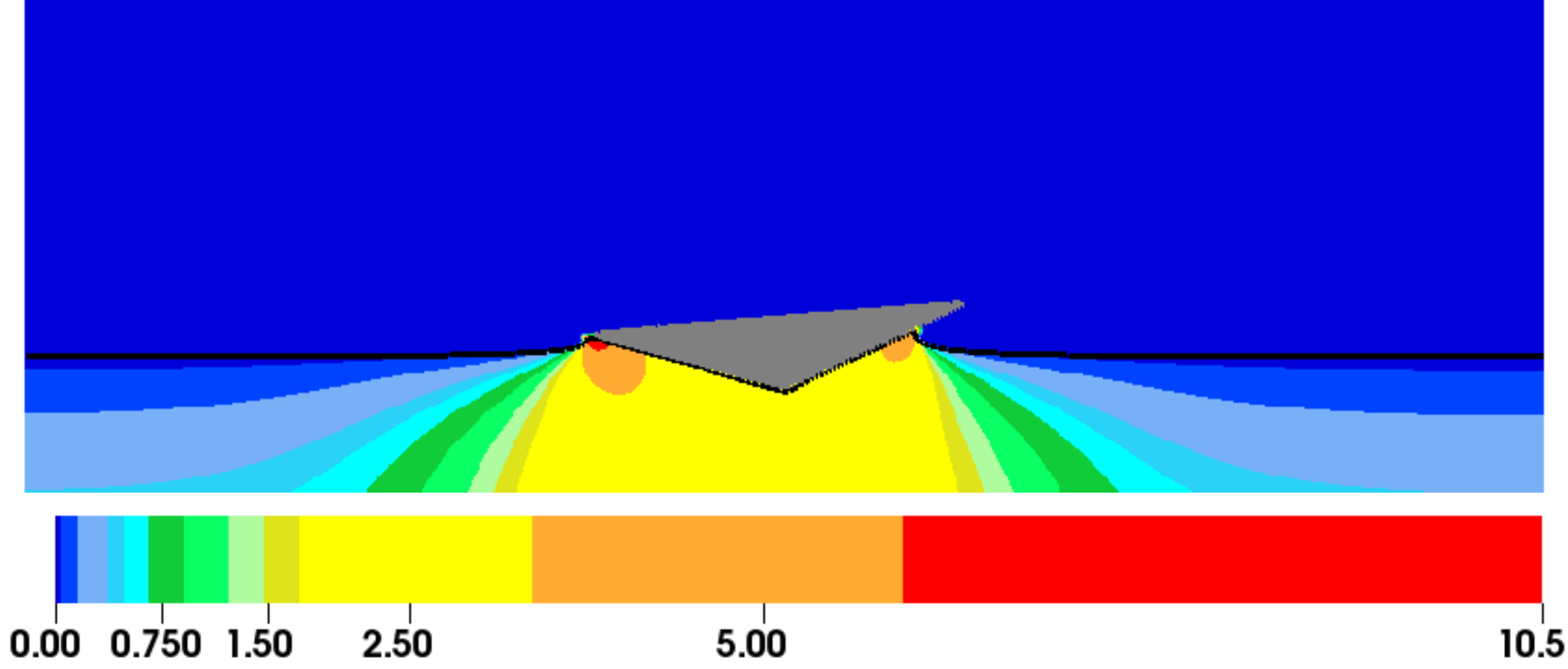}
    \label{iw_pressure_t0p3726}
  }
    \subfigure[\REVIEW{$t = 0.3826$ s}]{
    \includegraphics[scale = 0.3]{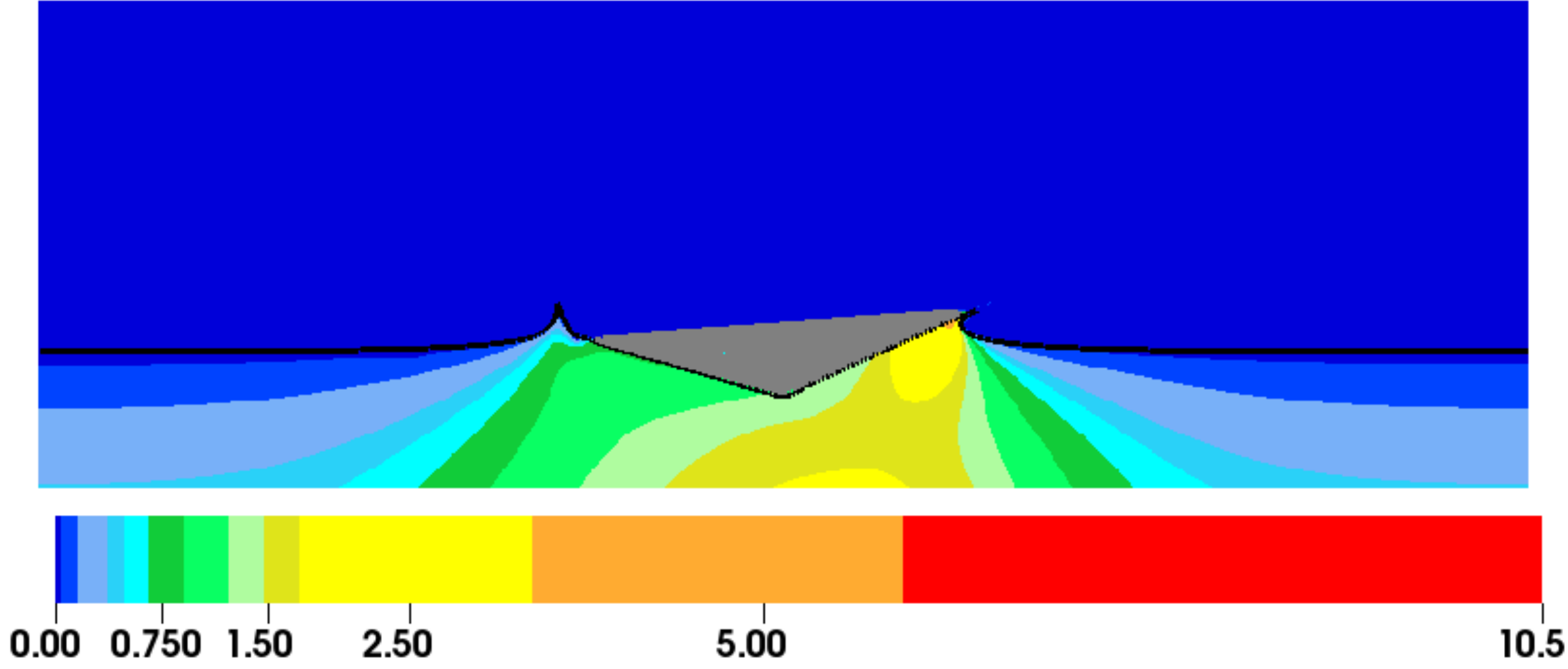}
    \label{iw_pressure_t0p3826}
  }
   \subfigure[\REVIEW{$t = 0.3950$ s}]{
    \includegraphics[scale = 0.3]{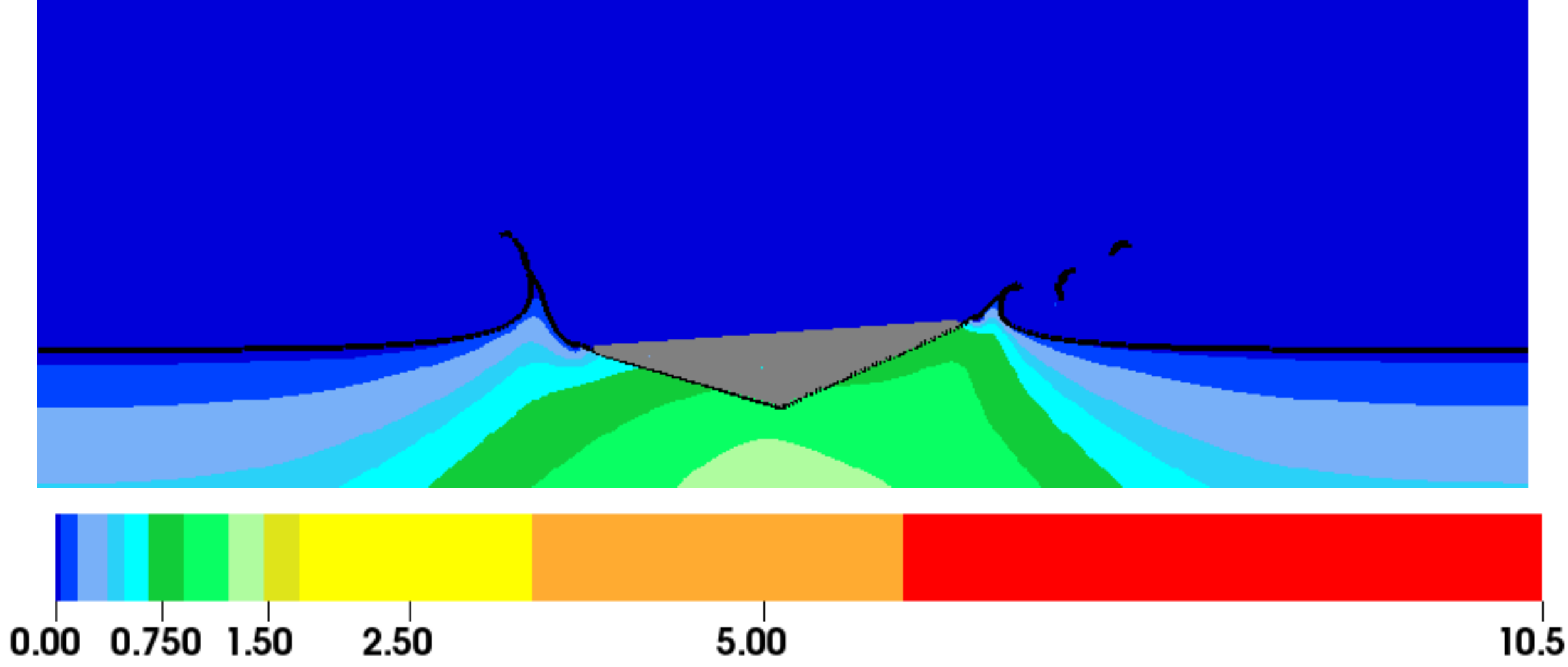}
    \label{iw_pressure_t0p395}
  }
    \subfigure[\REVIEW{$t = 0.4500$ s}]{
    \includegraphics[scale = 0.3]{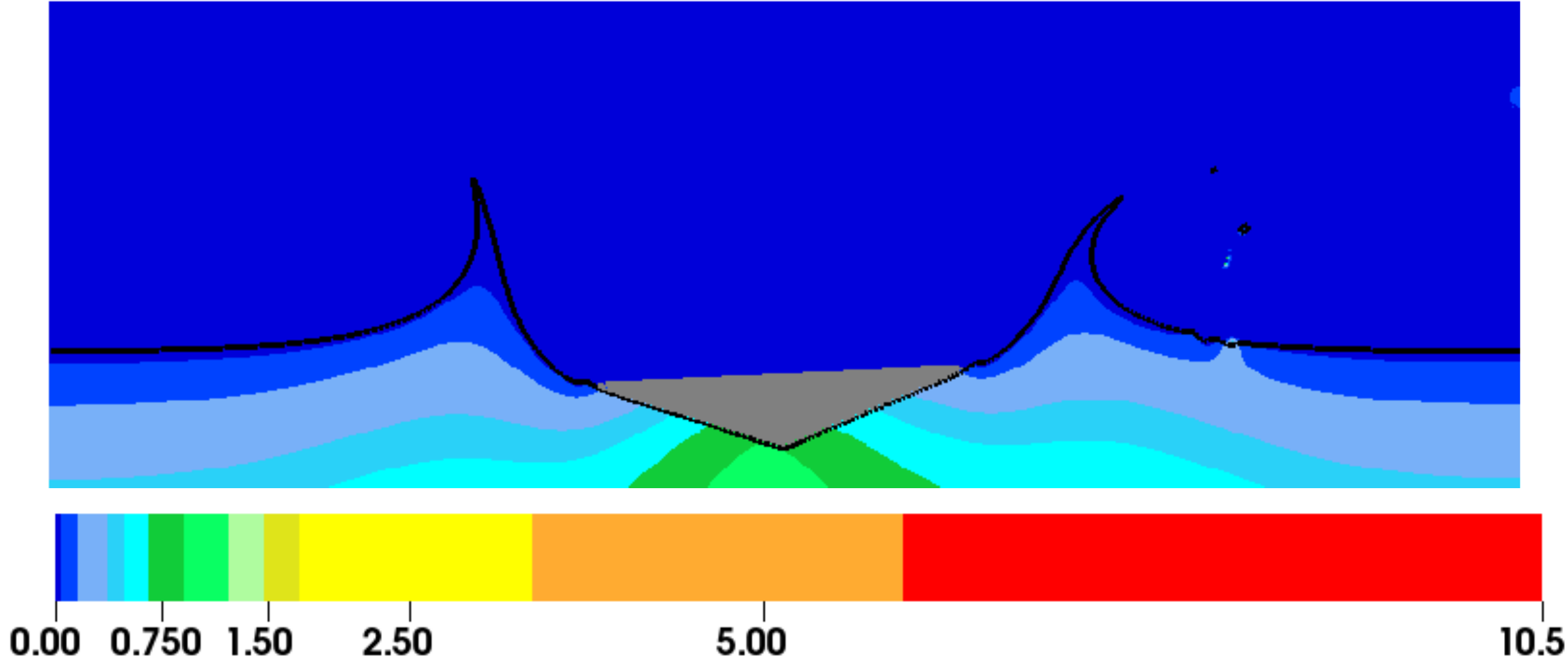}
    \label{iw_pressure_t0p45}
  }
  \caption{\REVIEW{Dimensionless pressure field $C_p = p/(\rhol g L)$ around a 2D inclined wedge freely falling into water at seven different time 
  instances using the FD/BP method on fine grid.}
}
  \label{fig_inclined_wedge2d_pressure}
\end{figure}
}

\subsection{Water-entry/exit of a free falling wedge}

\REVIEW{Next, we} consider the problem of a wedge-shaped object impacting
a pool of water. A 2D triangular body with top length $L = 1.2$ m is placed within
a computational domain of size $\Omega = [0, 10L] \times [0, 2.5L]$. The wedge
is oriented with one of its vertices pointing downwards, making a $25^\circ$ deadrise
angle with the horizontal. Water occupies the bottom third of the domain, while air occupies
the remainder of the tank. The bottom point of the wedge is placed with initial position
$(X_0, Y_0) = (5L, 23L/12)$ and the wedge has a density of $\rhos = 466.6$ kg/m$^3$. 
The free fall height of wedge is $\Delta s = 13L/12$. 

\begin{figure}[t!]
  \centering
  
   \subfigure[\REVIEW{Vertical position}]{
    \includegraphics[scale = 0.3]{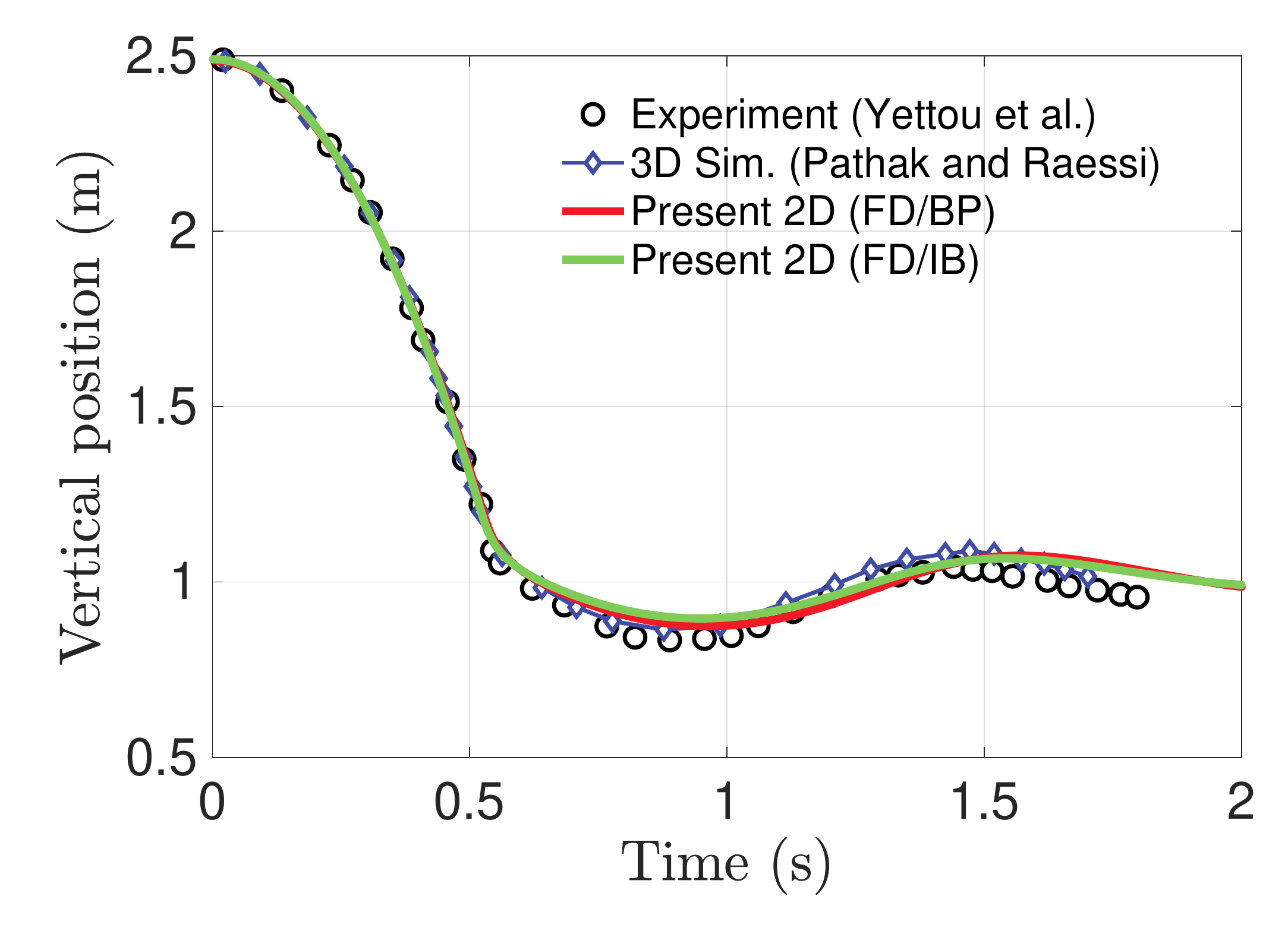}
    \label{wedge2d_posn}
  }
     \subfigure[\REVIEW{Vertical velocity}]{
    \includegraphics[scale = 0.3]{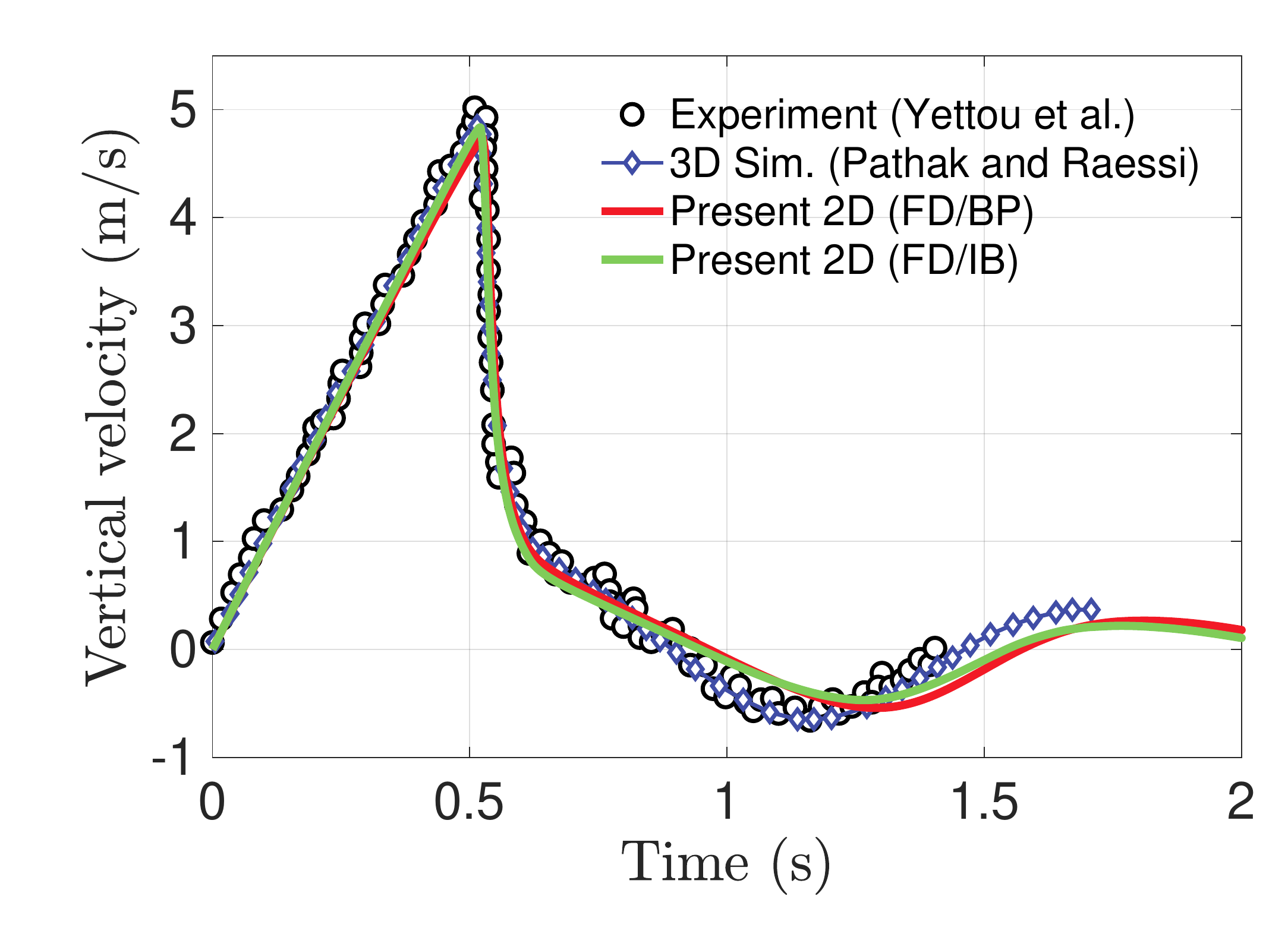}
    \label{wedge2d_vel}
  }
   \subfigure[\REVIEW{Vertical force}]{
    \includegraphics[scale = 0.3]{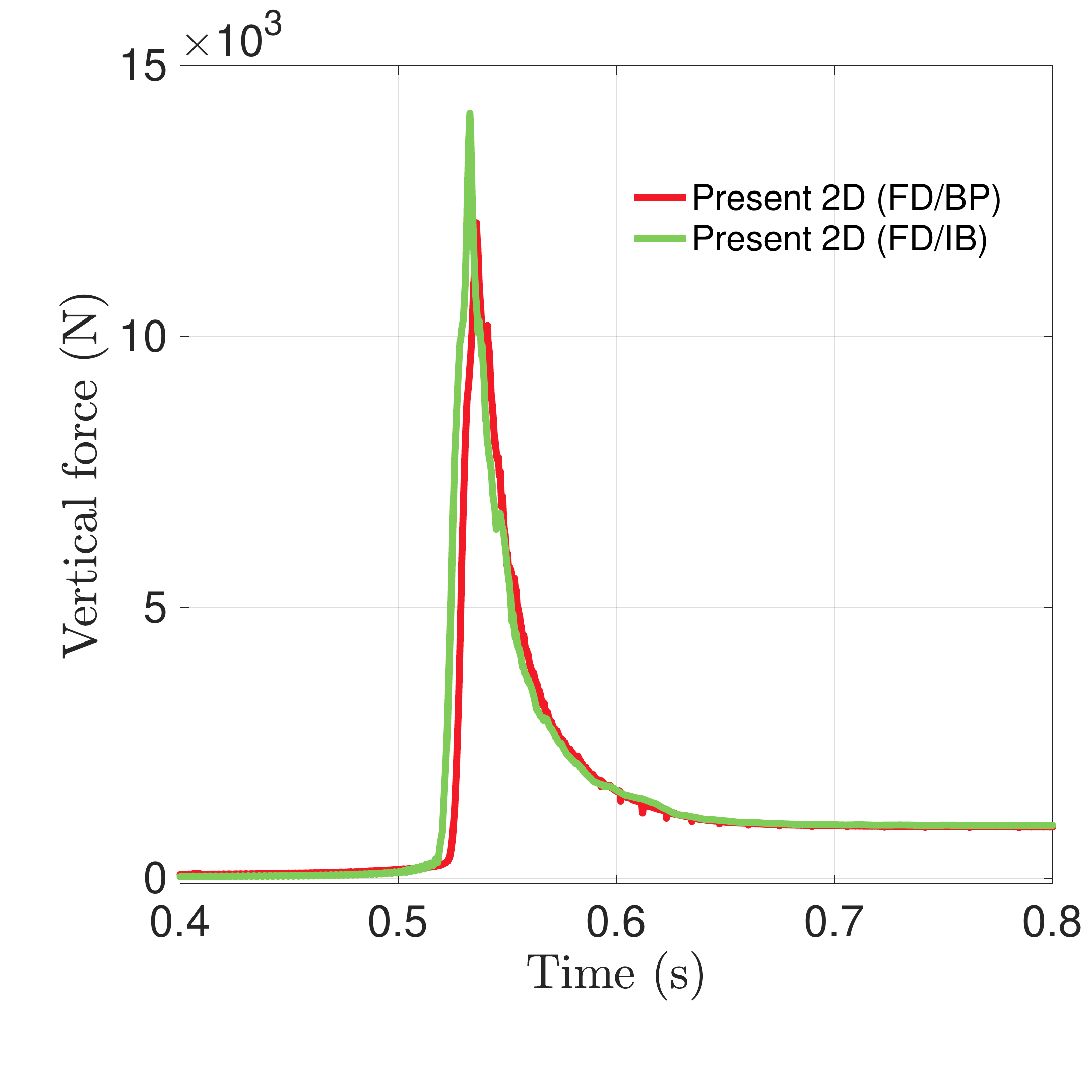}
    \label{wedge2d_force}
  }
  \caption{\REVIEW{Temporal evolution of \subref{wedge2d_posn} vertical position, \subref{wedge2d_vel} vertical velocity, and \subref{wedge2d_force} 
  vertical force on a 2D wedge freely falling in water. ($\circ$, black) experimental data from Yettou et al.~\cite{Yettou2006};
  ($\diamond$, blue) 3D simulation data from Pathak and Raessi~\cite{Pathak16}; (---, red) present FD/BP simulation data;
  (---, green) present FD/IB simulation data.}
}
  \label{fig_wedge2d_com}
\end{figure}

The 2D domain is discretized by a $1200 \times 300$ uniform grid, \REVIEW{which corresponds to 120 grid cells per wedge length. 
This grid resolution was found sufficient to resolve the FSI dynamics of a free-falling wedge in the previous section.} 
A constant time step size of $\dt = 6.25 \times 10^{-5}$ s is used. Fig.~\ref{fig_wedge2d_com} shows the time evolution 
of center of mass vertical position and velocity. The results are in good agreement with prior numerical~\cite{Pathak16} 
and experimental studies~\cite{Yettou2006}. The wedge reaches a peak velocity in the air phase just before impacting the water surface.
The vertical velocity keeps descending as it penetrates further into water. Eventually the buoyancy forces reverse
the wedge's velocity and it begins to exit the pool. Additionally, the hydrodynamic forces on the wedge in the vertical direction are compared. \REVIEW{As seen in Fig.~\ref{wedge2d_force}, the FD/IB method produces smooth forces compared to the FD/BP method. This is because evaluation of hydrodynamic forces for the FD/IB method is done in an \emph{extrinsic} manner (Eqs.~\eqref{eq_lmforce}-\eqref{eq_lmtorque}). 
In contrast, the FD/BP method computes the forces in an \emph{intrinsic} manner through direct stress evaluation involving 
derivatives of the velocity field. This is known (see Bergmann and Iollo~\cite{Bergmann11}) to produce spurious oscillations in the force evaluation. 
In our prior work (Nangia et al.~\cite{Nangia17}) we proposed a \emph{moving control volume approach} to obtain 
smooth forces by converting intrinsic integrals to extrinsic integrals over a moving Cartesian box. 
Others (Verma et al.~\cite{Verma2017} and Patel et al.~\cite{Namu2018}) have proposed to evaluate stress derivatives on a ``lifted" surface 
that is two cells distance away from the original surface to avoid small scale oscillations for the FD/BP method.}
   
\begin{figure}[]
  \centering
  
   \subfigure[FD/BP, $t = 0.45$ s]{
    \includegraphics[scale = 0.32]{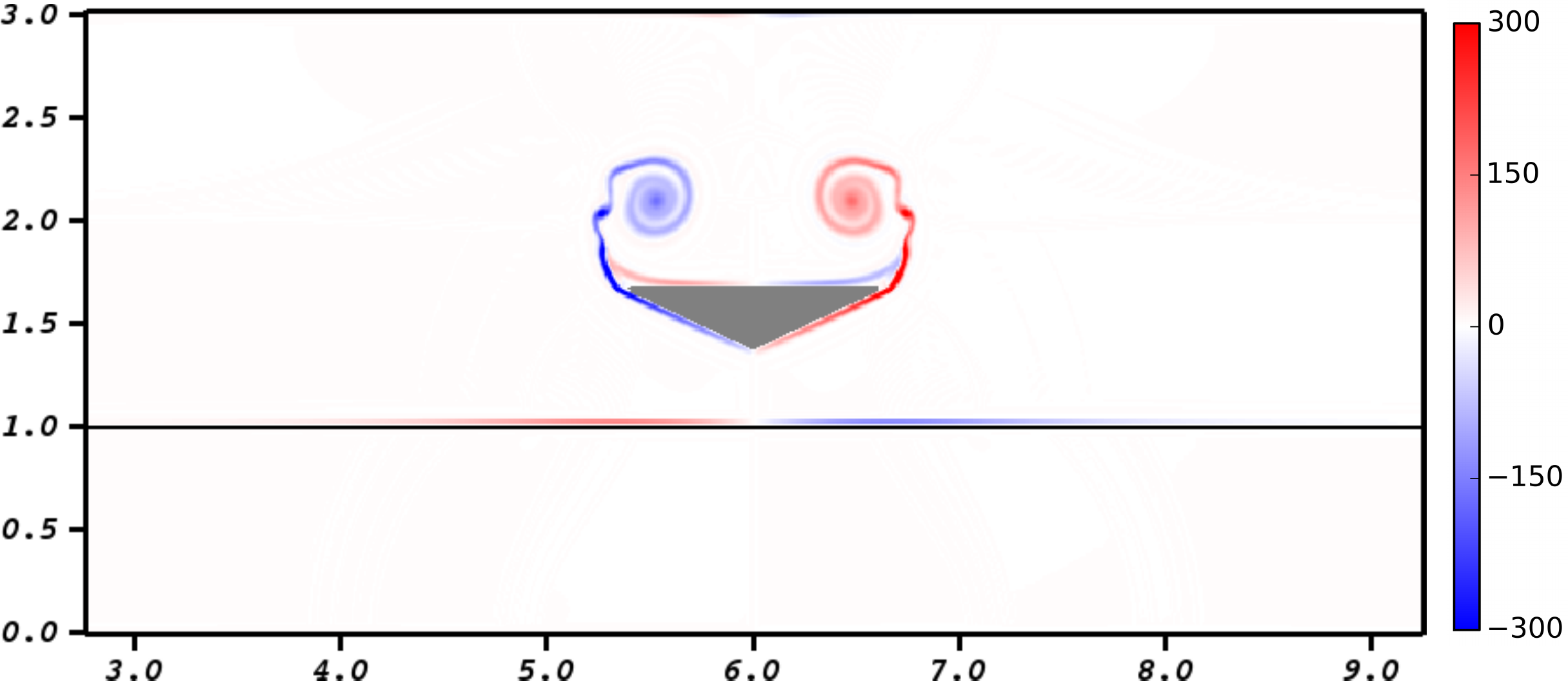}
    \label{ibls_wedge2d_t0p45}
  }
     \subfigure[FD/IB, $t = 0.45$ s]{
    \includegraphics[scale = 0.32]{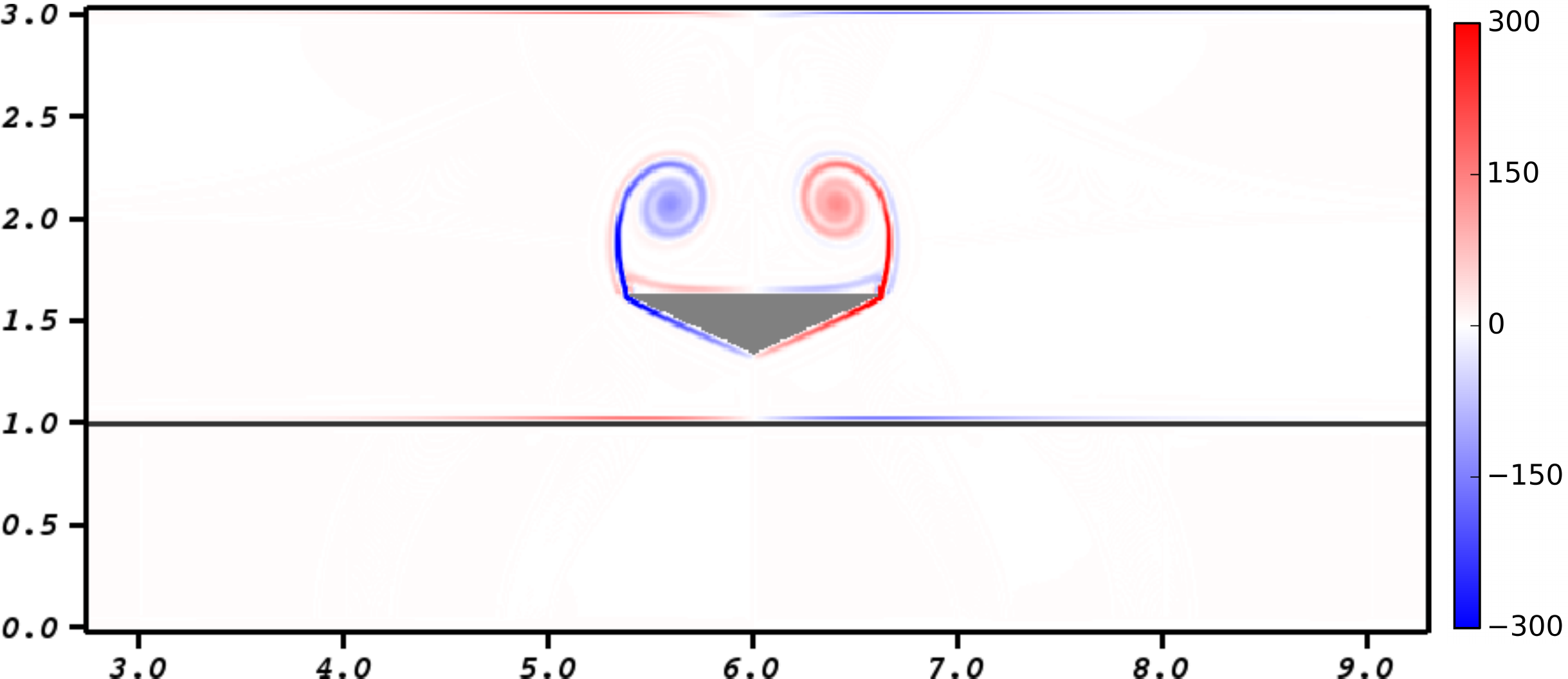}
    \label{cib_wedge2d_t0p45}
  }
   \subfigure[FD/BP, $t = 0.5625$ s]{
    \includegraphics[scale = 0.32]{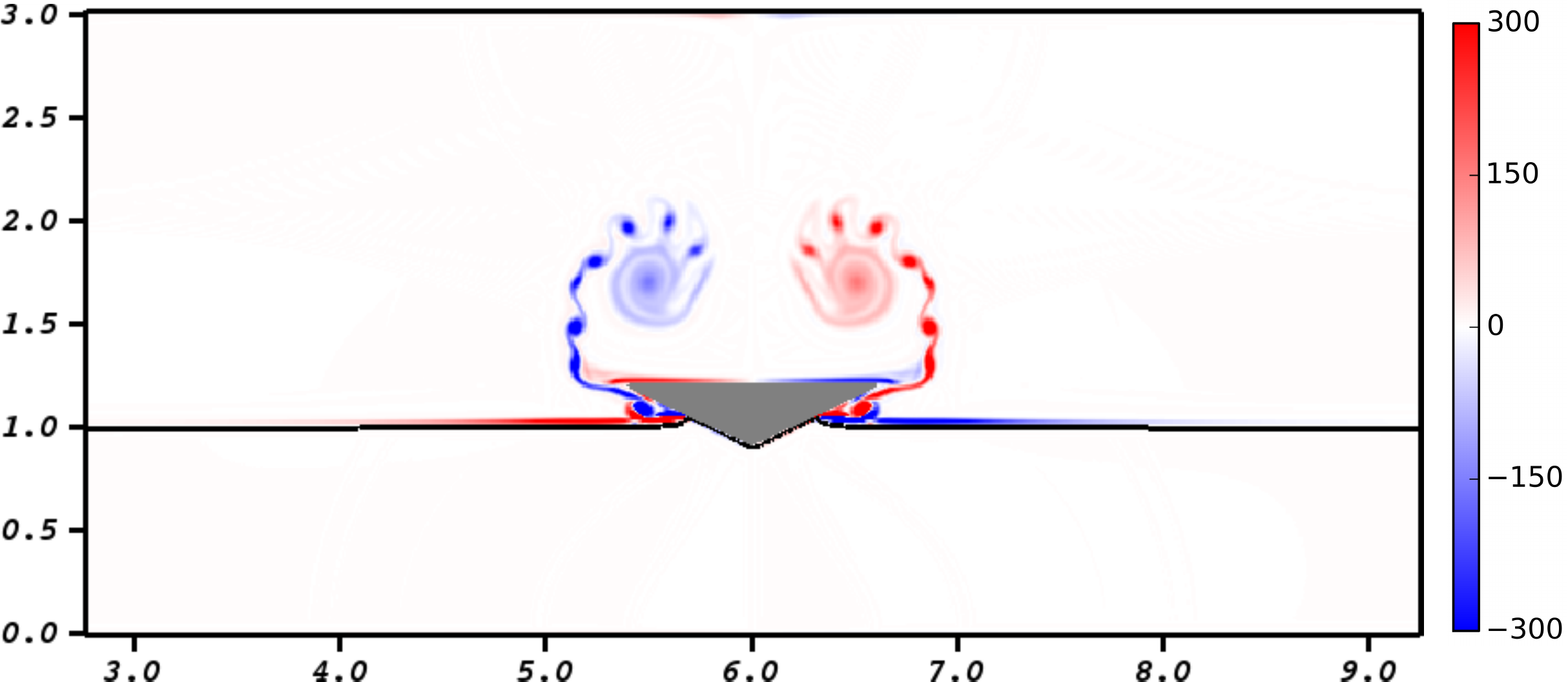}
    \label{ibls_wedge2d_t0p5625}
  }
    \subfigure[FD/IB, $t = 0.5625$ s]{
    \includegraphics[scale = 0.32]{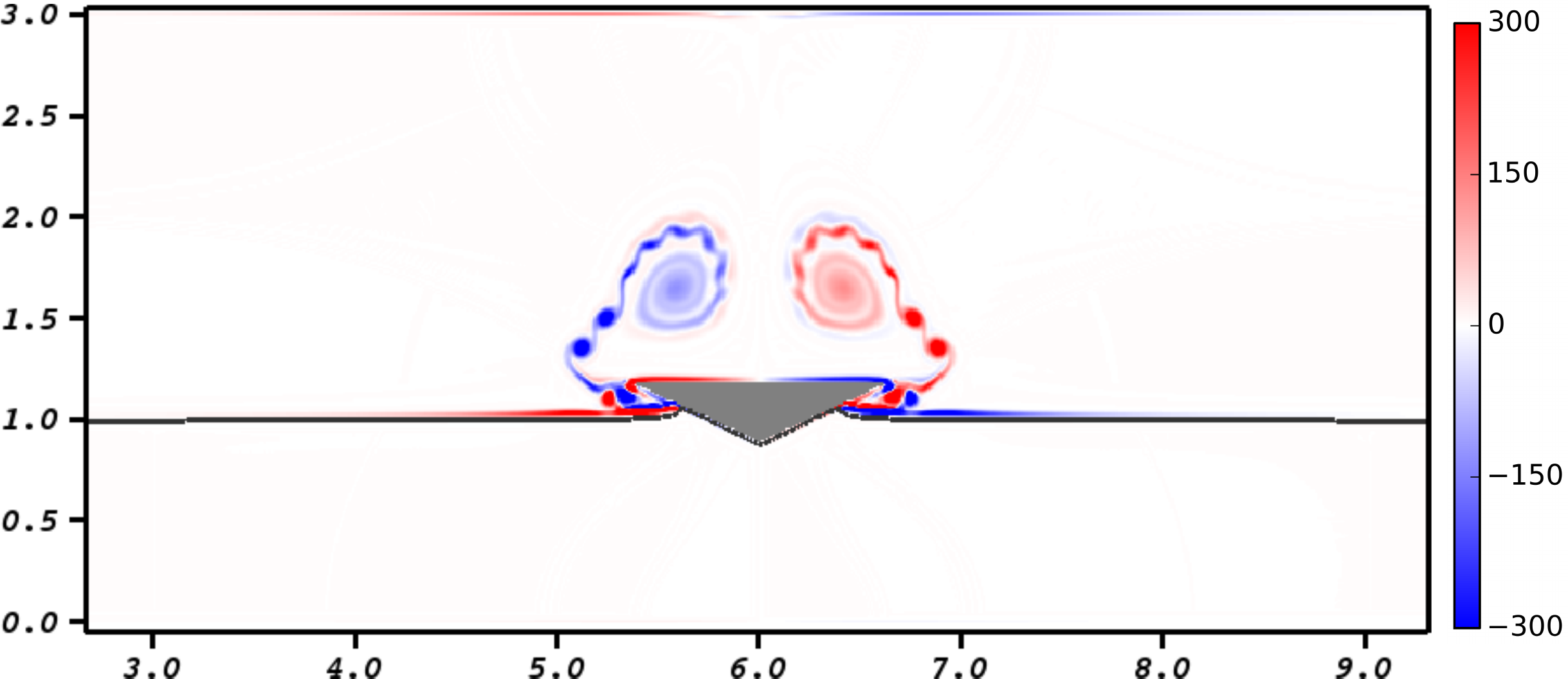}
    \label{cib_wedge2d_t0p5625}
  }
  \subfigure[FD/BP, $t = 0.875$ s]{
    \includegraphics[scale = 0.32]{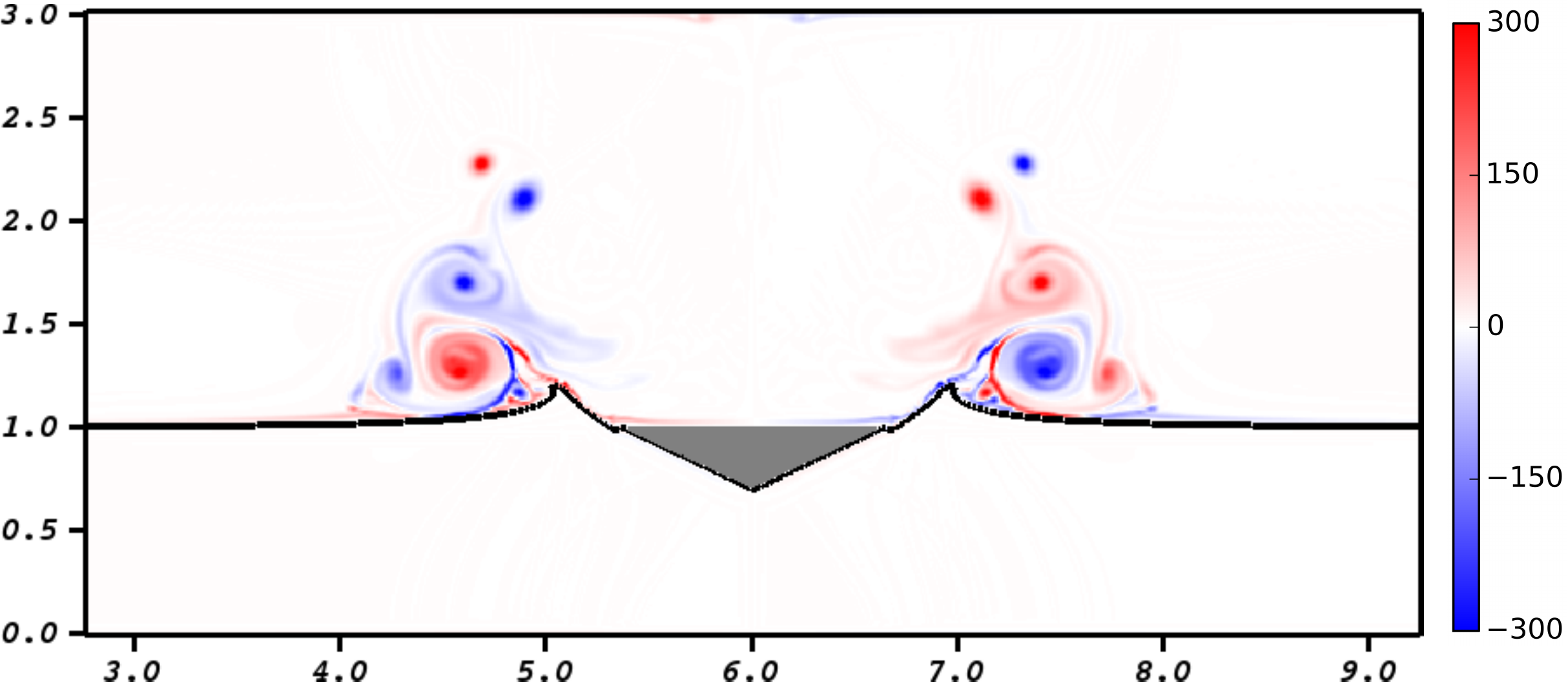}
    \label{ibls_wedge2d_t0p875}
  }
    \subfigure[FD/IB, $t = 0.875$ s]{
    \includegraphics[scale = 0.32]{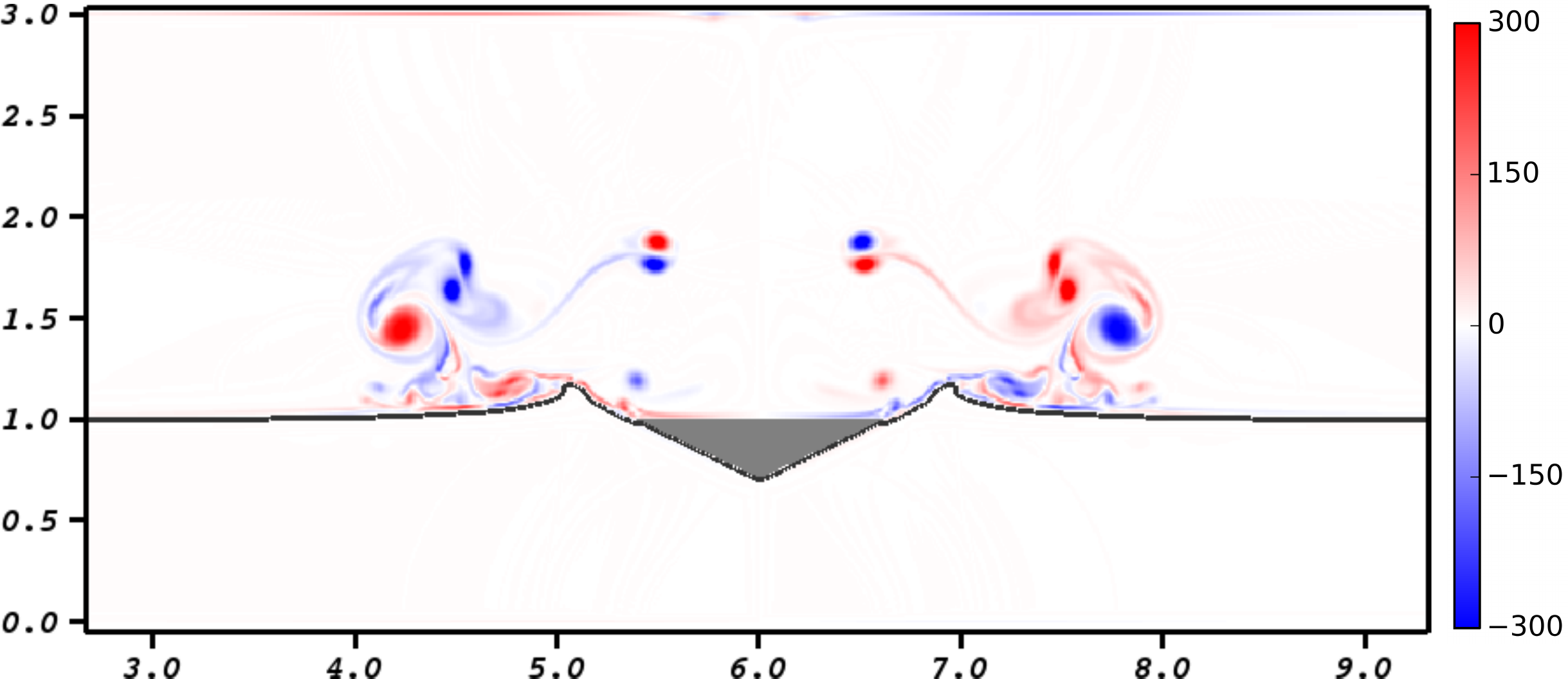}
    \label{cib_wedge2d_t0p875}
  }
    \subfigure[FD/BP, $t = 1.25$ s]{
    \includegraphics[scale = 0.32]{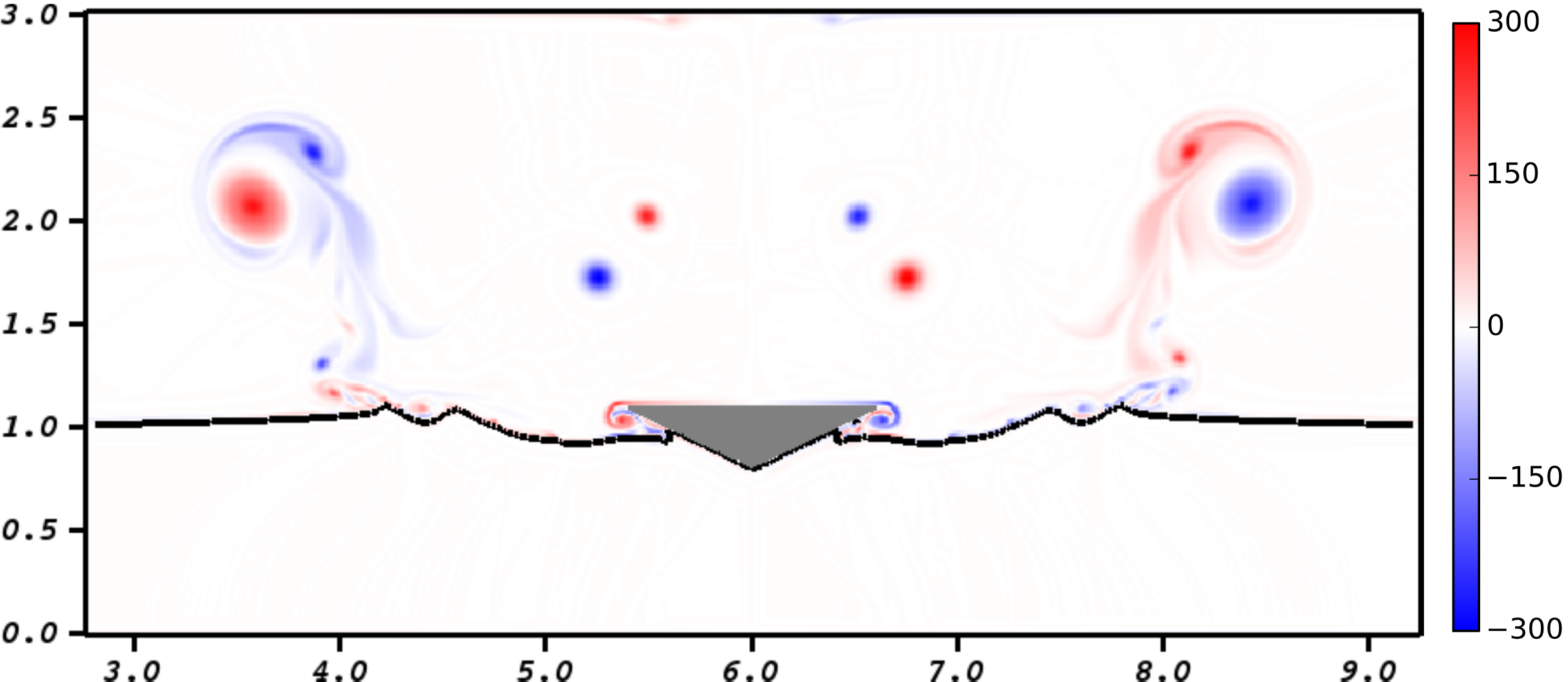}
    \label{ibls_wedge2d_1pt25}
  }
  \subfigure[FD/IB, $t = 1.25$ s]{
    \includegraphics[scale = 0.32]{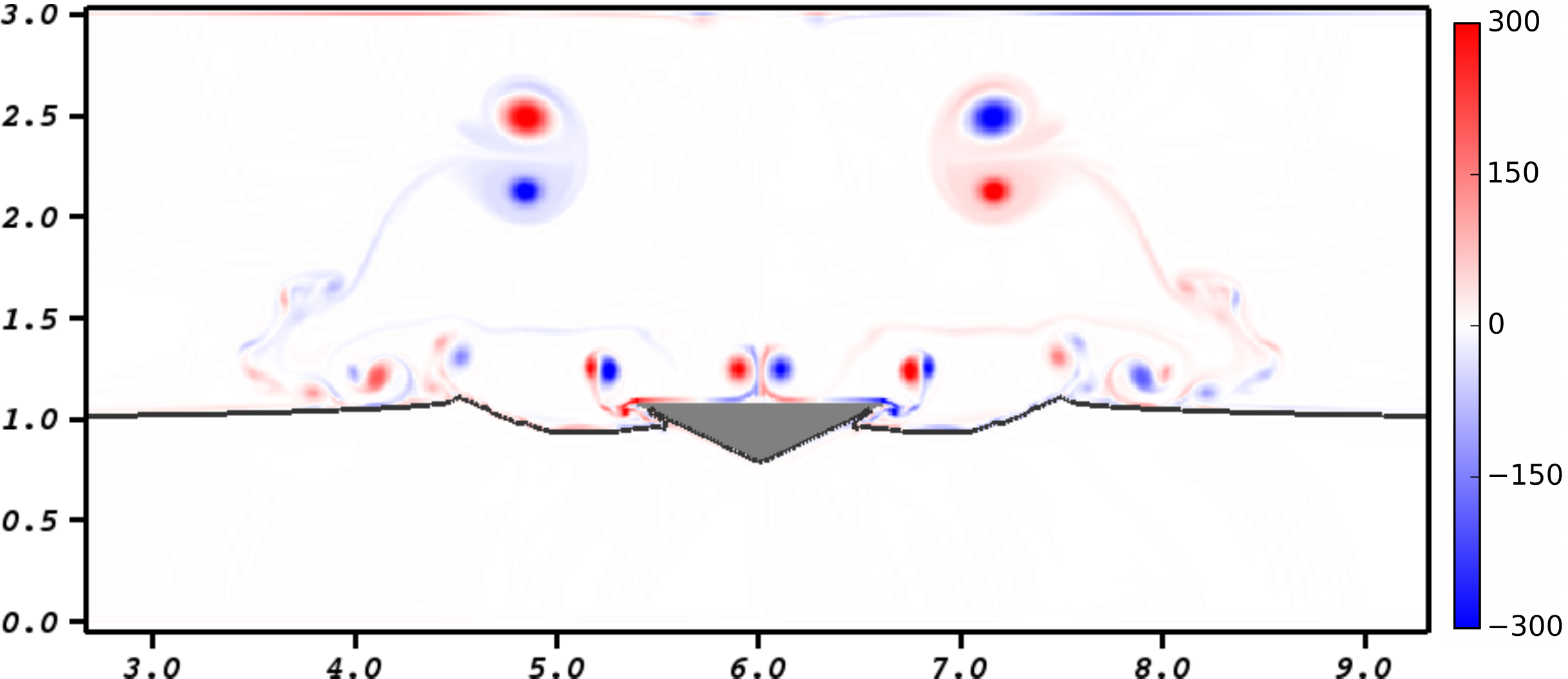}
    \label{cib_wedge2d_1pt25}
  }
  \caption{Vorticity generated by a free falling 2D wedge at four different time 
  instances using the FD/BP and FD/IB methods. The plotted vorticity is in the range $-300$ to $300$ s$^{-1}$.
}
  \label{fig_wedge2d_viz}
\end{figure}

\begin{figure}[]
  \centering
  
   \subfigure[\REVIEW{Experiment~\cite{Greenhow1983}}]{
    \includegraphics[width = 0.33\textwidth]{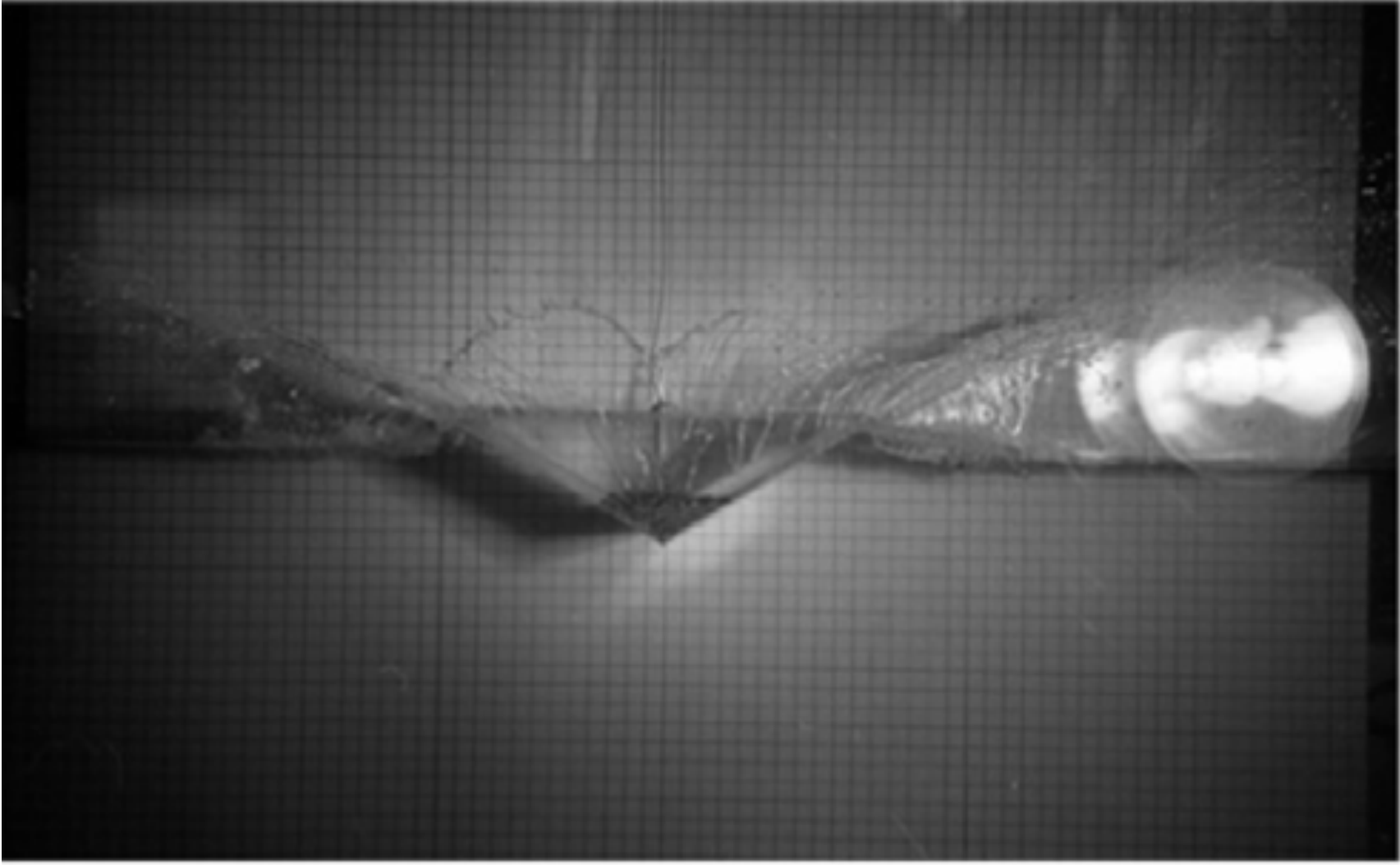}
    \label{wedge_expt0}
  }
    \subfigure[\REVIEW{Simulation~\cite{Nguyen2016}}]{
    \includegraphics[width = 0.3\textwidth]{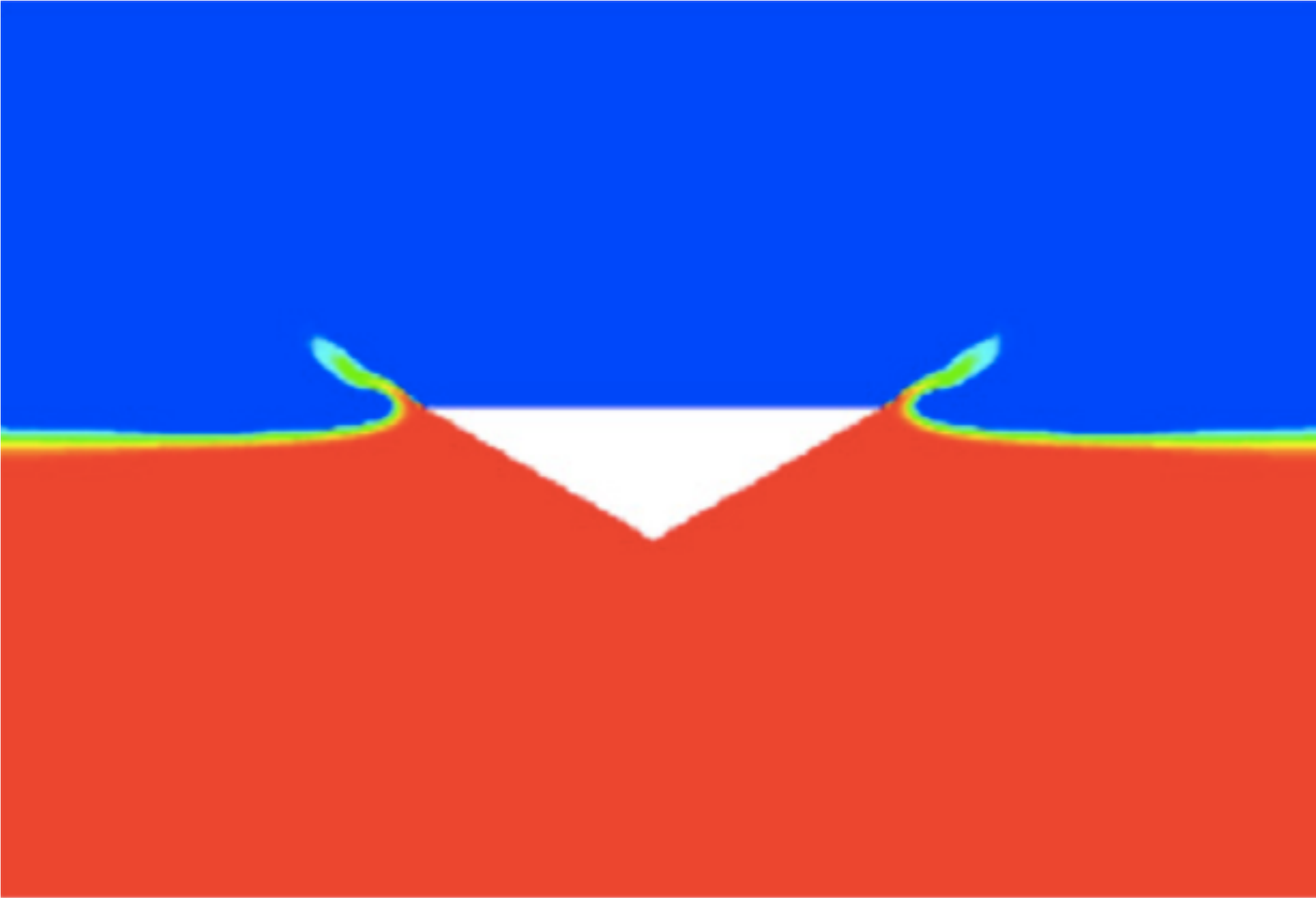}
    \label{wedge_nguyen0}
  }
    \subfigure[\REVIEW{Present FD/BP simulation}]{
    \includegraphics[width = 0.29\textwidth]{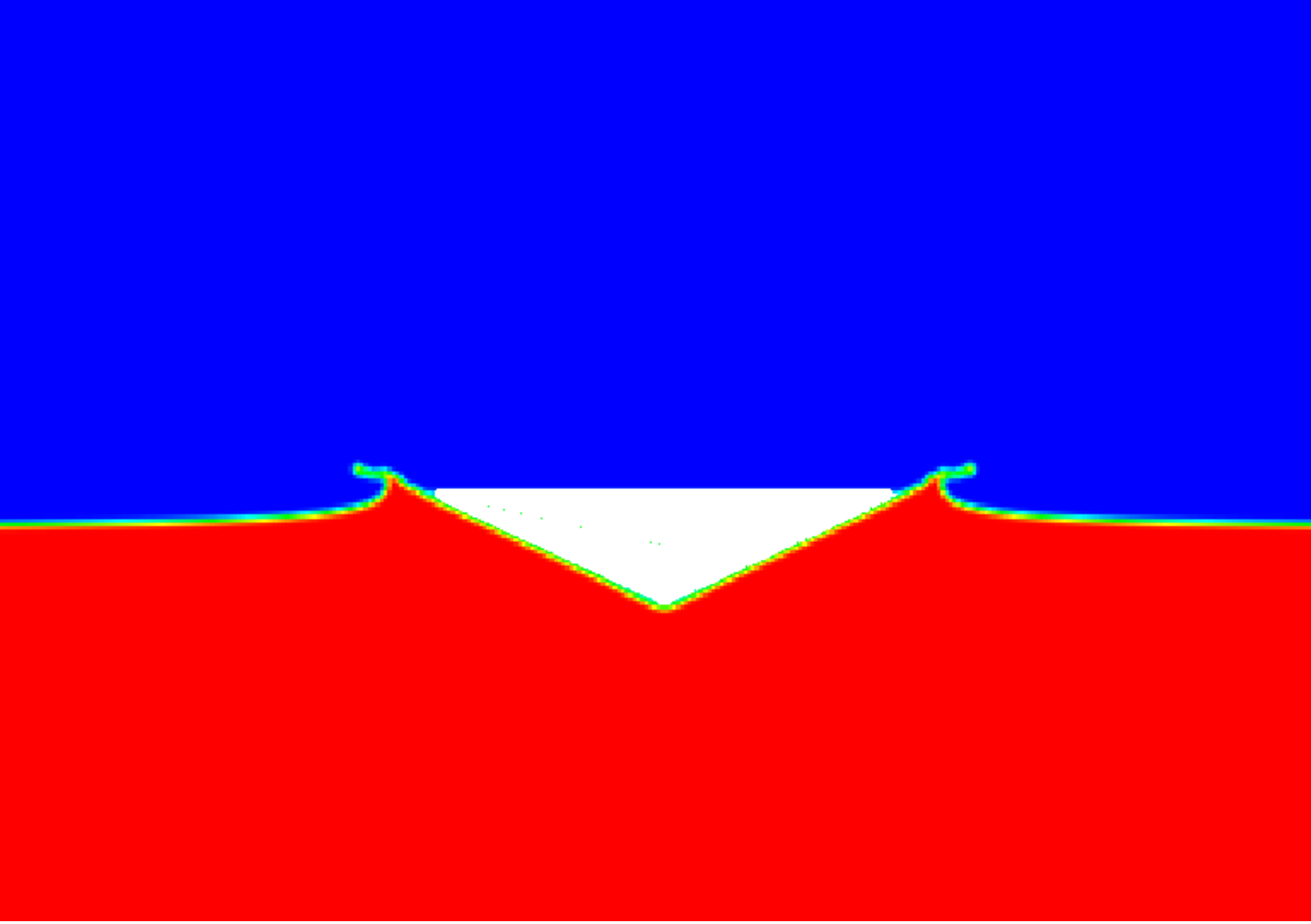}
    \label{wedge_present0}
  }
  \subfigure[\REVIEW{Experiment~\cite{Greenhow1983}}]{
    \includegraphics[width = 0.33\textwidth]{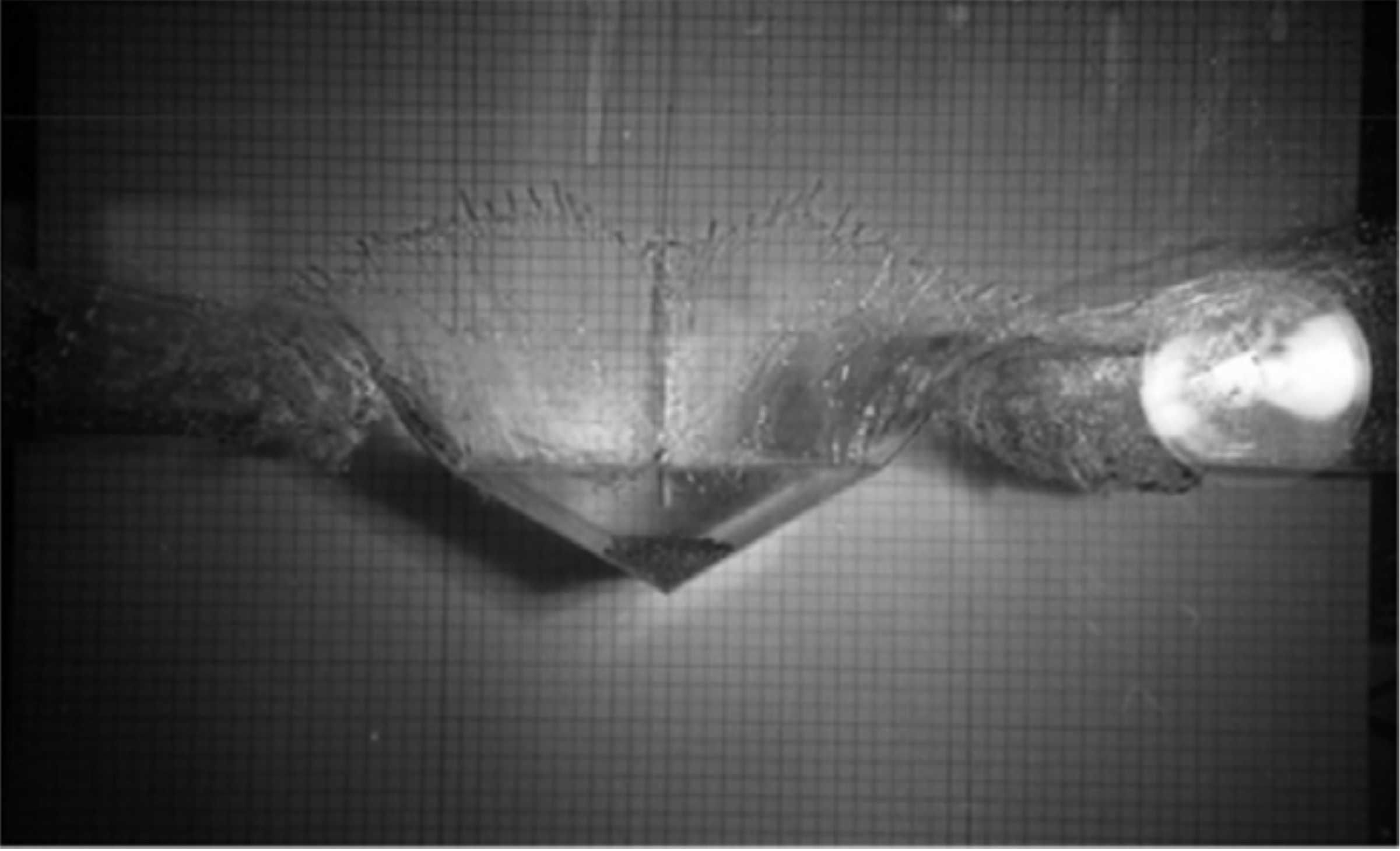}
    \label{wedge_expt1}
  }
    \subfigure[\REVIEW{Simulation~\cite{Nguyen2016}}]{
    \includegraphics[width = 0.3\textwidth]{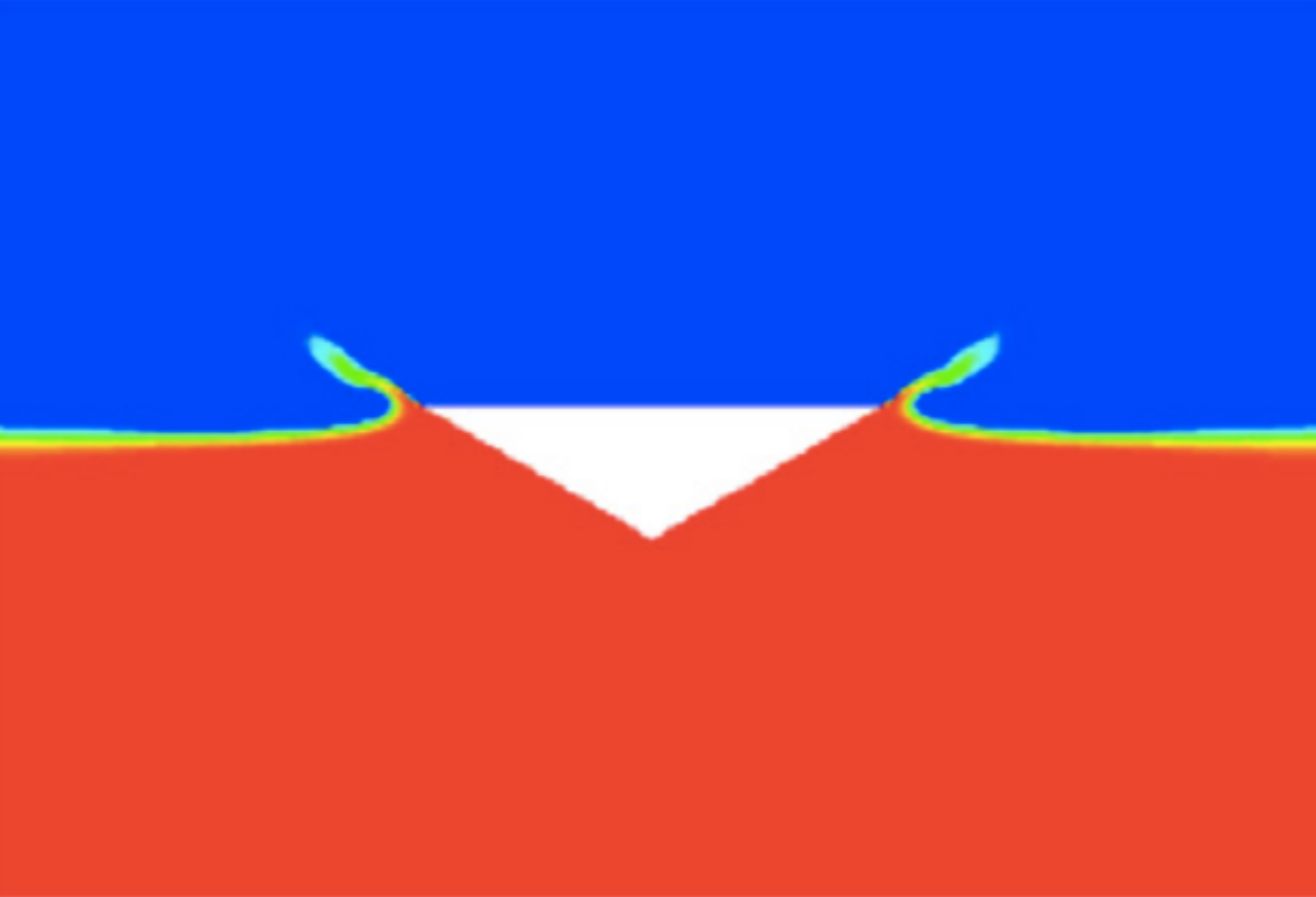}
    \label{wedge_nguyen1}
  }
    \subfigure[\REVIEW{Present FD/BP simulation}]{
    \includegraphics[width = 0.29\textwidth]{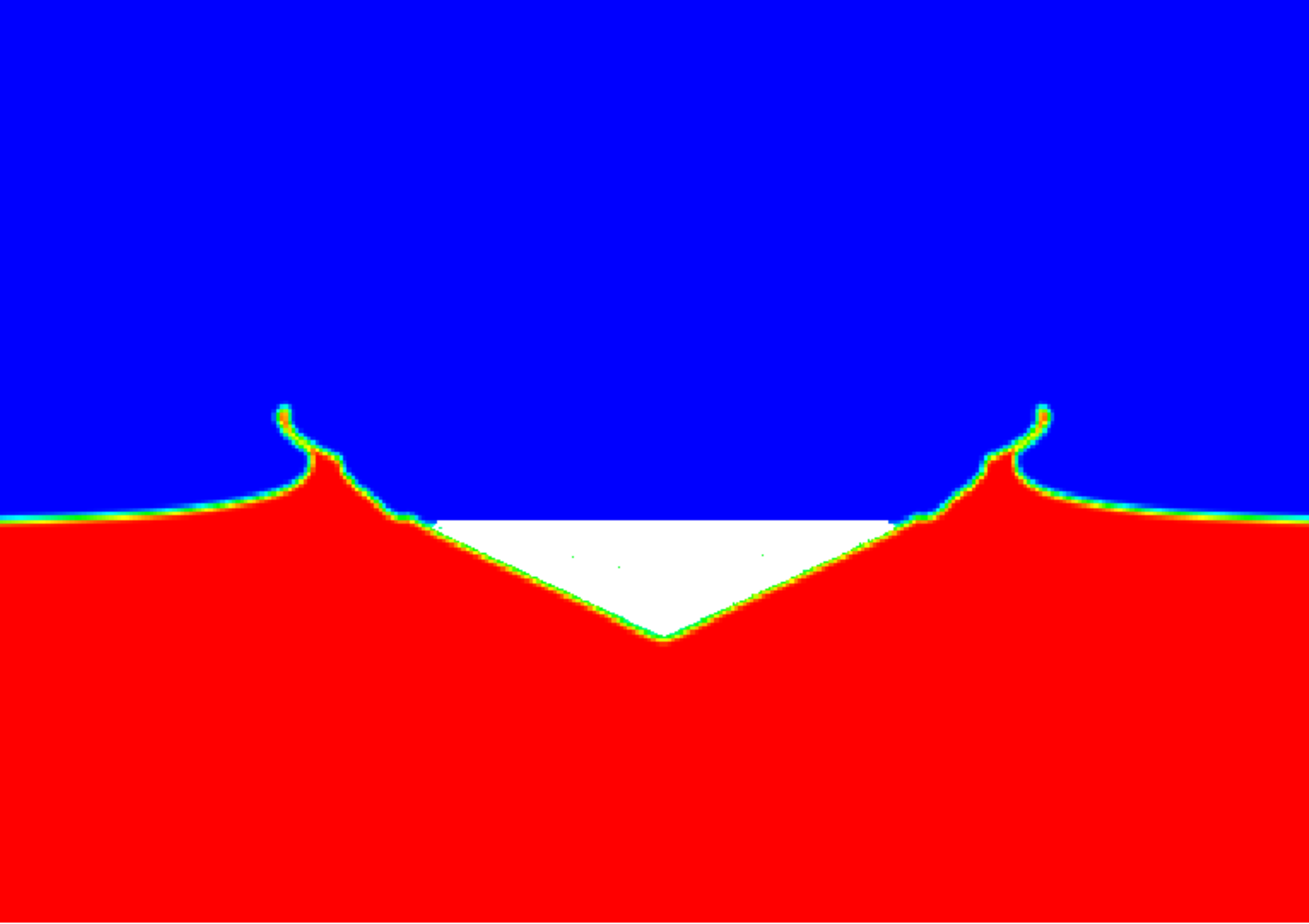}
    \label{wedge_present1}
  }
  \caption{\REVIEW{Visual comparisons of density and free surface: 
  \subref{wedge_expt0} and \subref{wedge_expt1} experimental photographs from Greenhow and Lin~\cite{Greenhow1983};
  \subref{wedge_nguyen0} and \subref{wedge_nguyen1} simulation snapshots from Nguyen et al.~\cite{Nguyen2016};
  \subref{wedge_present0} and \subref{wedge_present1} simulation snapshots from the present the FD/BP method. Results from the FD/IB method are similar. Image permissions from the journals~\cite{Greenhow1983} and~\cite{Nguyen2016} pending.}
}
  \label{fig_compare_wedge2d_viz}
\end{figure}

 Table~\ref{tab_wedge} compares the time and velocity of impact 
obtained from FD/BP, FD/IB and Newton's second law of motion. Aerodynamic air 
resistance is neglected from Newton's law of motion. Both methods are in reasonable agreement with each 
other and also agree with the analytically predicted impact time and velocity.  \REVIEW{Fig.~\ref{fig_compare_wedge2d_viz} compares the 
initial interfacial dynamics of wedge impact with prior experimental~\cite{Greenhow1983} and numerical studies~\cite{Nguyen2016};
decent agreement is seen.}  Fig.~\ref{fig_wedge2d_viz} shows the 
evolution of fluid-structure interaction along with the vorticity generated by the two methods. 
\REVIEW{Upon impact, the FD/IB method sheds two counter-rotating vortices that are oriented inwards, whereas the FD/BP method, upon impact, sheds them in a slightly outward orientation (see time panel $t = 0.5625$ s of Fig.~\ref{fig_wedge2d_viz}). This is attributed to the differences in the impact forces and velocities predicted by the two methods. Another difference is in the fluid-structure interface handling in the two methods. This can also have some minor effects on the vortex shedding dynamics at high Reynolds numbers.}
There is also a slight delay in the vortical dynamics of the FD/BP method compared to the FD/IB method which can be explained by considering 
the lag in the impact time predicted by the former method. Similarly, the higher impact force of the FD/IB method as seen 
in Fig.~\ref{wedge2d_force} can be attributed to a higher impact velocity as compared to the FD/BP method. 

\begin{table}
    \centering
    \caption{Water impact time and velocity of a free falling 2D wedge in air, computed using FD/BP, FD/IB, and Newton's second law of motion.}
    \rowcolors{2}{}{gray!10}
    \begin{tabular}{*6c}
        \toprule
        Method & $t_{\text{impact}}$ (s) & $v_{\text{impact}}$  (m/s) \\
        \midrule         
        FD/BP & $\REVIEW{0.5253}$ & $\REVIEW{4.7409}$ \\
        FD/IB  & $0.5189$ & $4.8347$  \\
        Newton's law & $ \sqrt{2 \Delta s/g} = 0.5148$ & $g t_{\text{impact}} = 5.0502$ \\
        \bottomrule
    \end{tabular}
    \label{tab_wedge}
\end{table}

\subsection{Water-entry/exit of a free falling cylinder}

In this section, we investigate the problem of a half-buoyant cylinder freely falling in water.
This case has been studied numerically by Sun et al.~\cite{Sun2015} using a weakly compressible 
smoothed particle hydrodynamics (SPH) method, \REVIEW{and by Patel and Natarajan~\cite{Patel2018} using an incompressible 
volume of fluid (VOF) solver}. A circular cylinder of diameter $D = 0.11$ m 
and density $\rhos = 500$ kg/m$^3$ is placed in a two dimensional computational 
domain of size $\Omega = [0,20D] \times [0, 12D]$ with initial center position $(X_0, Y_0) = (10D, 8.05D)$. 
The domain is filled from $y = 0$ to $y = 3D$ with water; the remainder of the tank from $y = 3D$ to 
$y = 12D$ is filled with air. The cylinder has a free fall height of $\Delta s = 4.55 D$. The domain is 
discretized using a $880 \times 528$ uniform grid, \REVIEW{which corresponds to 44 cells per diameter. This grid resolution was 
found sufficient in our prior work for similar water-impact cases~\cite{Nangia2019}.} A constant time step size $\dt = 10^{-5}$ s is used.

\begin{figure}[t!]
  \centering
  
   \subfigure[Vertical position]{
    \includegraphics[scale = 0.3]{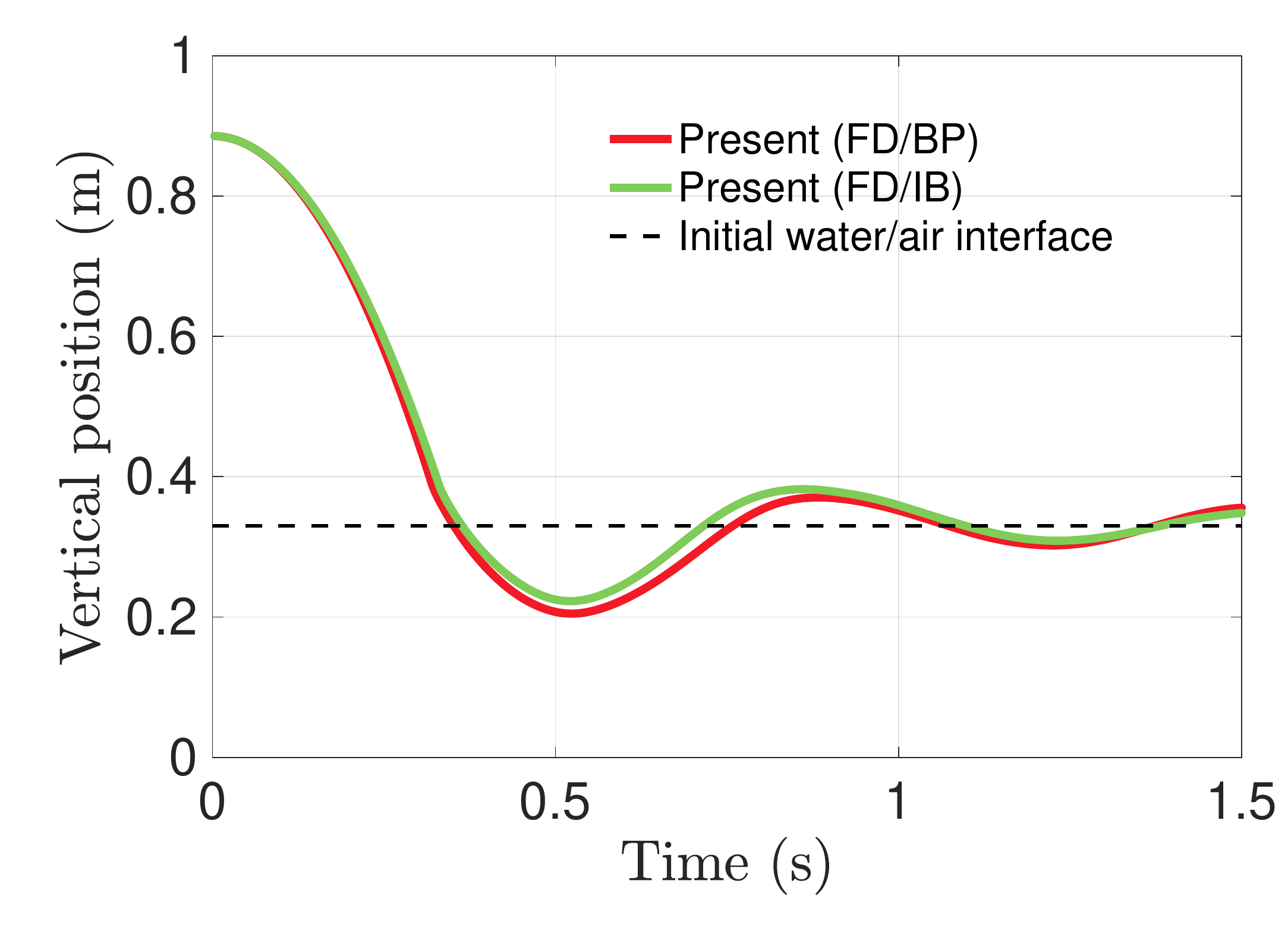}
    \label{cyl_posn}
  }
     \subfigure[Vertical velocity]{
    \includegraphics[scale = 0.3]{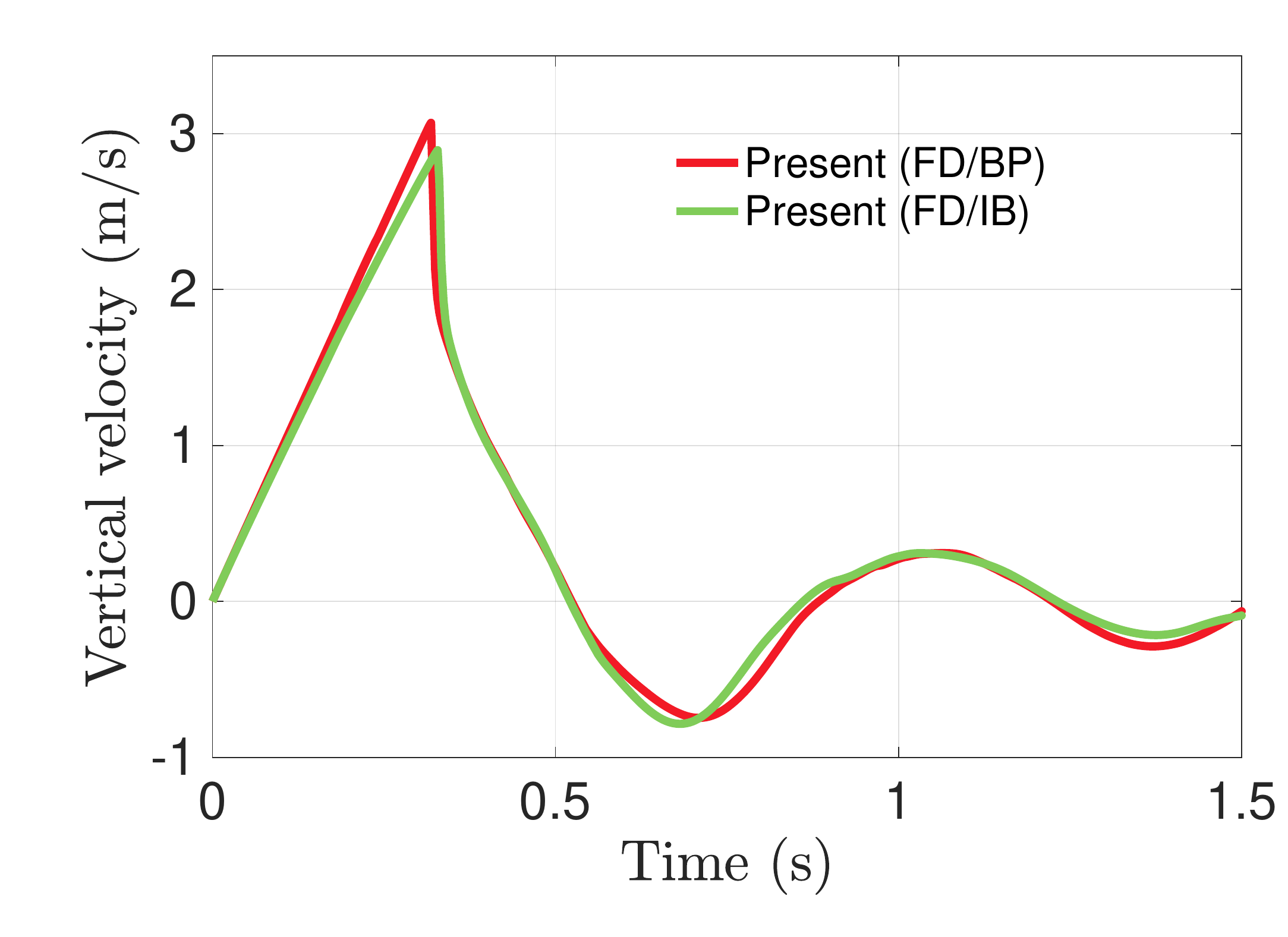}
    \label{cyl_vel}
  }
   \subfigure[Vertical force]{
    \includegraphics[scale = 0.3]{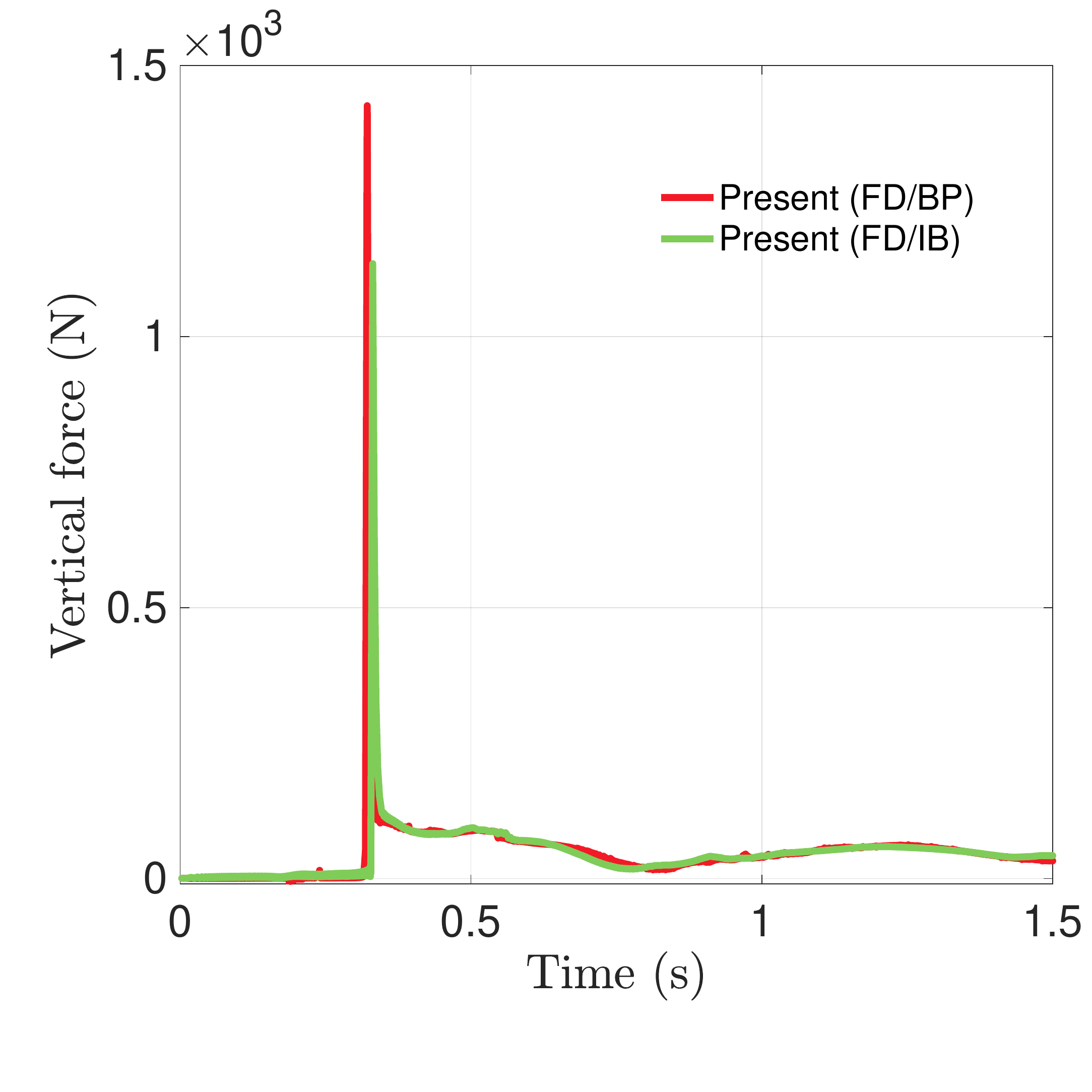}
    \label{cyl_force}
  }
     \subfigure[Penetration depth]{
    \includegraphics[scale = 0.3]{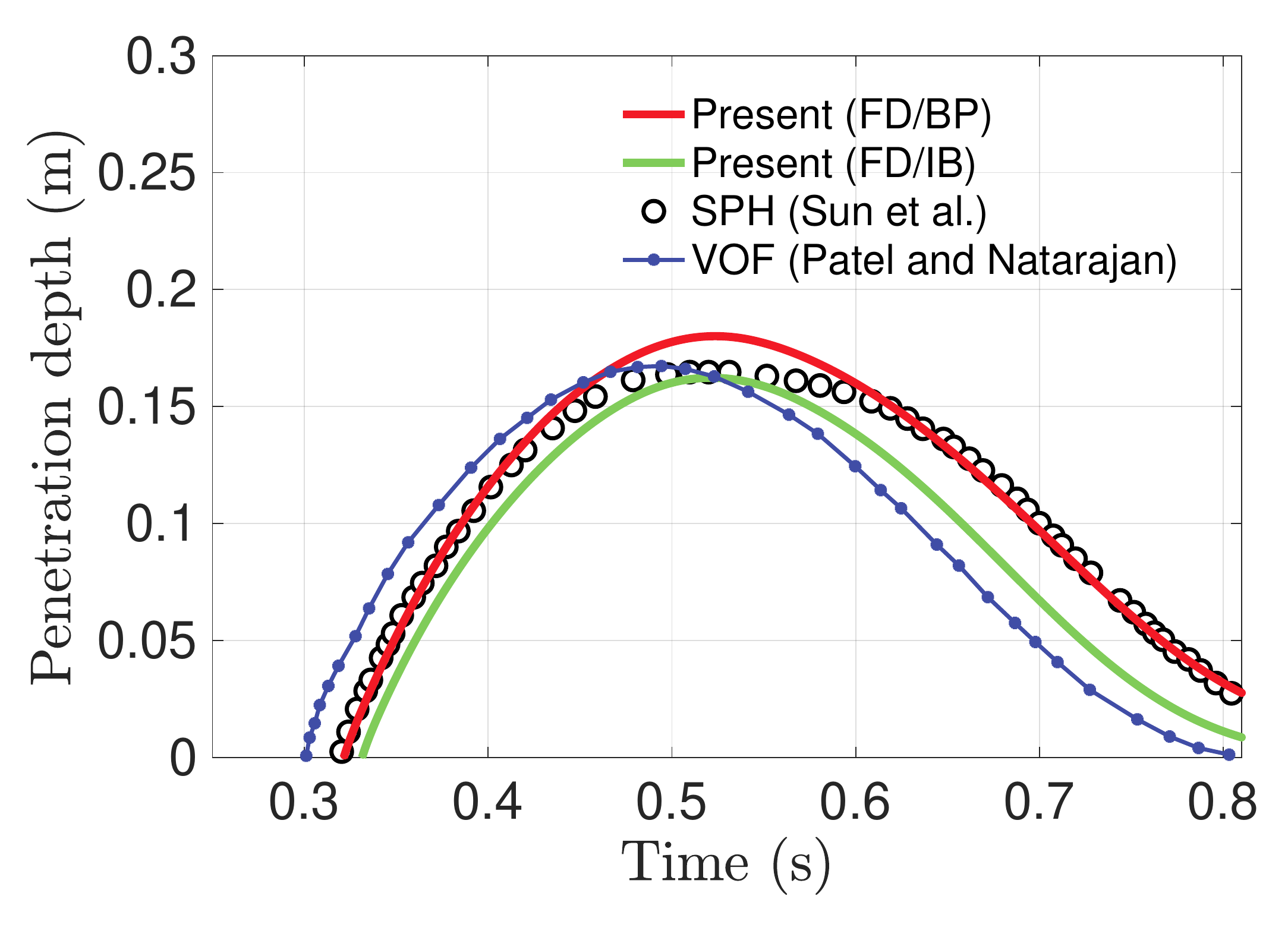}
    \label{cyl_sun}
  }
  \caption{Temporal evolution of \subref{cyl_posn} vertical position, \subref{cyl_vel} vertical velocity, and \subref{cyl_force} 
  vertical force on a half-buoyant cylinder freely falling in water. (---, red) present FD/BP simulation data; ; (---, green) present FD/IB simulation data;
  ($\circ$, black) SPH simulation data from Sun et al.~\cite{Sun2015}.
}
  \label{fig_cyl_com}
\end{figure}

We again compare the rigid body dynamics of the cylinder obtained from the two methods. Figs.~\ref{cyl_posn} 
and~\ref{cyl_vel} show the time evolution of the center of mass vertical position and velocity, respectively. 
The hydrodynamic forces in the vertical direction obtained from the two methods are plotted in Fig.~\ref{cyl_force}. 
Fig.~\ref{cyl_sun} compares the variation of the depth of penetration as time progresses with the 
prior \REVIEW{numerical studies~\cite{Sun2015,Patel2018}}. An excellent agreement is found between 
the FD/BP and SPH methods for most of the times. The FD/IB method gives slightly reduced penetration depth at later times 
compared to the FD/BP method. \REVIEW{However, FSI results obtained from both FD implementations fall in the range of prior numerical studies.}  

\begin{figure}[]
  \centering
  
   \subfigure[FD/BP, $t = 0.465$ s]{
    \includegraphics[scale = 1.15]{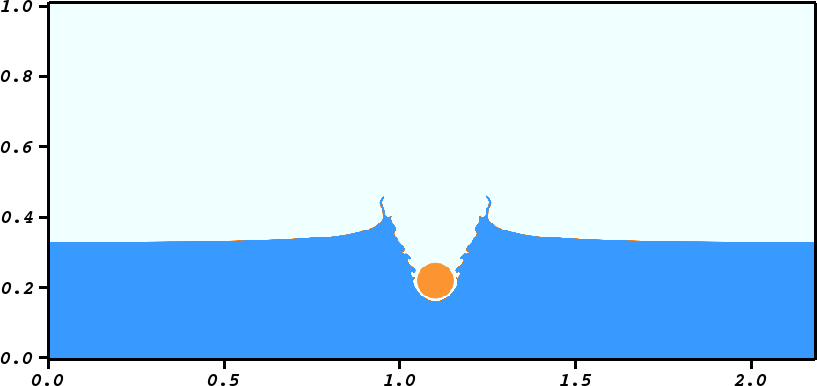}
    \label{ibls_cyl_t0p465}
  }
     \subfigure[FD/IB, $t = 0.465$ s]{
    \includegraphics[scale = 1.15]{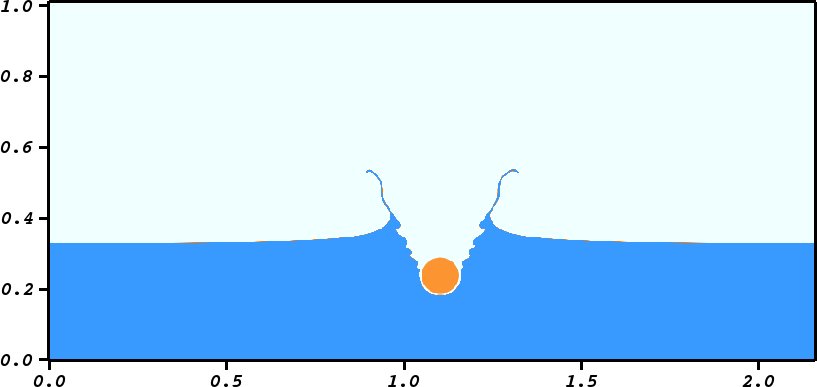}
    \label{cib_cyl_t0p465}
  }
   \subfigure[FD/BP, $t = 0.56$ s]{
    \includegraphics[scale = 1.15]{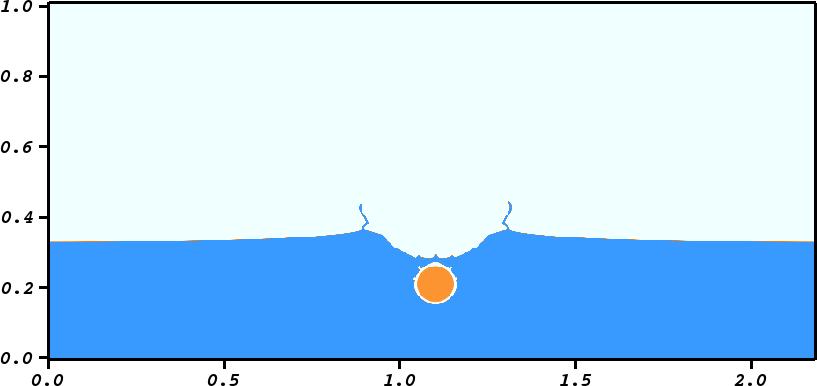}
    \label{ibls_cyl_t0p56}
  }
    \subfigure[FD/IB, $t = 0.56$ s]{
    \includegraphics[scale = 1.15]{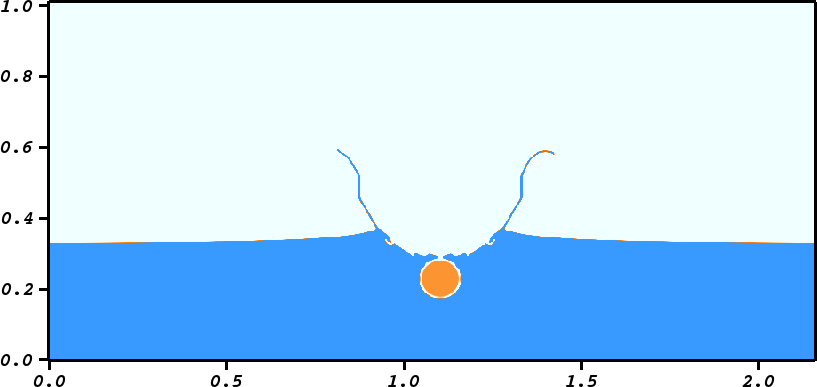}
    \label{cib_cyl_t0p56}
  }
  \subfigure[FD/BP, $t = 0.735$ s]{
    \includegraphics[scale = 1.15]{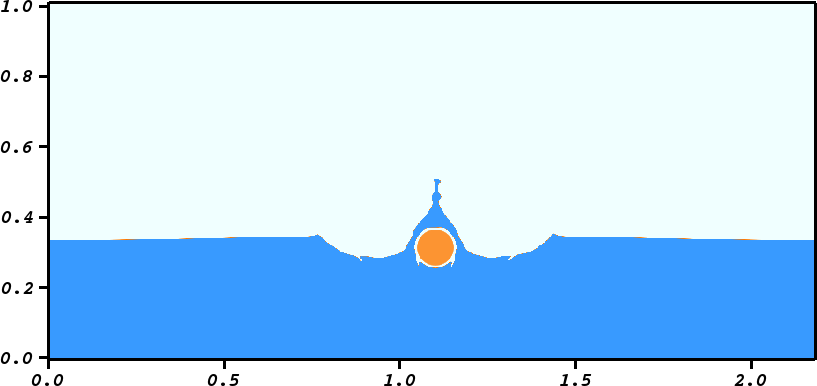}
    \label{ibls_cyl_t0p735}
  }
    \subfigure[FD/IB, $t = 0.735$ s]{
    \includegraphics[scale = 1.15]{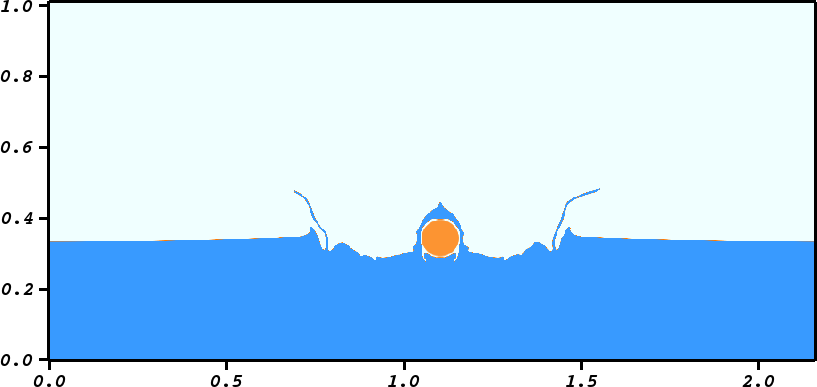}
    \label{cib_cyl_t0p735}
  }
    \subfigure[FD/BP, $t = 1.365$ s]{
    \includegraphics[scale = 1.15]{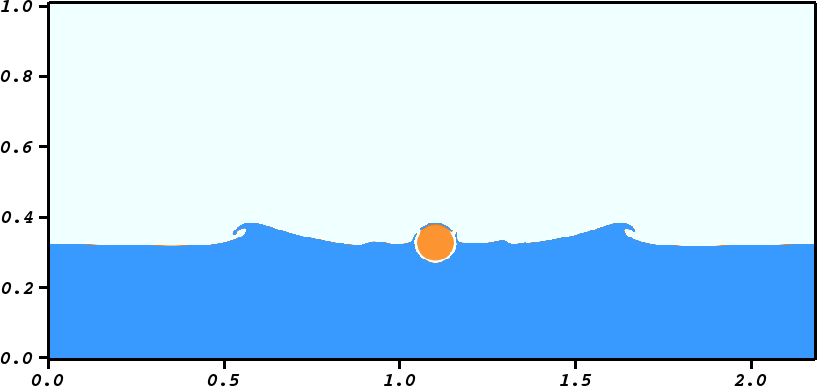}
    \label{ibls_cyl_1pt25}
  }
  \subfigure[FD/IB, $t = 1.365$ s]{
    \includegraphics[scale = 1.15]{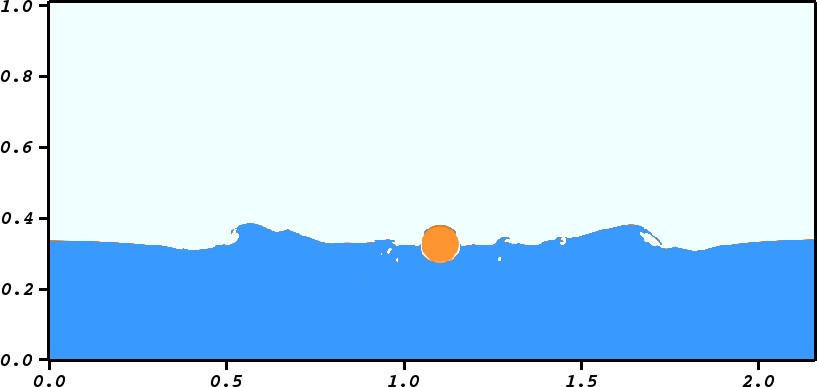}
    \label{cib_cyl_1pt25}
  }
  \caption{Splash dynamics generated by a half-buoyant cylinder freely falling in water at four different time 
  instances using the FD/BP and the FD/IB methods. Lighter blue color represents the air phase whereas darker blue color represents 
  the water phase. The solid phase is shown in orange shade. 
}
  \label{fig_cyl_viz}
 
\end{figure}

Table~\ref{tab_cyl} compares the time and velocity of impact obtained from FD/BP, FD/IB, prior numerical studies, and Newton's second law of motion. 
Again, both methods are in reasonable agreement with each other and also agree \REVIEW{well} with the analytically predicted impact time and velocity. \REVIEW{The impact time simulated by Patel and Natarajan underpredicts the expected
value of $t_{\text{impact}}$. Sun et al. release the cylinder with an initial velocity of $v_{\text{impact}}$ (as predicted by Newton's law) 
at the air-water interface and do not simulate the free-fall motion of the cylinder in the air phase.}
Fig.~\ref{fig_cyl_viz} shows the time evolution of interfacial dynamics using the two methods. The impacting 
cylinder produces distinct water jets while moving downward into the liquid. \REVIEW{The initial impact of the cylinder 
produces ragged and non-smooth deformations in the separated water layer, as seen distinctly at $t = 0.465$ s. These deformations are also 
observed in the weakly compressible SPH simulations of Sun et al.~\cite{Sun2018} and the incompressible VOF simulations of Patel and Natarajan~\cite{Patel2018}. Sun et al. attribute these deformations to negative pressure regions created in the receded 
water layer and numerically ``fix" them by zeroing out the negative pressure. They refer this fix as \emph{a numerical model of water 
repellent coating on the cylinder surface}.}  At around $t = 0.52$ s, the cylinder reverses its direction of 
motion and pushes a layer of liquid along its surface as it rises up at $t = 0.735$ s. Eventually, two 
opposite traveling waves on either side of the cylinder are formed when the 
cylinder enters back into the water for the second time. This can be distinctly seen at $t = 1.365$ s. The slightly stronger 
water jets produced by the FD/IB method (compared to the FD/BP method) at the initial impact affects the interfacial 
dynamics at even later times. They also cause slight asymmetries on opposite sides of the cylinder at later times. 
In contrast, the FD/BP method maintains interfacial symmetry for most of the times shown.    

\begin{table}
    \centering
    \caption{Water impact time and velocity of a free falling cylinder in air computed by FD/BP, FD/IB, Patel and Natarajan~\cite{Patel2018}, Sun et al.~\cite{Sun2015}, and Newton's second law of motion.}
    \rowcolors{2}{}{gray!10}
    \begin{tabular}{*6c}
        \toprule
        Method & $t_{\text{impact}}$ (s) & $v_{\text{impact}}$  (m/s) \\
        \midrule
        FD/BP & $0.3187$ & $3.0690$ \\
        FD/IB  & $0.3287$ & $2.8936$  \\
        \REVIEW{Patel and Natarajan~\cite{Patel2018}} &  \REVIEW{$0.2997$} & \REVIEW{N/A} \\
        \REVIEW{Sun et al.~\cite{Sun2015}} & \REVIEW{0.3194}  & \REVIEW{3.1337} \\
        Newton's law & $ \sqrt{2 \Delta s/g} = 0.3194$ & $g t_{\text{impact}} = 3.1337$ \\
        \bottomrule
    \end{tabular}
    \label{tab_cyl}
\end{table}

%Comparison with Natarajan and Sun for half-bouyant cylinder.

%For Natarajan:

%Non-dimensional time of impact T* = 2.83

%To get dimensional time we have

%R = 0.055;
%D = 2*R;
%G = 9.81;
%V_scale = sqrt(G*D).  = 1.0388
%T_scale = D/V         = 0.1059    

%T_impact = T_scale * (T*) = 0.1059*2.83 =  0.2997 s  

%Insert this time of impact in the Table and in the text mention that FD/IB and FD/BP lie in the range of the previously reported data and somewhat better than Natarajan. Sun?s SPH method does not simulate the free drop --- it just imposes v_impact at t = 0. I just shifted the time of impact by time predicted by Newton?s method.  

\REVIEW{\section{Comparison of computational costs} \label{sec_cost_comparison}
Finally, we briefly discuss and compare the computational costs associated with both methods. We consider
the same two-dimensional test case described in the previous section: a half buoyant cylinder falling into an
air-water interface. We discretize the problem on a coarser mesh of size $440 \times 263$ so that the problem
can be run relatively quickly on a single processor. A constant time step size of $\dt = 1 \times 10^{-5}$ is used.
This test problem is simulated for $60$ time steps with both the FD/IB and FD/BP methods. The wall-clock time is measured
for the final $50$ time steps of each simulation, and three experimental trials are conducted (and averaged) for each method. 
All of the solver options are identical across each trial. The total wall-clock time is broken up into four categories:
\begin{enumerate}
\item `INS Solver', which corresponds to the operations required to solve the discrete fluid flow equations described by 
Eqs.~\eqref{eq_c_discrete_momentum} and~\eqref{eq_c_discrete_continuity} in Sec.~\ref{sec_cons_ins}.
\item `Level Set Update', which corresponds to the operations required to discretely advect the level set variables
(Eqs.~\eqref{eq_dis_ls_fluid} and~\eqref{eq_dis_ls_solid} in Sec.~\ref{sec_scalar_adv}, and reinitialize them to signed
distances functions (see Appendix~\ref{app_reinit}).
\item  `BP FSI Correction' or `IB FSI Correction', which correspond to the operations required to correct the fluid velocity
in the domain occupied by the structure (described in Sec.~\ref{sec_ts_bp} for the FD/BP method and in Sec.~\ref{sec_ts_ib}
for the FD/IB method.
\item `Other', which corresponds to any operations not covered by the previous three descriptors. This includes allocation
and deallocation of data structures and various pre- and post-processing function calls. 
\end{enumerate}

Fig.~\ref{fig_comp_breakdown} shows the computation breakdown for each method. It is clear that discretely solving the fluid
flow equations is by far the costliest operation, taking over $95\%$ of the total computation time for both methods. In 
contrast, the level set advection and reinitialization routines take very little wall-clock time; this is unsurprising since there
are no Krylov iterative solvers employed in computing $\phi^{n+1}$ and $\psi^{n+1}$ -- the update is purely explicit in 
nature. Finally, we note that the FSI correction for each method is incredibly fast, each taking only a fraction of a percent
of total computation time. These results show that the present strong FSI coupling schemes are extremely efficient.

\begin{figure}[t!]
  \centering
   \subfigure[\REVIEW{FD/IB computation breakdown}]{
    \includegraphics[scale = 0.3]{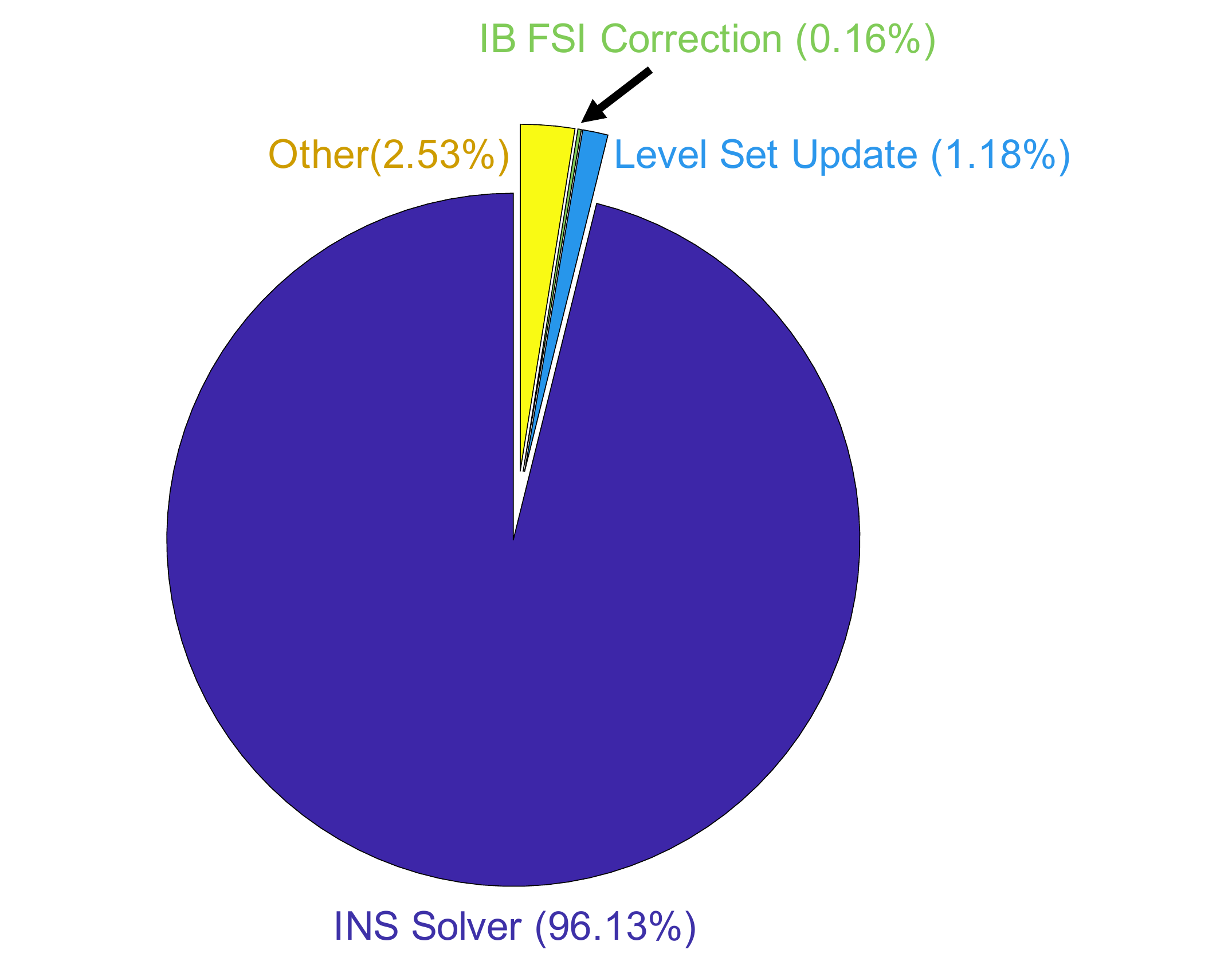}
    \label{fd_ib_timing}
  }
  \subfigure[\REVIEW{FD/BP computation breakdown}]{
    \includegraphics[scale = 0.3]{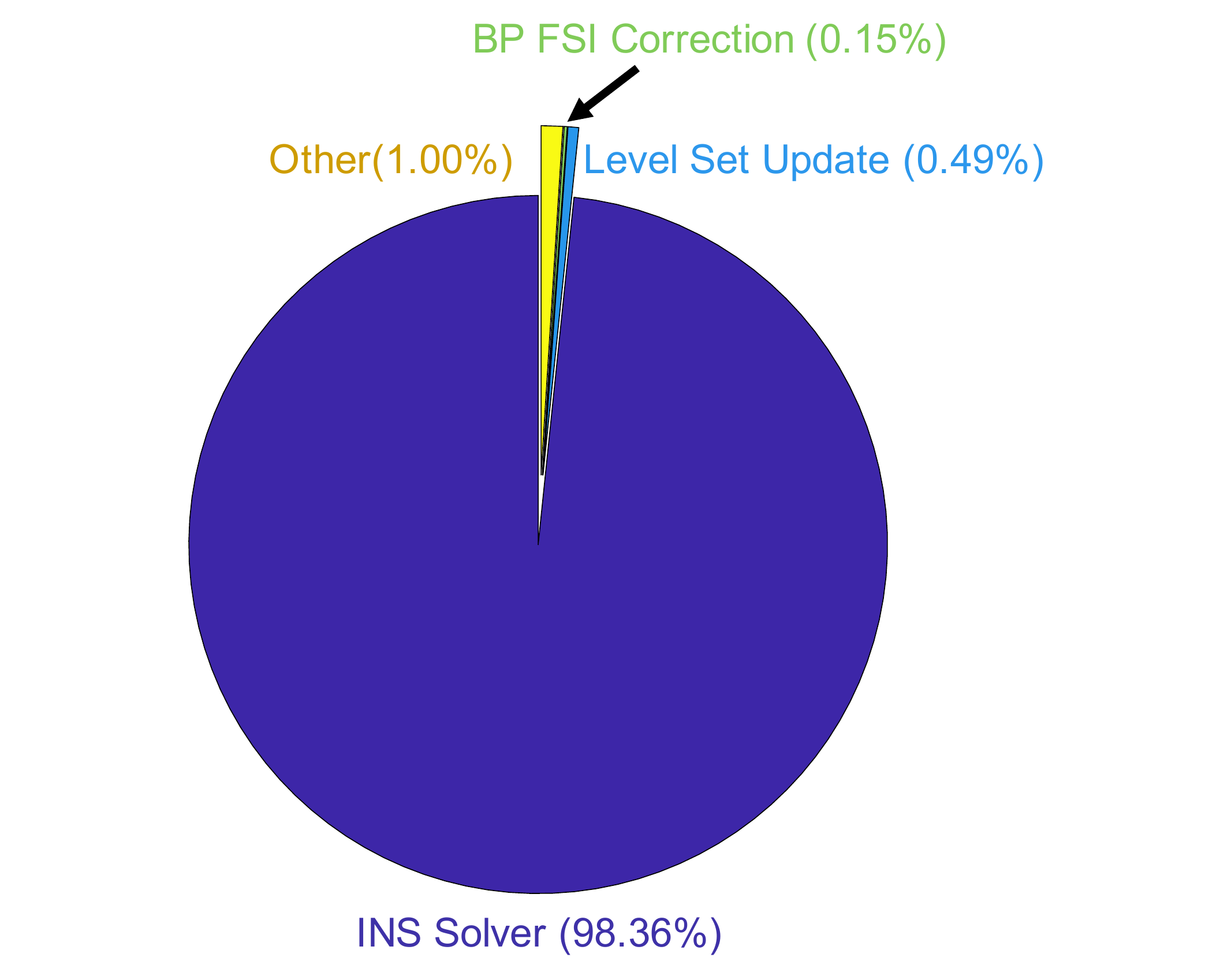}
    \label{fd_bp_timing}
  }
  \caption{\REVIEW{Breakdown of computing time spent for the
  \subref{fd_ib_timing} FD/IB and 
  \subref{fd_bp_timing} FD/BP methods.
  Note that the fraction of time spent for the FSI coupling computations is barely visible at this scale.
  All runs were carried out on a single processor.}
}
  \label{fig_comp_breakdown}
\end{figure}

Table~\ref{tab_wallclock} shows the average wall-clock time (in CPU units) required to compute each major component
of the FD/BP and FD/IB algorithms. Unsurprisingly, the level set update takes approximately the same amount of time
for both methods. Moreover, it is clear the the FSI coupling for the FD/BP method is more expensive than the FD/IB method.
This can be attributed to the fact that the Brinkman penalization approach requires additional computations of
hydrodynamic forces and torques (e.g. Eqs.~\eqref{eq_int_fh} and~\eqref{eq_int_mh}) and setting the matrix entries of the modified implicit fluid solver. 
Finally, it can be seen that the fluid solver for the FD/BP method is computationally more expensive than for the FD/IB method. This is because
of the additional term present in the discrete INS equations due to the Brinkman penalization formulation ($\thetaFD = 1$ in Eq.~\eqref{eq_c_discrete_momentum}). This term is treated implicitly, which changes the underlying discrete operator (and
hence the linear solver convergence properties). For this particular example, the fluid solver for the FD/IB method converged 
in $1-2$ iterations each time step, while the fluid solver for the FD/BP  (using the same solver settings) converged in $2-4$ iterations, 
explaining the disparity in wall-clock time.

We expect that the computational cost results would be similar for the other two cases described in this work. The FD/IB
method outperforms the FD/BP method because the two-dimensional Lagrangian meshes considered in the present work are 
relatively simple and small. Moreover, we should note that the implementation of the IB transfer operators (e.g. spreading 
and interpolation) is quite mature; these computations have 
been well optimized and vectorized over the past two decades or so~\cite{IBAMR-web-page}.
In contrast, the FD/BP has been recently implemented in IBAMR and we have not fully investigated speedup efforts.
In particular, we believe that development of special multigrid smoothers (or other solver settings) for the implicit Brinkman penalization equations
would help to reduce the number of Krylov iterations. These could make the computational cost of the FD/BP method
comparable to the FD/IB method even for these small-scale 2D problems.  We note that a more fair comparison of 
computational cost would be to compare timing data between the FD/IB method and an \emph{explicit} version of the FD/BP 
method. The multigrid smoothers and an explicit FD/BP method will be explored in future endeavors.

For full, 3D engineering geometries (such as a wave energy converter device), we expect a critical problem size
to exist at which the FD/BP method would outperform the FD/IB method. This is especially true when considering distributed
memory parallelism; when two separate meshes exist (i.e. Lagrangian and Eulerian), two different domain decompositions
must be maintained and mapped onto one another to achieve an efficient and scalable method. Even though the implicit 
FD/BP method will in general require more Krylov iterations, we expect that its purely Eulerian formulation will be advantageous
for large engineering problems requiring scalable load balancing.
The run time ratios and trends described in Fig.~\ref{fig_comp_breakdown} and Table~\ref{tab_wallclock} would not necessarily hold for large, multi-processor runs. This will be the subject of future investigations.

\begin{table}
    \centering
    \caption{\REVIEW{Average ($n = 3$) wall-clock time required to compute each major component of the FD/BP and FD/IB algorithms. All runs were carried out on a single processor.}}
    \rowcolors{2}{}{gray!10}
    \begin{tabular}{*6c}
        \toprule
        Method & INS Solver (CPU units) & Level Set Update (CPU units) & FSI Coupling (CPU units)  \\
        \midrule
        FD/BP & $541.79$ & $2.82$ & $0.79$ \\
        FD/IB  & $229.57$ & $2.73$ & $0.38$ \\
        \bottomrule
    \end{tabular}
    \label{tab_wallclock}
\end{table}

}

%%%%%%%%%%%%%%%%%%%%%%%%%%%%%%
\section{Conclusions}
In this study, we described two implementations of the fictitious domain method capable of simulating water-entry/exit
problems. One algorithm, based on the immersed boundary method, relied on an Eulerian description of the fluid variables
and a Lagrangian representation of the immersed structure. The second algorithm, based on Brinkman penalization, was
a purely Eulerian approach that imposed constraint forces implicitly rather than explicitly. We demonstrated that both methods 
can adequately resolve complex floating and splashing dynamics that are ubiquitous in practical marine engineering problems.
They provide a good alternative to overset mesh based methods for simulating complex FSI problems. 

For both methods, standard level set machinery was used to track air-water interfaces and the surface of the immersed body. 
The similarities and differences between the FD/BP and FD/IB methods were discussed, and several advantages and 
disadvantages of each technique were also described in Sec.~\ref{sec_comp_fd}.
An efficient method for computing hydrodynamic forces and torques was detailed as a part of the FD/BP solution algorithm.
We compared the techniques using two standard test cases: a half-buoyant, free-falling wedge and cylinder impacting an 
air-water interface. Both methods produced results \REVIEW{that are in reasonable agreement with each other} and compared favorably with prior 
results shown in literature.

The techniques described here have been implemented within the open-source IBAMR library. IBAMR is a flexible software infrastructure that 
provides support for several versions of the immersed boundary and fictitious domain methods. We are actively working
on extending these methods to solve more complex problems. Future work includes implementation
of a RANS or LES turbulence model and computational geometry algorithms to initialize and transport level set functions of more 
sophisticated solid geometries.

\section*{Acknowledgements}
A.P.S.B.~acknowledges research support provided by the San Diego State University
and the College of Engineering's Fermi high performance computing service.
N.N~acknowledges research support from the National Science Foundation Graduate Research Fellowship Program (NSF award DGE-1324585)

%%%%%%%%%%%%%%%%%%%%%%%%%%%%%%%%%%%%%%%%%%%
\appendix 
\renewcommand\thesection{\Alph{section}}
\REVIEW{
\section{Level set reinitialization}
\label{app_reinit}
It is well-known that the signed distance property of $\phi$ and $\psi$ is disrupted under linear advection,
Eqs.~\eqref{eq_ls_fluid_advection} and~\eqref{eq_ls_solid_advection}. Letting $\widetilde{\phi}^{n+1}$
denote the flow level set function following a single time step of transport through the interval $\left[t^{n}, t^{n+1}\right]$,
we need a procedure that \emph{reinitializes} the field into a signed distance function $\phi^{n+1}$.
This can be achieved by computing a steady-state solution to the Hamilton-Jacobi equation
\begin{align}
&\D{\phi}{\tau} + \sgn\left(\widetilde{\phi}^{n+1}\right)\left(\|\grad \phi \| - 1\right) = 0, \label{eq_eikonal} \\
& \phi(\x, \tau = 0) = \widetilde{\phi}^{n+1}(\x), \label{eq_eikonal_init}
\end{align}
which will yield a solution to the Eikonal equation $\|\grad \phi \|  = 1$ at the end of each time step. 
More details on the specific discretization of
Eqs.~\eqref{eq_eikonal} and~\eqref{eq_eikonal_init}, which employs second-order
ENO finite differences combined with a subcell-fix method described by Min~\cite{Min2010},
and an immobile interface condition described by Son~\cite{Son2005}, can be found in~\cite{Nangia2018}. 

Since we consider simple geometries in this work, the solid level set $\psi^{n+1}$ is analytically calculated 
by using the new location of center of mass at $t^{n+1}$. For more complex structures,
computational geometry techniques can be employed to compute the signed distance function.

\section{Discretization of the convective term: consistent mass/momentum transport}
\label{app_convective}
We use an explicit cubic upwind interpolation (CUI-limited) scheme~\cite{Roe1982,Waterson2007,Patel2015} 
to approximate the $\C^{n+1,k}$ nonlinear term in the momentum Eq.~\eqref{eq_c_discrete_momentum}.
A discretized mass balance equation is integrated directly on the faces of the staggered grid to obtain the newest 
approximation to density $\breve{\V \rho}^{n+1,k+1}$ in Eq.~\eqref{eq_c_discrete_momentum} from the previous time step 
and level set synchronized density field $\V \rho^{n}$ (obtained after averaging $\phi^{n}$ and $\psi^{n}$ onto faces):
	\begin{align}
	& \breve{\V \rho}^{(1)} = \V \rho^{n} - \dt \R\left(\u^{n}_\text{adv}, \V \rho^{n}_\text{lim}\right), \label{eq_rk1}\\
	& \breve{\V  \rho}^{(2)} = \frac{3}{4} \V \rho^{n} + \frac{1}{4} \breve{\V \rho}^{(1)} - \frac{1}{4} \dt \R\left(\u^{(1)}_\text{adv}, \breve{\V \rho}^{(1)}_\text{lim}\right), \label{eq_rk2} \\
	& \breve{\V \rho}^{n+1, k+1} = \frac{1}{3} \V \rho^n + \frac{2}{3} \breve{\V \rho}^{(2)} - \frac{2}{3} \dt \R\left(\u^{(2)}_\text{adv}, \breve{\V \rho}^{(2)}_\text{lim}\right) \label{eq_rk3}.
	\end{align}
	 A third-order accurate strong stability preserving Runge-Kutta (SSP-RK3) time integrator~\cite{Gottlieb2001} is employed for the above 
	 update. Here $\R\left(\u_{\text{adv}}, \V \varrho_{\text{lim}}\right) \approx \left[\left(\div \left(\u_\text{adv} \V \varrho_\text{lim}\right)\right)_{i-\half, j}, \left(\div \left(\u_\text{adv} \V \varrho_\text{lim}\right)\right)_{i, j-\half}\right]$  is an explicit CUI-limited approximation to the linear density 
	advection term; $\V \varrho$ is either $\V \rho$ or $\breve{\V \rho}$. To clarify the various approximations to the density field, we 
	make a distinction between $\breve{\V \rho}$, the density vector obtained 
	via the SSP-RK3 integrator, and $\V \rho$, the density vector that is set from the level set fields. 
	The subscript ``adv'' indicates the interpolated advective velocity on the faces of face-centered 
	control volume, and the subscript ``lim'' indicates the limited value.  Readers are referred to Nangia 
	et al.~\cite{Nangia2018} for details on obtaining  advective and flux-limited fields. Notice that this density update procedure 
	is occurring \emph{within} the overall fixed-point iteration scheme.
	
	In the SSP-RK3 update, we note that $\u^{(1)}$ is an approximation to ${\u}^{n+1}$, and $\u^{(2)}$ is an approximation to ${\u}^{n+\half}$. 	
	Moreover, $\breve{\V \rho}^{(1)}$ is an approximation to $\breve{\V \rho}^{n+1}$, and $\breve{\V \rho}^{(2)}$ is an approximation to $\breve{\V \rho}^{n+\half}$. We obtain these intermediate velocity and density approximations by using suitable interpolation and extrapolation 
	procedures.  For example, for the first cycle ($k = 0$), the velocities are
	\begin{align}
		& \u^{(1)} = 2 \u^n - \u^{n-1}, \\
		& \u^{(2)} = \3half \u^n - \half \u^{n-1}.
	\end{align}
	For all remaining cycles ($k > 0$), the velocities are
	\begin{align} 
		& \u^{(1)} = {\u}^{n+1,k}, \\
		& \u^{(2)} = \frac{3}{8} {\u}^{n+1,k} + \frac{3}{4} \u^{n} - \frac{1}{8} \u^{n-1}.
	\end{align}
	
	To ensure consistent transport of mass and momentum fluxes, the convective derivative in Eq.~\eqref{eq_c_discrete_momentum} 
	is given by
	\begin{equation}
	\C\left(\u^{(2)}_\text{adv}, \breve{\V \rho}^{(2)}_\text{lim}\u^{(2)}_\text{lim}\right) \approx 
	\begin{bmatrix}
		\left(\div \left(\u^{(2)}_\text{adv} \breve{\V \rho}^{(2)}_\text{lim} u^{(2)}_\text{lim}\right)\right)_{i-\half,j} \\ 
		\left(\div \left(\u^{(2)}_\text{adv} \breve{\V \rho}^{(2)}_\text{lim} v^{(2)}_\text{lim}\right)\right)_{i,j-\half}
	\end{bmatrix},
	\end{equation}
	which uses the same velocity $\u^{(2)}_{\text{adv}}$ and density $\breve{\V \rho}^{(2)}_{\text{lim}}$ used to 
	update $\breve{\V \rho}^{n+1}$ in Eq~\eqref{eq_rk3}. This is the key step required to strongly couple the mass and momentum 
	convective operators.

}

%%%%%%%%%%%%%%%%%%%%%%%%%%%%%%%%%
%%%%%%%%%%%%%%%%%%%%%%%%%%%%%%%%%
\section*{Bibliography}
\begin{flushleft}
 \bibliography{WaterEntry}
\end{flushleft}

\end{document}